\documentclass[12pt,reqno,psfig,openbib]{amsart}

\usepackage[pagebackref=true, colorlinks=true, citecolor=blue]{hyperref}
\usepackage{amssymb,amsmath,graphicx,amsfonts,euscript}
\usepackage{fancyhdr}
\usepackage{epsfig}
\usepackage{setspace}
\usepackage{dsfont}
\usepackage[all]{xy}
\usepackage{amsmath}
\usepackage{enumitem}
\graphicspath{ {./figures/} }

\usepackage{caption}
\usepackage{subcaption}

\usepackage{soul}

\usepackage{float}

 \usepackage{pst-grad} % For gradients
 \usepackage{pst-plot} % For axes

\newtheorem{theorem}{Theorem}[section]

\numberwithin{equation}{section}
\numberwithin{figure}{section}

 \usepackage[left=1.25in,top=1.15in,right=1.25in]{geometry}
\usepackage{verbatim}
\newcommand{\spacer}{\vspace{1.5mm}}

\date{\today}

\renewcommand{\epsilon}{\varepsilon}
\newcommand{\ep}{\varepsilon}

\newcommand{\ov}{\overline}
\newcommand{\pa}{\partial}
\newcommand{\tht}{\theta}

\newcommand{\sig}{\sigma}
\newcommand{\kap}{\kappa}
\newcommand{\blds}{\boldsymbol}


\begin{document}

\title[Singular layer PINN methods for PPF]
{Singular layer Physics Informed Neural Network method for Plane Parallel Flows}

\author[T.-Y. Chang, 
G.-M. Gie, 
Y. Hong, and 
C.-Y.  Jung]
{
Teng-Yuan Chang$^{1}$, 
Gung-Min Gie$^{2}$, 
Youngjoon Hong$^{3}$, and   
Chang-Yeol  Jung$^4$}
\address{$^1$ Department of Applied Mathematics, National Yang Ming Chiao Tung University, 1001 Ta Hsueh Road, Hsinchu 300, Taiwan}
\address{$^2$ Department of Mathematics, University of Louisville, Louisville, KY 40292}
\address{$^3$ Department of Mathematical Sciences, KAIST, Korea}
\address{$^4$ Department of Mathematical Sciences, Ulsan National Institute of Science and Technology, 
Ulsan 44919, Korea}
\email{tony60113890048@gmail.com}
\email{gungmin.gie@louisville.edu}
\email{hongyj@kaist.ac.kr}
\email{cjung@unist.ac.kr}

\begin{abstract}
We  construct in this article the
semi-analytic 
Physics Informed Neural Networks (PINNs), called {\em singular layer PINNs} (or {\em sl-PINNs}), that are suitable
to predict the stiff solutions of plane-parallel flows at a small viscosity. 
Recalling the boundary layer analysis,
we first find the corrector for the problem which describes the singular behavior of the viscous flow inside boundary layers. 
Then, using the components of the corrector and its curl, 
we build our new {\em sl-PINN} predictions for the velocity and the vorticity 
by either  
embedding the explicit expression of the corrector (or its curl) in the structure
of PINNs  
or 
by 
training the implicit parts of the corrector (or its curl) together with the PINN predictions. 
Numerical experiments confirm that our new {\em sl-PINNs} produce stable
and accurate predicted solutions for the plane-parallel flows at a small viscosity.  
\end{abstract}

\maketitle

\tableofcontents

\section{Introduction}\label{S. Intro}
We aim in this article to construct 
a version of Physics Informed Neural Networks (PINNs)  
to predict the motion of certain symmetric flows 
%solutions of 
%the (stationary) Navier-Stokes equations (NSE) under a certain symmetry, 
especially when the viscosity is very small.  
To this end, 
we consider the stationary Navier-Stokes equations (NSE) in a periodic channel domain $\Omega := (0, L)^2_{\text{per}} \times (0, 1)$  
%with boundary at $z=0,1$, 
when a symmetry imposed to the velocity vector field $\blds{u}^{\ep}$ (sol. to NSE at a small viscosity $\ep>0$) is given in the form,
\begin{equation}\label{e:PP velocity_NSE}
\blds{u}^{\ep} = (u^{\ep}_1 (z), \, u^{\ep}_2 (x, z), \, 0 ).
\end{equation}
% stationary plane-parallel flows, at a small viscosity, in a periodic channel domain $\Omega := (0, L)^2_{\text{per}} \times (0, 1)$
% with boundary at $z=0,1$. 
% Under a  symmetry given in the form,
% \begin{equation}\label{e:PP velocity_NSE}
% \blds{u}^{\ep} = (u^{\ep}_1 (z), \, u^{\ep}_2 (x, z), \, 0 ):
% \end{equation}
Then the stationary NSE are reduced to:
\begin{equation}\label{e:NSE_PP}
\left\{\begin{array}{rl}
                             \spacer
                                        u^{\ep}_1
                                        - \ep {d^2_z u^{\ep}_1}%{\pa z^2}
                                        & \hspace{-2mm}
                                        = f_1,
                                        \quad
                                            \text{in } \Omega,\\
                             \spacer
                                        u^{\ep}_2
                                        - \ep {\pa^2_x u^{\ep}_2}%{\pa x^2}
                                        - \ep {\pa^2_z u^{\ep}_2}%{\pa z^2}
                                        + u^{\ep}_1 \, {\pa_x u^{\ep}_2}%{\pa x}
                                        & \hspace{-2mm}
                                        = f_2,
                                        \quad
                                            \text{in } \Omega,\\
                              \spacer
                                        u^{\ep}_2 \text{ is periodic}
                                        & \hspace{-2mm}
                                        \text{in $x$ with period $L$},\\
                              %\spacer
                                        u^{\ep}_i
                                        & \hspace{-2mm}
                                        = 0,
                                        \quad
                                        %\text{ }
                                        i=1,2, \text{ on } \Gamma, \text{ i.e., at }
                                        z = 0, 1,
        \end{array}\right.
\end{equation}
provided that the smooth data $\blds{f}$, which is periodic in $x$,
satisfies the symmetry in (\ref{e:PP velocity_NSE}) as well. 
%\begin{equation}\label{e:PP data}
%            \blds{f} = (f_1 (z), \, f_2 (x, z), \, 0).
%\end{equation}
The fluid motion, associated with our model equations (\ref{e:NSE_PP}), occurs in the tangential directions in $x$ and $y$, i.e., 
$u^\ep_1$ and $u^\ep_2$, 
depending on the tangential variable $x$ and the normal variable $z$ in $\Omega$. 
Here and throughout this article,
we use the notations, $d_z = d/dz$, $\pa_x = \pa/\pa_x$, and $\pa_z = \pa/\pa_z$.
\smallskip

We now introduce the vorticity of the Navier-Stokes solutions,
\begin{equation}\label{e:vorticity}
	\blds{\omega}^{\ep}
		= ( \omega^{\ep}_1(x, z), \omega^\ep_2(z), \omega^\ep_3(x, z))
		:= \textnormal{curl }\blds{u}^{\ep}
		= \big(
			-\pa_z u^\ep_2, \,
			d_z u^\ep_1, \,
			\pa_x u^\ep_2
			\big).
\end{equation}
By computing 
$-\pa_z$(\ref{e:NSE_PP})$_2$, {$d_z$(\ref{e:NSE_PP})$_1$}, and $\pa_x$(\ref{e:NSE_PP})$_2$, 
we write the vorticity formulation of NSE (\ref{e:NSE_PP}) in the form, 
\begin{equation}\label{e:NSE_PP_Vor_eqn}
\left\{\begin{array}{rl}
			    \spacer
                                        \omega^{\ep}_1
                                        - \ep {\pa^2_x \omega^{\ep}_1}%{\pa x^2}
                                        - \ep {\pa^2_z \omega^{\ep}_1}%{\pa z^2}
                                        - \omega^{\ep}_2 \, \omega^{\ep}_3
                                        + u^{\ep}_1 \, {\pa_x \omega^{\ep}_1}%{\pa x}
                                        & \hspace{-2mm}
                                        = -\pa_z f_2,
                                        \quad
                                            \text{in } \Omega,\\
                             \spacer
                                        \omega^{\ep}_2
                                        - \ep {d^2_z \omega^{\ep}_2}%{\pa z^2}
                                        & \hspace{-2mm}
                                        = d_z f_1,
                                        \quad
                                            \text{in } \Omega,\\
                             \spacer
                                        \omega^{\ep}_3
                                        - \ep {\pa^2_x \omega^{\ep}_3}%{\pa x^2}
                                        - \ep {\pa^2_z \omega^{\ep}_3}%{\pa z^2}
                                        + u^{\ep}_1 \, {\pa_x \omega^{\ep}_3}%{\pa x}
                                        & \hspace{-2mm}
                                        = \pa_x f_2,
                                        \quad
                                            \text{in } \Omega,\\
%                              \spacer
                                        \omega^{\ep}_i \text{ is periodic}
                                        & \hspace{-2mm}
                                        \text{in $x$ with period $L$, $i=1,3$}.        \end{array}\right.
\end{equation}

{
To derive a set of proper boundary conditions for the vorticity, 
we apply the
Lighthill principle 
by restricting (\ref{e:NSE_PP})$_{1, 2}$ on $\Gamma$, and using (\ref{e:NSE_PP})$_4$; 
see, e.g., \cite{GKLMN} as well. 
Then we obtain the boundary conditions for $\blds{\omega}^{\ep}$
as  
}
\begin{equation}\label{e:NSE_PP_Vor_BC}
\left\{\begin{array}{rl}
			    \spacer
                                        \pa_z \omega^{\ep}_1
                                        & \hspace{-2mm}
                                        = \dfrac{1}{\ep} f_2,
                                        \quad
                                            \text{on } \Gamma,\\
                             \spacer
                                        d_z \omega^{\ep}_2
                                        & \hspace{-2mm}
                                        = -\dfrac{1}{\ep} f_1,
                                        \quad
                                            \text{on } \Gamma,\\
                             \spacer
                                        \omega^{\ep}_3
                                        & \hspace{-2mm}
                                        = 0,
                                        \quad
                                            \text{on } \Gamma.
\end{array}\right.
\end{equation}

Setting $\ep = 0$ in the NSE (\ref{e:NSE_PP}) and 
(\ref{e:NSE_PP_Vor_eqn}),
we find  that 
the inviscid limits $\blds{u}^0 = (u^0_1 (z), \, u^0_2 (x, z), \, 0 )$ of $\blds{u}^{\ep}$ 
and 
$\blds{\omega}^0 = \text{curl }\blds{u}^0$ of $\blds{\omega}^{\ep}$
satisfy the following equations: 
% \begin{equation}\label{e:PP velocity_EE}
% \blds{u}^0 = (u^0_1 (z), \, u^0_2 (x, z), \, 0 ),
% \end{equation}
% that satisfies the system:
\begin{equation}\label{e:EE_PP}
\left\{\begin{array}{rl}
                             \spacer
                                        u^0_1
                                        & \hspace{-2mm}
                                        = f_1,
                                        \quad
                                            \text{in } \Omega,\\
                             \spacer
                                        u^0_2 + u^0_1 \, {\pa_x u^0_2}%{\pa x}
                                        & \hspace{-2mm}
                                        = f_2,
                                        \quad
                                            \text{in } \Omega,\\
                              %\spacer
                                        u^0_2 \text{ is periodic}
                                        & \hspace{-2mm}
                                        \text{in $x$ with period $L$}.
        \end{array}\right.
\end{equation}
\begin{equation}\label{e:EE_PP_Vor}
\left\{\begin{array}{rl}
			    \spacer
                                        \omega^{0}_1
                                        - \omega^{0}_2 \, \omega^{0}_3
                                        + u^{0}_1 \, {\pa_x \omega^{0}_1}%{\pa x}
                                        & \hspace{-2mm}
                                        = -\pa_z f_2,
                                        \quad
                                            \text{in } \Omega,\\
                             \spacer
                                        \omega^{0}_2
                                        & \hspace{-2mm}
                                        = d_z f_1,
                                        \quad
                                            \text{in } \Omega,\\
                             \spacer
                                        \omega^{0}_3
                                        + u^{0}_1 \, {\pa_x \omega^{0}_3}%{\pa x}
                                        & \hspace{-2mm}
                                        = \pa_x f_2,
                                        \quad
                                            \text{in } \Omega,\\
%                              \spacer
                                        \omega^{0}_i \text{ is periodic}
                                        & \hspace{-2mm}
                                        \text{in $x$ with period $L$, $i=1,3$}.        \end{array}\right.
\end{equation}

Under a sufficient regularity assumption on the data $f$, 
the well-posedness of (\ref{e:NSE_PP}),  (\ref{e:NSE_PP_Vor_eqn}), 
(\ref{e:EE_PP}) or (\ref{e:EE_PP_Vor}), 
and 
the asymptotic behavior of solutions 
to (\ref{e:NSE_PP}) and (\ref{e:NSE_PP_Vor_eqn}) 
are fully investigated in \cite{GJL1} and 
\cite{GKLMN} (for the time dependent case).  
More precisely, 
following the methodology of boundary layer analysis from, e.g., \cite{Book16, Ho95, SK87}, 
it is verified in those earlier works that, as the viscosity $\ep$ vanishes,  
the Navier-Stokes solution $\blds{u}^\ep$ converges strongly in $L^2$ to the Euler solution $\blds{u}^0$ 
and that 
the corresponding Navier-Stokes vorticity $\blds{\omega}^\ep$ 
converges to the Euler vorticity $\blds{\omega}^0$ in a certain weak sense, up to a positive measure supported on the boundary. 
Here we briefly recall this result from \cite{GJL1}:
\begin{theorem}\label{t:PPF}
        Under a sufficient regularity assumption on the data  $\blds{f}$, e.g.,
        $\blds{f} \in H^6(\Omega)$, %for $k \geq 6$,
  %       the difference between the plane-parallel viscous solution and
  %       its asymptotic expansion vanishes as the viscosity parameter tends to zero in the sense that
  %       \begin{equation}\label{e:conv_corrected_PPF}
  %       \left\{
  %       \begin{array}{l}
  %       		\spacer
  %                       \| \ov{{u}}_1 \|_{L^2(\Omega)}
  %                       +
  %                       \ep^{\frac{1}{2}}
  %                       \| d_z \ov{{u}}_1 \|_{L^2 (\Omega)}
  %                       		\leq \kap \ep,\\
		% \spacer
  %                       \| \ov{{u}}_2 \|_{L^2(\Omega)}
  %                       +
  %                       \ep^{\frac{1}{2}}
  %                       \| \nabla \ov{{u}}_2 \|_{L^2 (\Omega)}
  %                       		\leq \kap \ep,\\
		% \spacer
		% 	{
		% 	\| \ov{\omega}_1 \|_{L^2(\Omega)}
  %                       +
  %                       \ep^{\frac{1}{2}}
  %                       \| \nabla \ov{\omega}_1 \|_{L^2 (\Omega)}
  %                       		\leq \kap \ep^{\frac{1}{2}},
		% }\\
		% \spacer
		% 	\| \ov{\omega}_2 \|_{L^2(\Omega)}
  %                       +
  %                       \ep^{\frac{1}{2}}
  %                       \| d_z \ov{\omega}_2 \|_{L^2 (\Omega)}
  %                       		\leq \kap \ep^{\frac{3}{4}},\\
		% 	{
		% 	\| \ov{\omega}_3 \|_{L^2(\Omega)}
  %                       +
  %                       \ep^{\frac{1}{2}}
  %                       \| \nabla \ov{\omega}_3 \|_{L^2 (\Omega)}
  %                       		\leq \kap \ep.
		% }
  %       \end{array}
  %       \right.
  %       \end{equation}
%Moreover, as the viscosity $\ep$ tends to zero, 
        the viscous plane-parallel solution $\blds{u}^{\ep}$  converges to the corresponding inviscid solution  $\blds{u}^0$  in the sense that
\begin{equation}\label{e:VVL_PPF}
        \| \blds{u}^{\ep} - \blds{u}^0 \|_{L^2(\Omega)}
                        \leq \kap \ep^{\frac{1}{4}}, 
        \quad
        \text{ as } \ep \rightarrow 0.
\end{equation}
In addition, as $\ep \rightarrow 0$, we have 
\begin{equation}\label{e:vorticity_accumulation_1}
	\lim_{\ep \rightarrow 0} \big( \blds{\omega}^{\ep}, \, \blds{\psi} \big)_{L^2(\Omega)}
		=
		\big( \blds{\omega}^{0}, \, \blds{\psi} \big)_{L^2(\Omega)}
		+
		\big( \blds{u}^{0}|_{\Gamma} \times \blds{n}, \, \blds{\psi} \big)_{L^2(\Gamma)},
	\quad
	\forall \blds{\psi} \in C(\ov{\Omega}),
\end{equation}
which expresses the fact that
\begin{equation}\label{e:vorticity_accumulation_2}
	\lim_{\ep \rightarrow 0} \blds{\omega}^{\ep}
		=
			\blds{\omega}^{0}
		+
			( \blds{u}^{0}|_{\Gamma} \times \blds{n}) \delta_{\Gamma},
	% \quad
	% \forall \blds{\varphi} \in C(\ov{\Omega}),
\end{equation}
in the sense of weak$^*$ convergence of bounded measures on $\ov{\Omega}$.
\end{theorem}

In proving Theorem \ref{t:PPF} in \cite{GJL1} and 
\cite{GKLMN}, 
the authors employed the method of matching asymptotic expansions and 
introduce an expansion of $\blds{u}^\ep$ (hence of $\blds{\omega}^\ep$) at a small $\ep$ in the form, 
\begin{equation*}\label{e:assympt exp_PPF_temp}
        \blds{u}^{\ep}
                \simeq
                    \blds{u}^{0} 
                    + 
                    \blds{\varphi},
        \qquad
        \blds{\omega}^{\ep}
                \simeq
                    \blds{\omega}^{0} 
                    + 
                   \text{curl } \blds{\varphi}, 
\end{equation*}
by constructing an artificial function  $\blds{\varphi}$, called corrector, which  describes the singular behavior of $\blds{u}^\ep - \blds{u}^0$ near the boundary $\Gamma$; see below Section \ref{S:BL} for more information about the asymptotic expansion above and the corrector $\blds{\varphi}$. 
In fact, performing involved analysis, 
the authors verify the validity at  a small viscosity of the asymptotic expansion in the sense that 
\begin{equation}\label{e:conv_corrected_PPF}
        \left\{
        \begin{array}{l}
        		\spacer
                        \| \ov{{u}}_1 \|_{L^2(\Omega)}
                        +
                        \ep^{\frac{1}{2}}
                        \| d_z \ov{{u}}_1 \|_{L^2 (\Omega)}
                        		\leq \kap \ep,\\
		\spacer
                        \| \ov{{u}}_2 \|_{L^2(\Omega)}
                        +
                        \ep^{\frac{1}{2}}
                        \| \nabla \ov{{u}}_2 \|_{L^2 (\Omega)}
                        		\leq \kap \ep,\\
		\spacer
			{
			\| \ov{\omega}_1 \|_{L^2(\Omega)}
                        +
                        \ep^{\frac{1}{2}}
                        \| \nabla \ov{\omega}_1 \|_{L^2 (\Omega)}
                        		\leq \kap \ep^{\frac{1}{2}},
		}\\
		\spacer
			\| \ov{\omega}_2 \|_{L^2(\Omega)}
                        +
                        \ep^{\frac{1}{2}}
                        \| d_z \ov{\omega}_2 \|_{L^2 (\Omega)}
                        		\leq \kap \ep^{\frac{3}{4}},\\
			{
			\| \ov{\omega}_3 \|_{L^2(\Omega)}
                        +
                        \ep^{\frac{1}{2}}
                        \| \nabla \ov{\omega}_3 \|_{L^2 (\Omega)}
                        		\leq \kap \ep, 
		}
        \end{array}
        \right.
\end{equation}
where the difference between the viscous solution $\blds{u}^\ep$ and the proposed expansion $\blds{u}^0 + \blds{\varphi}$ is defined by 
$$
\ov{\blds{u}}
        =
        \blds{u}^{\ep}
                -
                (    \blds{u}^{0} 
                    + 
                    \blds{\varphi}),
        \qquad
        \ov{\blds{\omega}}
        =
        \text{curl } \ov{\blds{u}}.
        % \blds{\omega}^{\ep}
        %         \simeq
        %             \blds{\omega}^{0} 
        %             + 
        %             \nabla \times \blds{\varphi}
$$
Then, thanks to the properties of the corrector, 
Theorem \ref{t:PPF} follows from (\ref{e:conv_corrected_PPF}).
\smallskip

It is well-understood that a singular perturbation  problem, such as our model (\ref{e:NSE_PP}) or (\ref{e:NSE_PP_Vor_eqn}) - (\ref{e:NSE_PP_Vor_BC}), 
generates a thin layer near the boundary, 
{{called {\em  boundary layer}}, 
where  
a sharp transition of the solution occurs;}    
see, e.g.,   
\cite{
Book16, 
Ho95, 
OM08, 
SK87}. 
Generally speaking, 
when we approximate a solution  to any 
singularly perturbed boundary value problem, 
% by using the classical numerical methods,  
%or the Neural Networks, e.g., Physics Informed Neural Networks (PINNs),   
a very large computational error is created near the boundary, 
because of the stiffness nature of the solution inside boundary layers.  
Hence,  
to achieve a sufficiently accurate approximate solution, especially near the boundary,  
certain additional treatments are proposed and successfully employed, 
e.g., 
mesh refinements inside boundary layers \cite{RST96, St05} 
and 
semi-analytic methods using the boundary layer corrector 
\cite{%CJL19, 
GJL1, 
%GJL2, 
GJL22, 
HK82, 
%HJL13, 
HJT14}.  
The
main idea of the semi-analytic method 
or 
the mesh refinement method is to 
provide more information about the solution 
or 
introduce more finer approximation 
of the 
solution inside the boundary layer. 
These 
enriched methods have proven to be highly efficient to approximate the solutions of singularly perturbed boundary value problems.

Neural networks, 
in particular  
the Physics Informed Neural Networks (PINNs), have developed as an important and effective tool 
for predicting solutions to partial differential equations.  
One special feature of PINNs is that 
the loss function is defined based on 
a residual from the differential equation under consideration where the collocation points in
the time-space domain are used as the input for the neural network.  
This property makes PINNs  well-suited particularly for solving time-dependent, multi-dimensional equations in a domain involved in complex geometries, see, e.g., 
\cite{KZK19, 
LMMK21, 
LDKRKS20, 
RPK19, 
SDBSK20, 
YMK21, KDJH21, 
MZK20, 
KKLPWY21,  LLF98, SS18, 
XZRW21, 
WTP21}.   
Although the PINNs have become a popular tool
in scientific machine learning, 
%enabling the integration of physics-based and data-driven modeling within a deep learning framework,  However, despite their success, the robustness
there is an well-known issue on their robustness 
for certain types of problems, including singularly perturbed boundary value problems. 
% remains an ongoing issue. For instance, PINNs have
% limited accuracy when dealing with complex and highly nonlinear flow patterns, such as
% turbulence, vortical structures, and boundary layers. Addressing this challenge is a priority
% in scientific machine learning research, as developing robust and reliable models is essential
% for advancing the field. 
As investigated in some recent works, e.g., \cite{ACD, Tony9, FSF22,  FT20,  pinn_lim, PINN_bl01, PINN_bl02, PINN_bl03, Tony3, LMMK21}, 
the  
PINNs experience some limitations to capture a  singular behavior of solutions, e.g., inside a boundary/interior layer.  
Hence it is of great interest to 
overcome this challenging task by 
developing a methodology to build a version of robust and reliable PINN methods for singular perturbations.
 
In a series of recent papers \cite{GHJL, GHJM, GHJ}, 
the authors suggested a new semi-analytic approach to enrich the conventional PINN methods applied to boundary layer problems. 
In fact, using the so-called corrector (which describes the singular behavior of the solution inside boundary/interior layers), the authors construct a version of PINNs enriched by the corrector and this new semi-analytic PINNs are proven to outperform the conventional PINNs for the 1 or 2D singularly perturbed elliptic/parabolic equations tested in the papers.

Following the methodology in \cite{GHJL, GHJM, GHJ}, 
we aim in this article 
to construct the semi-analytic PINNs, 
called {\em singular layer PINNs} (or {\em sl-PINNs}), 
that are suitable to predict the stiff solutions of 
plane-parallel flows at a small viscosity. 
To this end, in Section \ref{S:BL} below,   
we first recall from \cite{GKLMN, Book16, GJL1} the asymptotic expansions  for 
(\ref{e:NSE_PP}) and (\ref{e:NSE_PP_Vor_eqn}) - (\ref{e:NSE_PP_Vor_BC}).  
Then, using the corresponding boundary layer analysis,   
we construct our new {\em singular layer PINNs} ({\em sl-PINNs}) for the velocity vector field in Sections \ref{comp_u1} - \ref{comp_u2}, 
and ({\em sl-PINNs}) for the vorticity in Section \ref{comp_vor}. 
The numerical computations for our {\em sl-PINNs} with a comparison with those from the conventional PINNs appear in Section \ref{comparion_velocity} for the velocity and 
\ref{comparion_vor} for the vorticity. 
It is verified numerically that our novel {\em sl-PINNs} 
capture well the sharp transition 
and hence produce good predicted solutions for the stiff plane-parallel flows.

\section{Boundary layer analysis}\label{S:BL}
 
We recall briefly from \cite{GJL1}  
the boundary layer analysis for our model systems (\ref{e:NSE_PP}) and (\ref{e:NSE_PP_Vor_BC}): 

We propose an asymptotic expansion of $\blds{u}^\ep$, solution of (\ref{e:NSE_PP}), as
\begin{equation}\label{e:assympt exp_PPF}
        \blds{u}^{\ep}
                \simeq
                    \blds{u}^{0} 
                    + 
                    \blds{\varphi},
\end{equation}
where $\blds{\varphi}$ is the corrector, which will be  constructed below, in the form,
\begin{equation}\label{e:Tht_PPF}
        \blds{\varphi}
                = \big(\varphi_1 (z), \, \varphi_2 (x, z), \, 0 \big).
\end{equation}

Inserting $\blds{\varphi}$ into the difference between the equations (\ref{e:NSE_PP}) and (\ref{e:EE_PP}),
we collect all the terms of dominant order $\ep^0$ and write the asymptotic equation for $\blds{\varphi}$ as the weakly coupled system below:  
\begin{equation}\label{e:Eq_Tht}
        \left\{
                \begin{array}{rl}
                        \spacer \displaystyle
                                \varphi_1 - \ep {d^2_z \varphi_1}%{\pa z^2}
                        &
                                \hspace{-2mm}
                                = 0,
                                \quad
                                \text{in } \Omega ,\\
                        \spacer \displaystyle
                                \varphi_2
                                - \ep \Delta  \varphi_2
                                %- \ep {\pa^2_z \tht_2}%{\pa z^2}
                                +
                                    %\textcolor{blue}{
                                    \big(
                                        \varphi_1
                                        + u^0_1%\big|_{\Gamma}
                                    \big)
                                    {\pa_x \varphi_2}%{\pa x}
                                    %}
                                    +
                                    \varphi_1 {\pa_x u^0_2}
                        & \hspace{-2mm}
                                =
                                0, 
                                % - \tht_1 {\pa_x u^0_2}, %\big|_{\Gamma},
                        \quad
                                \text{in } \Omega,\\
                        %\spacer \displaystyle
                                \varphi_i
                        & \hspace{-2mm}
                        = - u^0_i,
                        \quad
                                \text{on } \Gamma, \text{ } i=1,2.
                \end{array}
        \right.
\end{equation}

To construct the corrector $\blds{\varphi}$ as an (approximate) solution  of (\ref{e:Eq_Tht}) above,
we introduce a $C^\infty$ truncation function $\sigma_L$ (and $\sigma_R$) near the boundary at $z = 0$ (and $z = 1$) such that
\begin{equation}\label{e:Ch7:sig_LR}
\begin{array}{ccc}
                \sig_L(z) =
                        \left\{
                        \begin{array}{l}
                                \spacer
                                        1, \text{ } 0 \leq z \leq 1/4,\\
                                        0, \text{ } z  \geq 1/2,
                        \end{array}
                        \right.
        & \quad &
                \sig_R(z) = \sig_L(1-z).
\end{array}
\end{equation}
Using the truncations in (\ref{e:Ch7:sig_LR}),
we define the first component $\varphi_1$ of the corrector as
\begin{equation}\label{e:Ch7:Tht_1}
        \varphi_1(z)
                = \sig_L(z) \, \varphi_{1, \, L} (z )
                    + \sig_R(z) \, \varphi_{1, \, R} (z ),
\end{equation}
where
\begin{equation}\label{e:Ch7:Tht_1_L_R}
%\left\{
%\begin{array}{l}
	%\spacer
        \varphi_{1, \, L} (z )
                =
                        - u^0_1(0) \, e^{ - z / \sqrt{\ep}}, %\\
    \qquad \qquad
        \varphi_{1, \, R} (z )
                =
                        - u^0_1(1) \, e^{ - (1-z) / \sqrt{\ep}}.
%\end{array}\right.
\end{equation}

One can verify that ${\varphi}_1$ enjoys the estimates,
\begin{equation}\label{e:Lp_est_Tht_1}
\| {d_z^{m} \varphi_1} \|_{L^p(\Omega)}
                                            \leq
                                                    \kap \ep^{\frac{1}{2p} - \frac{m}{2}},
                                            \quad m \geq 0,
                                            \quad 1 \leq p \leq \infty,
\end{equation}
as well as the proposed asymptotic equation in (\ref{e:Eq_Tht}),
\begin{equation}\label{e:Ch7:Eq_Tht_1}
        \left\{
                \begin{array}{rl}
                        \spacer \displaystyle
                                {\varphi_1}%{\pa t}
                                - \ep {d^2_z \varphi_1}%{\pa z^2}
                        &
                                =
                                    e.s.t.,
                                    %-\ep \Big\{
                                    %            2 \sig_L^{\pri} \, \dfrac{\pa \tht_{1, \, L}}{\pa z}
                                    %            + \sig_L^{\pri \pri} \, \tht_{1, \, L}
                                    %            - 2 \sig_R^{\pri} \, \dfrac{\pa \tht_{1, \, R}}{\pa z}
                                    %            + \sig_R^{\pri \pri} \, \tht_{1, \, R}
                                    %        \Big\},
                                \quad
                                \text{in } \Omega,\\
                        %\spacer
                        \displaystyle
                                \varphi_1
                        &
                                = - u^0_1,
                                \quad
                                \text{on } \Gamma,
                \end{array}
        \right.
\end{equation}
where
$e.s.t.$ denotes a term that is exponentially small with respect to the small viscosity parameter $\ep$
in the spaces, e.g., $C^s(\Omega)$ or $H^s(\Omega)$, $0 \leq s \leq \infty$.

Using the first component $\varphi_1$,
we define $\varphi_{2, \, L}$ and $\varphi_{2, \, R}$ as the solutions of
\begin{equation}\label{e:Ch7:Eq_Tht_2_L}
        \left\{
                \begin{array}{rl}
                        \spacer \displaystyle
                                { \varphi_{2, \, L}}%{\pa t}
                                - \ep \Delta \varphi_{2, \, L}
                                + \big( \varphi_{1, \, L} + u^0_1 \big)
                                {\pa_x \varphi_{2, \, L}}
                                &
                                =
                                - \varphi_{1, \, L} \,
                                {\pa_x u^0_2} ,
                        %& \displaystyle
                        		\quad
                                z > 0,\\
                        \spacer
			\varphi_{2, \, L} \text{ is}
			&  \text{$L$-periodic in } x,\\
                        \spacer \displaystyle
                                \varphi_{2, \, L}
                          & = - u^0_2,
                          \quad
                                z = 0, \\
                        %\spacer
			\varphi_{2, \, L}
			&
			\rightarrow 0,
			\quad
			\text{as } z \rightarrow \infty,			
	\end{array}
        \right.
\end{equation}
and
\begin{equation}\label{e:Ch7:Eq_Tht_2_R}
        \left\{
                \begin{array}{rl}
                        \spacer \displaystyle
                                {\varphi_{2, \, R}}%{\pa t}
                                -\ep \Delta \varphi_{2, \, R}
                                + \big( \varphi_{1, \, R} + u^0_1  \big)
                                {\pa_x \varphi_{2, \, R}}%{\pa x}
                         &
                         =
                                - \varphi_{1, \, R} \, {\pa_x u^0_2} ,
                        \quad
			z < 1,\\
                        \spacer
			\varphi_{2, \, R} \text{ is}
			&
			\text{$L$-periodic in } x,\\
                        \spacer \displaystyle
                                \varphi_{2, \, R}
                          & = - u^0_2,
                          %\text{ }
                          \quad
                          z = 1, \\
                        %\spacer
			\varphi_{2, \, R}
			&
			\rightarrow 0
			\quad
			\text{as } z \rightarrow -\infty.
                \end{array}
        \right.
\end{equation}

The elliptic equation (\ref{e:Ch7:Eq_Tht_2_L}) (or (\ref{e:Ch7:Eq_Tht_2_R})) is well-defined as
all the coefficients are of class $C^{\infty}(0, \infty)$ (or $C^{\infty}(- \infty, 1)$).
Moreover, by performing the energy estimates,
one can verify that the ${ \varphi_{2, \, L}}$ (or ${ \varphi_{2, \, R}}$)
behaves like an exponentially decaying function with respect to the stretched variable $z/\sqrt{\ep}$ (or $(1-z)/\sqrt{\ep}$) in the sense that, for $* = L, R$,
\begin{equation}\label{e:Lp_est_Tht_2_LR}
\left\{
\begin{array}{l}
	\spacer
	\| \pa_x^k \varphi_{2, \, *} \|_{L^p(\Omega)}
                                            \leq
                                                    \kap \ep^{\frac{1}{2p}},
                                            \quad k \geq 0,
                                            \quad 1 \leq p \leq \infty,	\\
        \spacer
	\| \pa_x^k \pa_z \varphi_{2, \, *} \|_{L^1(\Omega)}
                                            \leq
                                                    \kap,
                                            \quad k \geq 0, \\
	\spacer
	\| \pa_x^k \pa_z \varphi_{2, \, *} \|_{L^2(\Omega)}
                                            \leq
                                                    \kap \ep^{- \frac{1}{4}},
                                            \quad k \geq 0, \\
	%\spacer
	\| \pa_x^k \pa^2_z \varphi_{2, \, *} \|_{L^2(\Omega)}
                                            \leq
                                                    \kap \ep^{- \frac{3}{4}},
                                            \quad k \geq 0;
\end{array}
\right.
\end{equation}
see, e.g., \cite{GKLMN, Book16} 
for the detailed verification of the estimates above.  
%here, see \cite{GKLMN} where
%the (more involved) time-dependent problem is fully analyzed.

Now, using the truncations in (\ref{e:Ch7:sig_LR}),
we define the second component $\varphi_2$ of the corrector as
\begin{equation}\label{e:Ch7:Tht_2}
        \varphi_2(x, z)
                = \sig_L(z) \, \varphi_{2, \, L} (x, z)
                    + \sig_R(z) \, \varphi_{2, \, R} (x, z),
\end{equation}
so that
\begin{equation}\label{e:Lp_est_Tht_2}
	\text{
	$\varphi_2$ satisfies the estimates in (\ref{e:Lp_est_Tht_2_LR})
	with $\varphi_{2, \, *}$ replaced by $\varphi_2$.
	}
\end{equation}
In addition, one can verify that
the second component $\varphi_2$
satisfies the proposed asymptotic equation in(\ref{e:Eq_Tht}) up to an exponentially small error: % $E(\tht_2)$:
\begin{equation}\label{e:Ch7:Eq_Tht_2}
        \left\{
                \begin{array}{rl}
                        \spacer \displaystyle
                                {\varphi_2}%{\pa t}
                                -\ep \Delta \varphi_2
                                + (\varphi_1 + u^0_1) {\pa_x \varphi_2}%{\pa x}
                                                                %\\
                        %\spacer \qquad

				+ \varphi_1 {\pa_x u^0_2}%{\pa x}
				%E(\tht_2)
    &
    =
				e.s.t.,
				%-
				%\ep \Big\{
                                  %              2 \sig_L^{\pri} \, \dfrac{\pa \tht_{2, \, L}}{\pa z}
                                  %              + \sig_L^{\pri \pri} \, \tht_{2, \, L}
                                  %              - 2 \sig_R^{\pri} \, \dfrac{\pa \tht_{2, \, R}}{\pa z}
                                  %              + \sig_R^{\pri \pri} \, \tht_{2, \, R}
                                  %          \Big\},
                        %& \displaystyle
                        \quad
                                \text{in } \Omega,\\
                        %\spacer \displaystyle
                                \varphi_2
                                &
                                = - u^0_2,
                        %&
                        \quad
                                \text{on } \Gamma.
                \end{array}
        \right.
\end{equation}

When we construct our {\em sl-PINN} approximations below, 
we will use the profile of the corrector $\blds{\varphi}$ constructed in this section.

% {\red \bf
% ADD properties of the corrector... 
% SMALL and it is approx. of identity...
% }
% Hence, 
% for our convenience in the sections below, we introduce the normalized corrector $\ov{\blds{\varphi}}$ by setting:
% \begin{equation}\label{e:cor_normalized}
%     \ov{\varphi}_{i, k} 
%         =
%            \dfrac{1}{{\varphi}_{i, k}|_{z = 0}} {\varphi}_{i, k}, 
%     \qquad
%     i = 1, 2, 
%     \,\,
%     k = L, R. 
% \end{equation}

\section{Numerical experiments for the velocity}\label{sec_numerics}

In this section, 
we present  
the framework of our new {\em singular layer PINN} ({\em sl-PINN}) methods for solving the equations (\ref{e:NSE_PP}) with a small viscosity $0<\ep \ll 1$. 
Then we compare the computational results of our {\em sl-PINNs} and those obtained by the conventional PINNs (PINNs).
For this purpose,  
we first recall the conventional PINN (PINN) structure below: 

An $L$-layer feed-forward neural network (NN) is recursively defined by 
\begin{align*}
    \text{input layer: }& \mathcal{N}^{0}(\blds{x}) = \blds{x} \in \mathbb{R}^{N_{0}}, \\
    \text{hidden layers: }& \mathcal{N}^{l}(\blds{x}) = \sigma(\blds{W}^{l} \mathcal{N}^{l-1}(\blds{x}) + \blds{b}^{l}) \in \mathbb{R}^{N_{l}}, \quad 1\leq l \leq L-1, \\
    \text{output layer: }& \mathcal{N}^{L}(\blds{x}) = \blds{W}^{L} \mathcal{N}^{L-1}(\blds{x}) + \blds{b}^{L} \in \mathbb{R}^{N_{L}}.
\end{align*}
Here 
$\blds{x}$ is an input vector, 
and 
$\blds{W}^{l} \in \mathbb{R}^{N_{l} \times N_{l-1}}$ and $\blds{b}^{l} \in \mathbb{R}^{N_{l}}$, $l=1,\cdots,L$, are the (unknown) weights and bias respectively. 
The $N^{l}$, $l=0,\cdots,L$, is the number of neurons in the $l$-th layer,  
and $\sigma$ is an activation function.
 
When we employ the conventional PINNs to solve (\ref{e:NSE_PP}), 
we construct first a Neural Network (NN) to  predict $u_{1}^{\ep}$, and 
then introduce another NN to predict $u_{2}^{\ep}$ with using the  predicted solution of $u_{1}^{\ep}$   
because (\ref{e:NSE_PP}) can be solved sequentially first for $u_{1}^{\ep}$  and then for $u_{2}^{\ep}$. 
Let $\hat{u}_{1}=\hat{u}_{1}(z; \blds{\theta}_{1})$ 
and 
$\hat{u}_{2}=\hat{u}_{2}(x, z; \blds{\theta}_{2})$ 
denote the outputs   
of the first and second NNs 
where $\blds{\theta}_{1}, \blds{\theta}_{2}$ are the collections of the corresponding weights and bias for the outputs $\hat{u}_{1}$ and $\hat{u}_{2}$ respectively.   
Recalling  that 
the conventional PINNs rely on minimizing the loss function, which is composed of equation residuals and boundary conditions, 
we define the loss functions for the 
conventional PINN predictions $\hat{u}_1$ and $\hat{u}_2$  
as 
\begin{align} \label{loss1PINN}
    \mathcal{L}_{1}(\blds{\theta}_{1}; \mathcal{T}_{1}) = \frac{1}{|\mathcal{T}_{1}|} \sum_{z\in\mathcal{T}_{1}} \left| \hat{u}_{1} - \ep d_z^2 \hat{u}_{1} - f_1 \right|^{2} + \left| \hat{u}_{1}(z=0) \right|^{2} + \left| \hat{u}_{1}(z=1) \right|^{2},
\end{align}
where $\mathcal{T}_{1} \subset [0, 1]$ is the set of training points, and
\begin{equation} \label{loss2PINN}
\begin{split}
    \mathcal{L}_{2} (\blds{\tht}_{2}; \mathcal{T}_{2}, \mathcal{T}_{\Gamma}, \mathcal{T}_{B}) &= \frac{1}{|\mathcal{T}_{2}|} \sum_{(x, z)\in \mathcal{T}_{2}} \left| \hat{u}_{2} - \ep \Delta \hat{u}_{2} + \hat{u}_{1} \partial_x \hat{u}_{2} - f_2 \right|^2 \\
    &+ \frac{1}{|\mathcal{T}_{\Gamma}|} \sum_{(x, z)\in \mathcal{T}_{\Gamma}} \left| \hat{u}_{2}(x, z=0) \right|^{2} + \left| \hat{u}_{2}(x, z=1) \right|^{2} \\
    &+ \frac{1}{|\mathcal{T}_{B}|} \sum_{(x, z)\in \mathcal{T}_{B}} \left| \hat{u}_{2}(x=0, z) - \hat{u}_{2}(x=L, z) \right|^2,
\end{split}
\end{equation}
where $\mathcal{T}_{2} \subset \Omega$ is the training data set chosen in the domain, $\mathcal{T}_{\Gamma} \subset \Gamma$ is the training data set chosen at $z=0, 1$ and $\mathcal{T}_{B} \subset \partial\Omega \setminus \Gamma$ is the training data set chosen on the periodic boundary. 

As mentioned in the introduction above, 
excessive computational costs are needed for PINNs to predict the solution which changes rapidly in a region such as boundary layer,   
and 
in many cases, the conventional PINNs happen to  
produce incorrect predictions for this stiff problem.  
% The numerical experiments of PINNs based on (\ref{loss1PINN}) and (\ref{loss2PINN}) will be shown in Section \ref{comparion}. 
In order to train a network to learn the solution's behavior accurately inside the boundary layers,  
we will enrich our training solutions by leveraging the corrector, which is similar to our early studies \cite{GHJL, GHJM, GHJ}. The detail construction of our new training solutions are in the following Sections \ref{comp_u1} and \ref{comp_u2}.

\subsection{{\em sl-PINN} predictions  for $u_1^\ep$}\label{comp_u1}

To obtain an accurate prediction of $u_1^\ep$ at a small $\ep$,  
we enrich the NN output by incorporating 
a proper approximation of the corrector. 
More precisely, 
following the boundary layer analysis studied in Section \ref{S:BL}, 
we modify the traditional PINN structure and 
construct 
the training solution in the form, 
\begin{equation} \label{tilde u1}
    \tilde{u}_{1}(z; \blds{\tht}_{1}) = \hat{u}_{1}(z; \blds{\tht}_{1}) + \tilde{\varphi}_{1,L}(z; \blds{\tht}_{1}) + \tilde{\varphi}_{1,R}(z; \blds{\tht}_{1}). 
\end{equation}
Here 
$\hat{u}_{1}$ is the neural network output 
and, 
using the explicit profile of the corrector   ${\varphi}_1$ in (\ref{e:Ch7:Tht_1_L_R}),  
we set 
\begin{align}
    \tilde{\varphi}_{1,L}(z; \blds{\tht}_{1}) &= -\hat{u}_{1}(0; \blds{\tht}_{1}) e^{-z/\sqrt{\ep}},
    \label{e:cor_PINN_1L}
    \\
    \tilde{\varphi}_{1,R}(z; \blds{\tht}_{1}) &= - \hat{u}_{1}(1; \blds{\tht}_{1}) e^{-(1-z)/\sqrt{\ep}}.
    \label{e:cor_PINN_1R}
\end{align}
%are the training version of the correctors.
% Because of the exponentially decaying functions $e^{-z/\sqrt{\ep}}$ and $e^{-(1-z)/\sqrt{\ep}}$,
Note that the predicted solution $\tilde{u}_{1}$ in (\ref{tilde u1}) satisfies  the zero boundary condition at $z=0, 1$, up to negligible $e.s.t.$ with respect to $\ep$, 
thanks to the fact that the function 
$e^{-z/\sqrt{\ep}}$ (or $e^{-(1-z)/\sqrt{\ep}}$) 
decays exponentially fast away from the boundary at $z = 0$ (or $z=1$).   
Hence  the loss function for this {\em sl-PINN} for $u^\ep_1$ is defined by
\begin{equation} \label{loss1}
    \mathcal{L}_{1}(\blds{\tht}_{1}; \mathcal{T}_{1}) = \frac{1}{|\mathcal{T}_{1}|} \sum_{z\in \mathcal{T}_{1}} \left| \tilde{u}_{1} - \ep d_z^2 \tilde{u}_{1} - f_1 \right|^2,
\end{equation}
where $\mathcal{T}_{1}$ is the set of training points chosen in $\{z \ | \  0\leq z \leq 1 \}$, 
i.e., the boundary condition for $u^\ep_1$ is already enforced to be satisfied (up to $e.s.t.$) in the structure of $\tilde{u}_{1}$. 

\begin{figure}[H]
    \centering
    
    \begin{subfigure}{\textwidth}
        \centering
        \includegraphics[height=6cm]{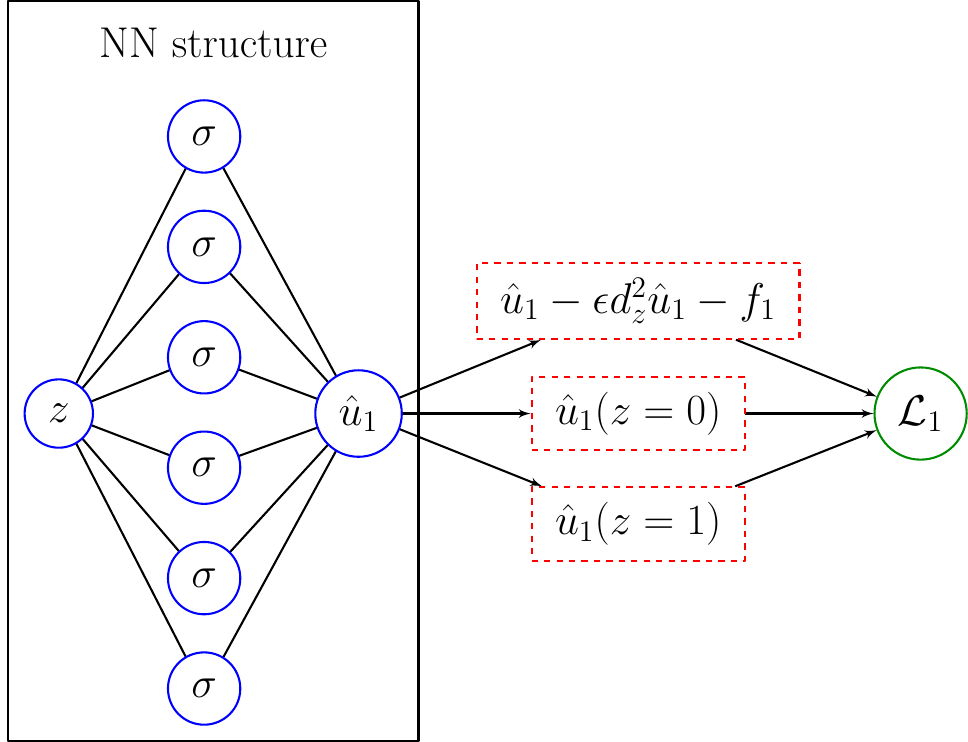}
        \caption{traditional PINN for $u^\ep_1$}
    \end{subfigure}
    
    \begin{subfigure}{\textwidth}
        \centering
        \includegraphics[height=6cm]{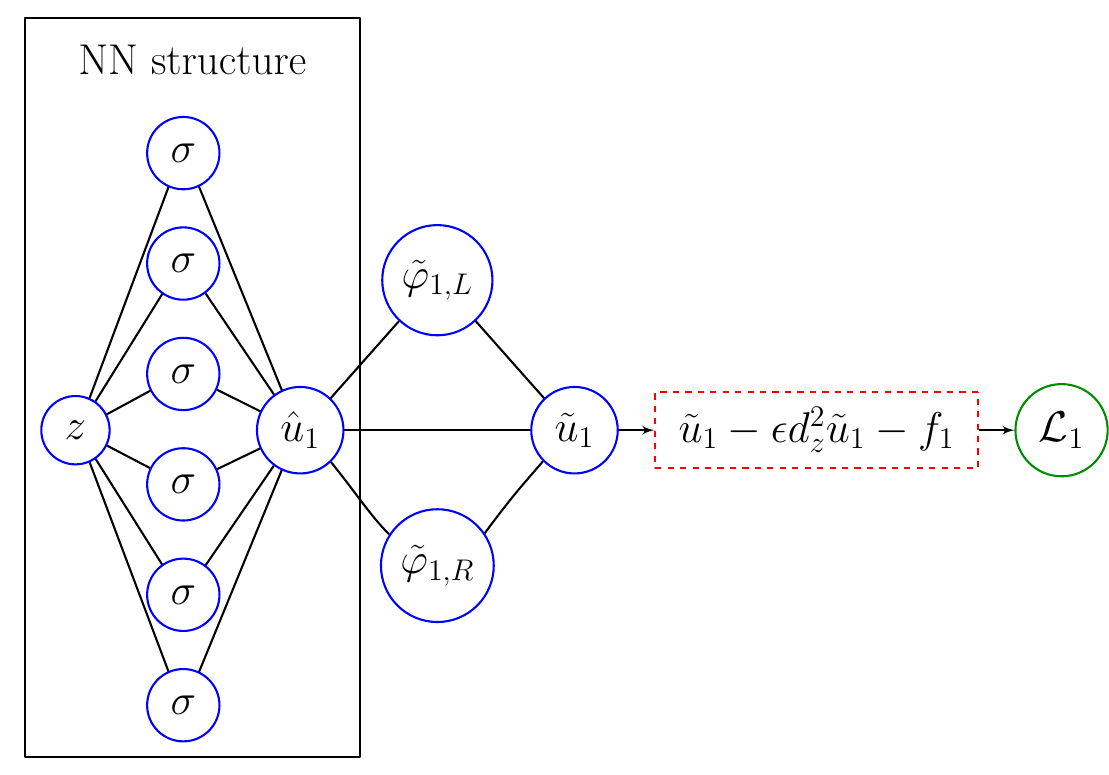}
        \caption{{\em sl-PINN} for $u^\ep_1$}
    \end{subfigure}
    
    \caption{Schematic difference between PINN and {\em sl-PINN} for $u_1^\ep$. }
    \label{schematic_u1}
\end{figure}

Thanks to (\ref{e:Ch7:Eq_Tht_1}), 
we notice from the direct computation with (\ref{tilde u1}) that 
\begin{align*}
    \tilde{u}_{1} - \ep d_z^2 \tilde{u}_{1} - f_1 = -\ep (z-1) d_z^2 \hat{u}_{1}(z; \blds{\tht}_{1}) -  \hat{u}_{1}(z; \blds{\tht}_{1}) - f_1 
    + 
    e.s.t.
\end{align*}
Thus the loss function (\ref{loss1}) remains bounded as $\ep\to 0$. 
This important fact makes the training of our {\em sl-PINN} for $u^\ep_1$ reliable and 
it leads to a solid minimization and accurate prediction. 
The schematics are illustrated in Figure \ref{schematic_u1} for both PINN  and {\em sl-PINN}. 
The numerical experiments appear in Section \ref{comparion_velocity}.

% is not considered in the loss of EPINNs since it is automatically satisfied by way of the construction of $\tilde{u}_{1}$ (\ref{tilde u1}). In addition, it is worth noting that the loss function (\ref{loss1}) remains bounded as $\ep\to 0$. That is, 
% \begin{align*}
%     \tilde{u}_{1} - \ep d_z^2 \tilde{u}_{1} - f_1 = -\ep (z-1) d_z^2 \hat{u}_{1}(z; \blds{\tht}_{1}) -  \hat{u}_{1}(z; \blds{\tht}_{1}) - f_1,
% \end{align*}
% this implies $\mathcal{L}_{1}$ (\ref{loss1}) is bounded as $\ep\to 0$. This makes training reliable and leads to a solid minimization and accurate approximation. The schematics are illustrated in Figure \ref{schematic_u1} for both PINNs and EPINNs. The numerical experiments can be found in Section \ref{comparion_velocity}.

\subsection{{\em sl-PINN} predictions for $u_2^\ep$}\label{comp_u2}

Here we construct our {\em sl-PINN} for the second component $u^\ep_2$, solution of (\ref{e:NSE_PP})$_2$. 
 
\begin{figure}[H]
    \centering
    
    \begin{subfigure}{\textwidth}
        \centering
        \includegraphics[height=6cm]{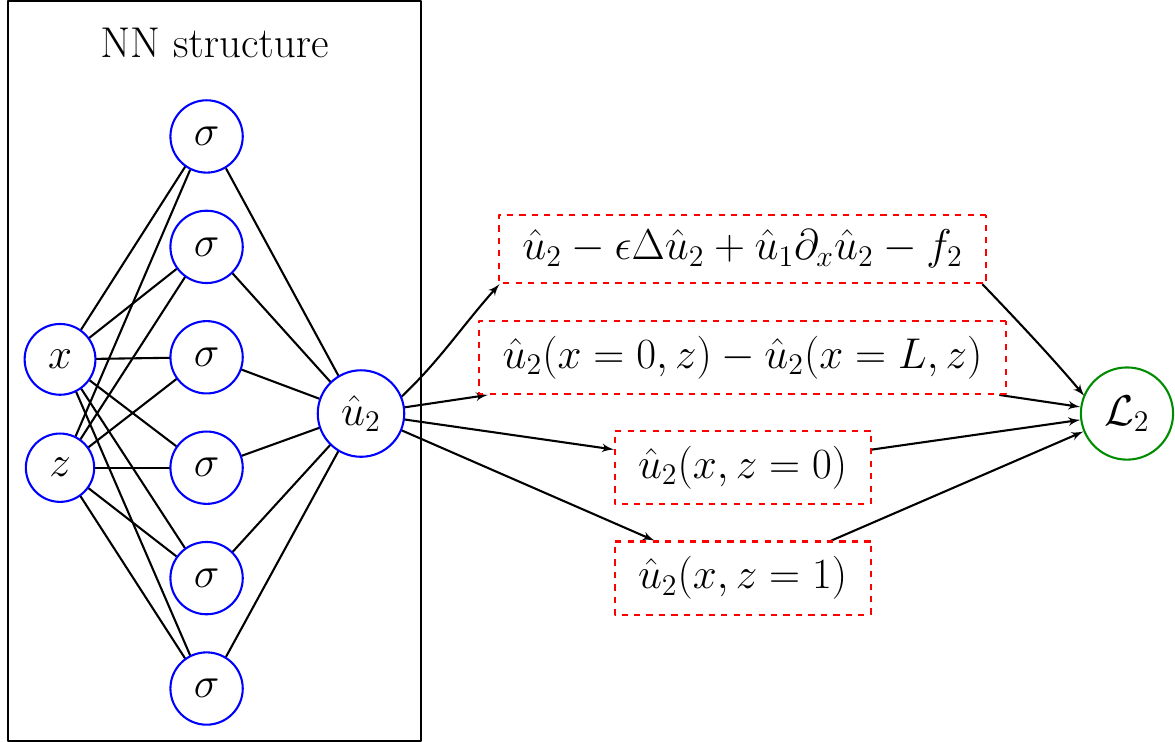}
        \caption{traditional PINN for $u^\ep_2$}
    \end{subfigure}
    
    \begin{subfigure}{\textwidth}
        \centering
        \includegraphics[height=6cm]{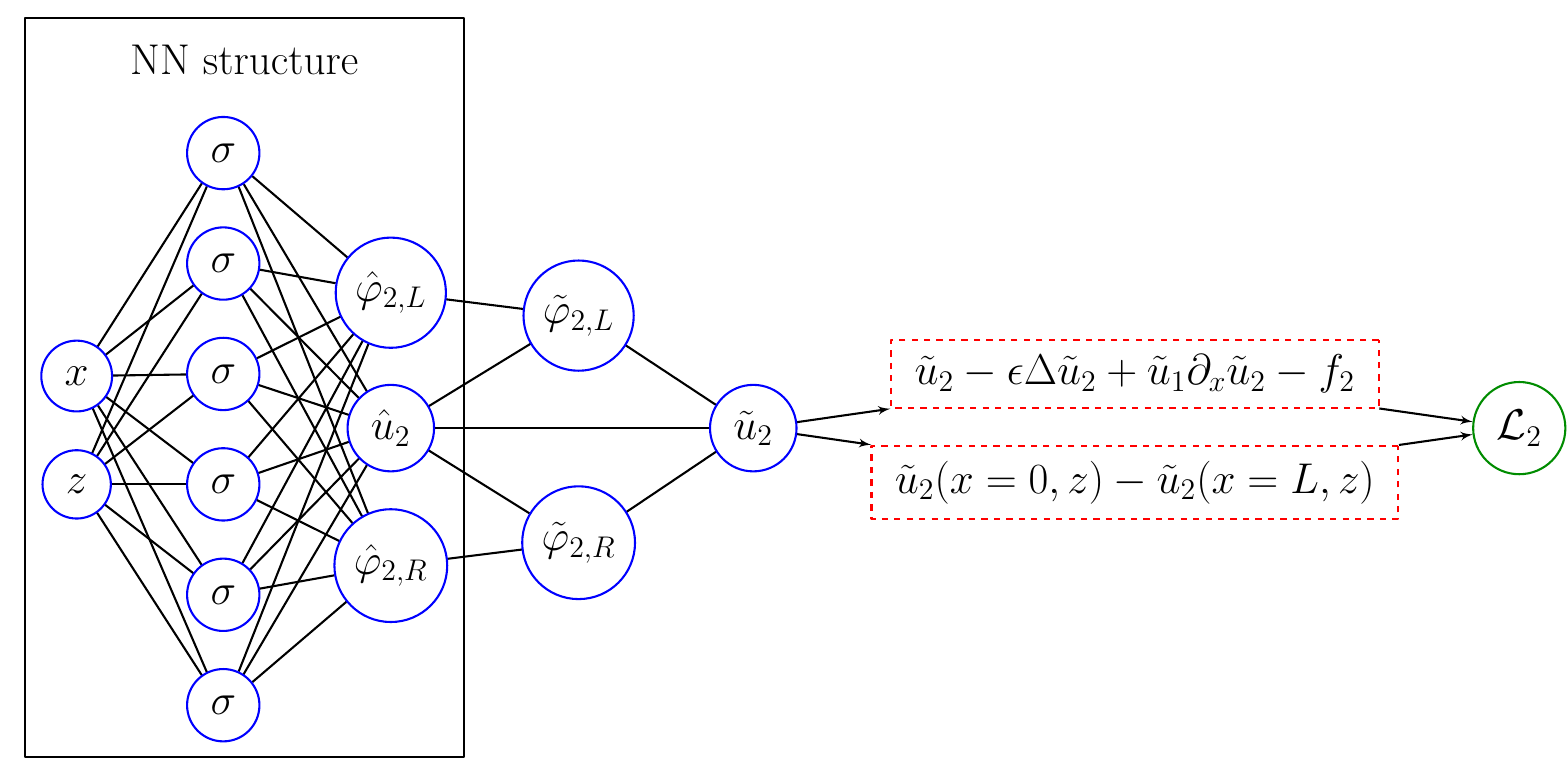}
        \caption{{\em sl-PINN} for $u^\ep_2$}
    \end{subfigure}
    
    \caption{Schematic difference between PINN and {\em sl-PINN} for  $u_2^\ep$.}
    \label{schematic_u2}
\end{figure}

Because the explicit expression of the corrector $\varphi_2$, defined in (\ref{e:Ch7:Tht_2}), 
is not available, we are no longer able to 
construct our (singular layer) PINN prediction of $u^\ep_2$ as simple as that for $u^\ep_1$ in (\ref{tilde u1}).  
Instead, 
we train the viscous solution $u^\ep_2$ and the corrector $\varphi_2$ at the same time. That is, 
we 
construct one neural network, 
in part to train $u^\ep_2$ and also to train $\varphi_2$ as well,  
and hence they interact each other in the form,
% To approximate $u_2^\ep$ by enrich method, we first construct the training solutions for corrector $\tht_{2}$. Since the explicit form of the corrector of the second component, $\tht_{2, L}$ and $\tht_{2, R}$ defined in (\ref{e:Ch7:Eq_Tht_2_L}) and (\ref{e:Ch7:Eq_Tht_2_R}), is not well-known, a properly designed trainable solution is considered. 
\begin{equation} \label{tilde u2}
    \tilde{u}_{2}(x, z; \blds{\Theta}_{2}) = \hat{u}_{2}(x, z; \blds{\theta}_{2}) 
    + 
    \tilde{\varphi}_{2, L}(x, z; \blds{\theta}_{2,s,L}) 
    + 
    \tilde{\varphi}_{2, R}(x, z; \blds{\theta}_{2,s,R}), 
\end{equation}
and  
\begin{align}
    \spacer
    \tilde{\varphi}_{2, L}(x, z; \blds{\theta}_{2}, \blds{\theta}_{2,s,L}) 
        &= 
        -\hat{u}_{2}(x, 0; \blds{\theta}_{2}) e^{-z/\sqrt{\ep}} 
        + 
        z e^{-z/\sqrt{\ep}} \hat{\varphi}_{2, L}(x, z; \blds{\theta}_{2,s,L}), \label{tilde tht2L}\\
    \tilde{\varphi}_{2, R}(x, z; \blds{\theta}_{2}, \blds{\theta}_{2,s,R}) 
        &= 
        -\hat{u}_{2}(x, 1; \blds{\theta}_{2}) e^{-(1-z)/\sqrt{\ep}} 
        + 
        (1-z) e^{-(1-z)/\sqrt{\ep}} \hat{\varphi}_{2, R}(x, z; \blds{\theta}_{2,s,R}), \label{tilde tht2R}
\end{align}
where $\hat{u}_{2}$, $\hat{\varphi}_{2, L}$ and $\hat{\varphi}_{2, R}$ are the network outputs 
and 
$\blds{\Theta}_{2} = \{ \blds{\theta}_{2}, \blds{\theta}_{2,s,L}, \blds{\theta}_{2,s,R}  \}$ with 
$\blds{\theta}_{2}$, 
$\blds{\theta}_{2,s,L}$, and $\blds{\theta}_{2,s,R}$ are 
the unknown parameters of the outputs $\hat{u}_{2}$, $\hat{\varphi}_{2, L}$, and $\hat{\varphi}_{2, R}$ respectively.
% The output layer of the second NN architecture is three-dimensional where additional two neurons are dedicated to representing $\hat{\tht}_{2, L}$ and $\hat{\tht}_{2, R}$. 
We refer readers to see the picture of NN structure in Figure \ref{schematic_u2} with a comparison to the conventional PINN structure.

As analyzed in (\ref{e:Lp_est_Tht_2_LR}); see \cite{GKLMN, GJL2} as well, 
it is verified that $\varphi_{2, L}$ and $\varphi_{2, R}$ 
behave like exponentially decaying functions 
from the boundaries at $z = 0$ and $z = 1$ (at least in the $H^k$ sense for $k = 0, 1, 2$). 
Using this property, 
we construct our {\em sl-PINN} predicted solutions for $\varphi_{2, L}$ and $\varphi_{2, R}$ 
in part by using the exponentially decaying functions with respect to 
the stretched variables $z/\sqrt{\ep}$ and $(1-z)/\sqrt{\ep}$;  
the first terms on the right-hand side of (\ref{tilde tht2L}) and (\ref{tilde tht2R}). 
In addition, in order to implement the (implicit) interaction between the correctors and the equation, 
we add the training parts in the structure, the second terms on 
on the right-hand side, in (\ref{tilde tht2L}) and (\ref{tilde tht2R}). 

By our construction, the 
{\em sl-PINN} predictions of the correctors satisfies 
\begin{itemize}
    \item 
        $\tilde{\varphi}_{2,L} = 
                -\hat{u}_{2}$ at $z=0$  
        \text{ and }
        $\tilde{\varphi}_{2,R} = 
                -\hat{u}_{2}$ at $z=1$,   \vspace{2mm}
    \item   
        $\tilde{\varphi}_{2,L} \to 0$ as $z \to \infty$ 
        and 
        $\tilde{\varphi}_{2,R} \to 0$ as $z \to -\infty$,
\end{itemize}
and hence we notice  that 
\begin{itemize}
    \item
        $\tilde{u}_{2} = 0$, up to an $e.s.t.$, at $z = 0, 1$.
\end{itemize}

Now, because the zero boundary condition is 
ensured in the structure of $\tilde{u}_2$ (up to an $e.s.t$),  
we define the loss function as 
\begin{equation} \label{loss2}
\begin{split}
    \mathcal{L}_{2} (\blds{\Theta}_{2}; \mathcal{T}_{2}, \mathcal{T}_{B}) &= \frac{1}{|\mathcal{T}_{2}|} \sum_{(x, z)\in \mathcal{T}_{2}} \left| \tilde{u}_{2} - \ep \Delta \tilde{u}_{2} + \tilde{u}_{1} \partial_x \tilde{u}_{2} - f_2 \right|^2 \\
    &+ \frac{1}{|\mathcal{T}_{B}|} \sum_{(x, z)\in \mathcal{T}_{B}} \left| \tilde{u}_{2}(x=0, z) - \tilde{u}_{2}(x=L, z) \right|^2,
\end{split}
\end{equation}
where $\mathcal{T}_{2}$ is the training set chosen in the domain $\Omega$ and $\mathcal{T}_{B}$ is the training set chosen on the periodic boundary $\partial \Omega \setminus \Gamma$. Since there is no boundary layer occurring in $x$ direction, the periodic boundary condition is simply enforced by penalizing the equivalent in the second term of the loss function. 
It is straightforward 
(by direct computations with (\ref{tilde u2}) - (\ref{tilde tht2R})) 
to verify that the loss function $\mathcal{L}_{2}$ defined in (\ref{loss2}) is bounded as $\ep \to 0$. Therefore, the training process is reliable enough to reach an expected minimizer. The numerical experiments can be found in the next subsection \ref{comparion_velocity}.

\subsection{Comparison between the conventional PINNs and the {\em sl-PINNs} for the velocity $\blds{u}^\ep$}\label{comparion_velocity}

To examine the performance of the singular layer PINNs ({\em sl-PINNs}) enriched by the correctors,  
we introduce an exact solution $\blds{u}^{\ep}$ 
for the problem (\ref{e:NSE_PP}) 
where the periodicity is set to be $L=1$ in the $x$-direction.  
By choosing the data $f_{1}(z) = 1+z(1-z)$, we find the exact solution $u_{1}^{\ep}$  as 
\begin{align} \label{exactu1}
    u_{1}^{\ep}(z) = (1-2\ep) \left( 1- \frac{1-e^{-\frac{1}{\sqrt{\ep}}}}{1-e^{-\frac{2}{\sqrt{\ep}}}} \left( e^{-\frac{z}{\sqrt{\ep}}} + e^{-\frac{1-z}{\sqrt{\ep}}} \right)\right).
\end{align}
By sequentially choosing the data $f_{2}$ as
\begin{align*}
    f_{2}(x, z) = f_{1}(z) (1+\sin(2\pi x)) + 4\pi^2 \ep u_{1}^{\ep}(z) \sin(2\pi x) + 2\pi (u_{1}^{\ep}(z))^{2} \cos(2\pi x),
\end{align*}
we find that $u_{2}^{\ep}$ as 
\begin{align} \label{exactu2}
    u_{2}^{\ep}(x, z) = u_{1}^{\ep}(z) (1+\sin(2\pi x)).
\end{align}

\begin{figure}
    \begin{subfigure}{0.5\textwidth}
        \centering
        \includegraphics[width=\linewidth]{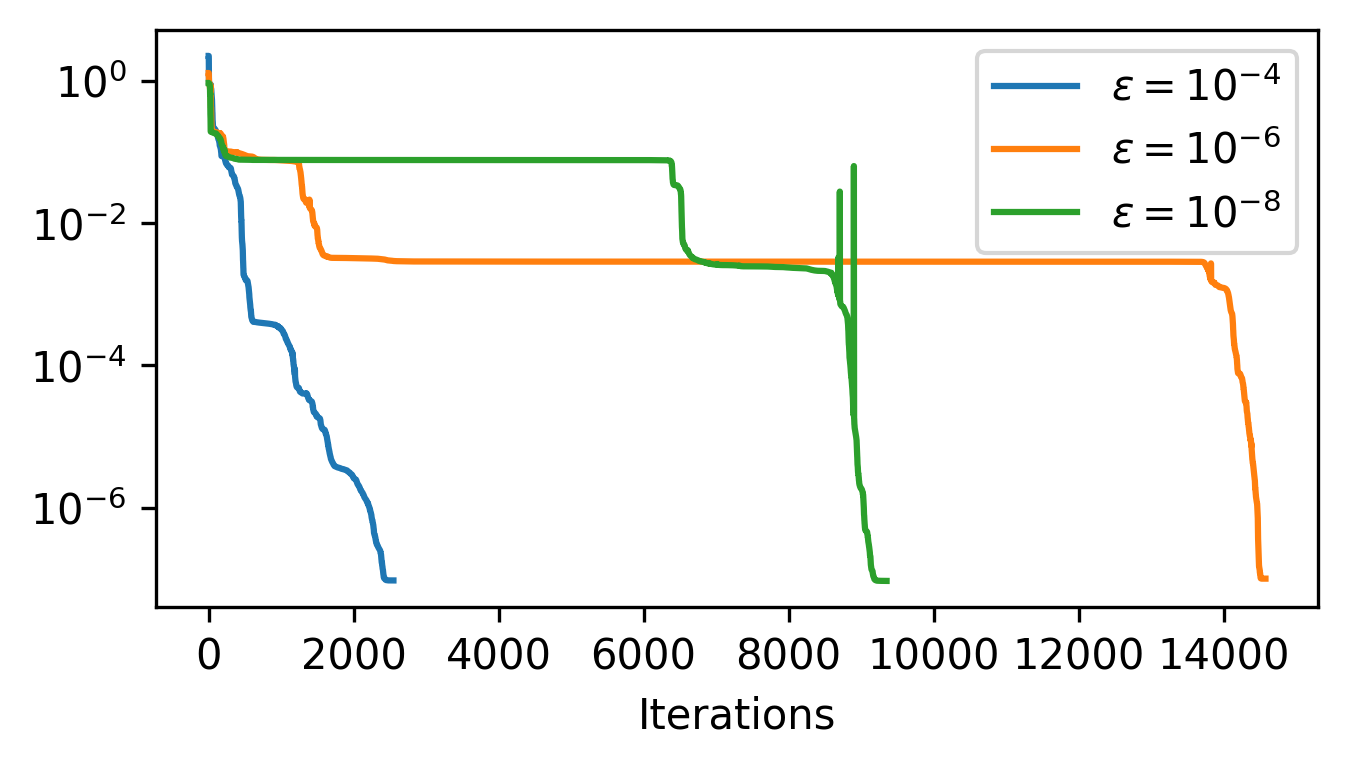}
        \caption{PINNs}
    \end{subfigure}%
    \begin{subfigure}{0.5\textwidth}
        \centering
        \includegraphics[width=\linewidth]{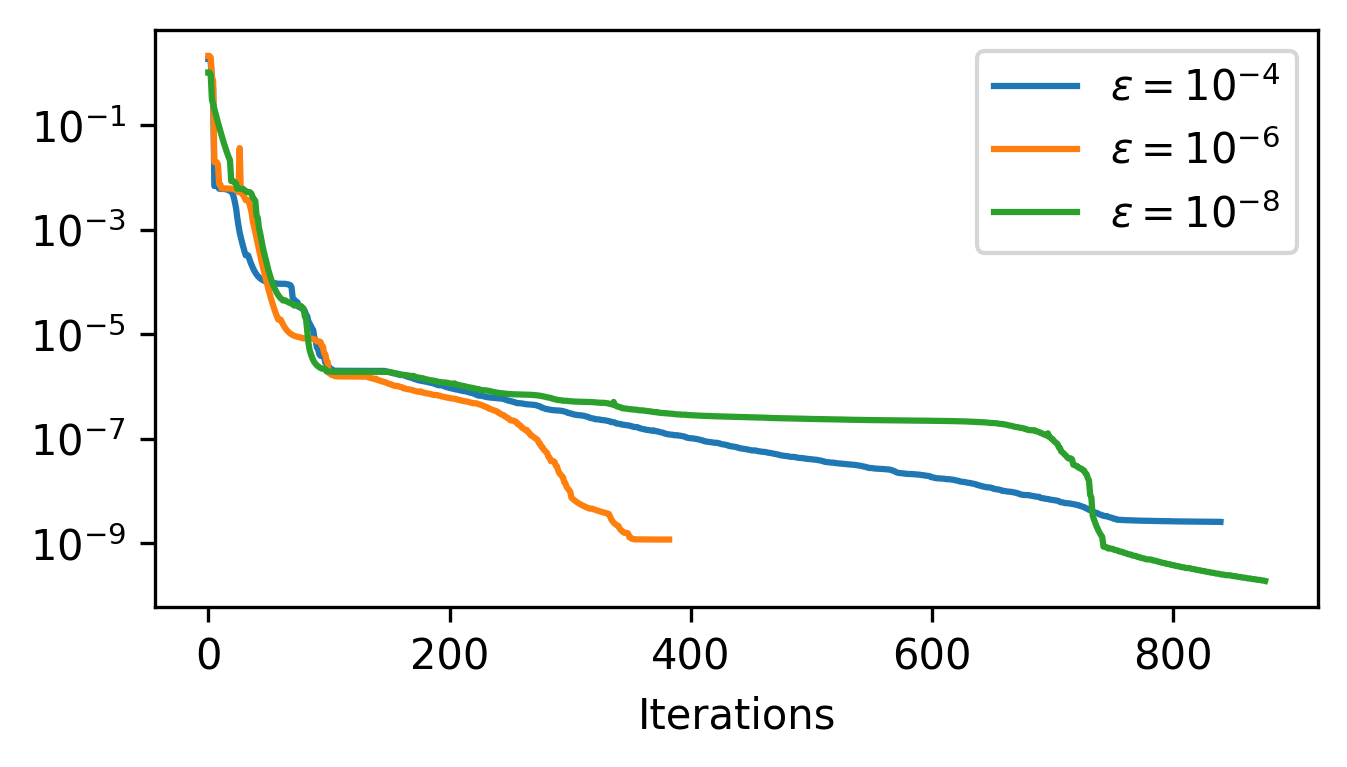}
        \caption{{\em sl-PINNs}}
    \end{subfigure}
    \caption{The loss values during the training process for $u_{1}^{\ep}$.}
    \label{plotloss1}
\end{figure}

\begin{figure}
    \begin{subfigure}{0.5\textwidth}
        \centering
        \includegraphics[width=\linewidth]{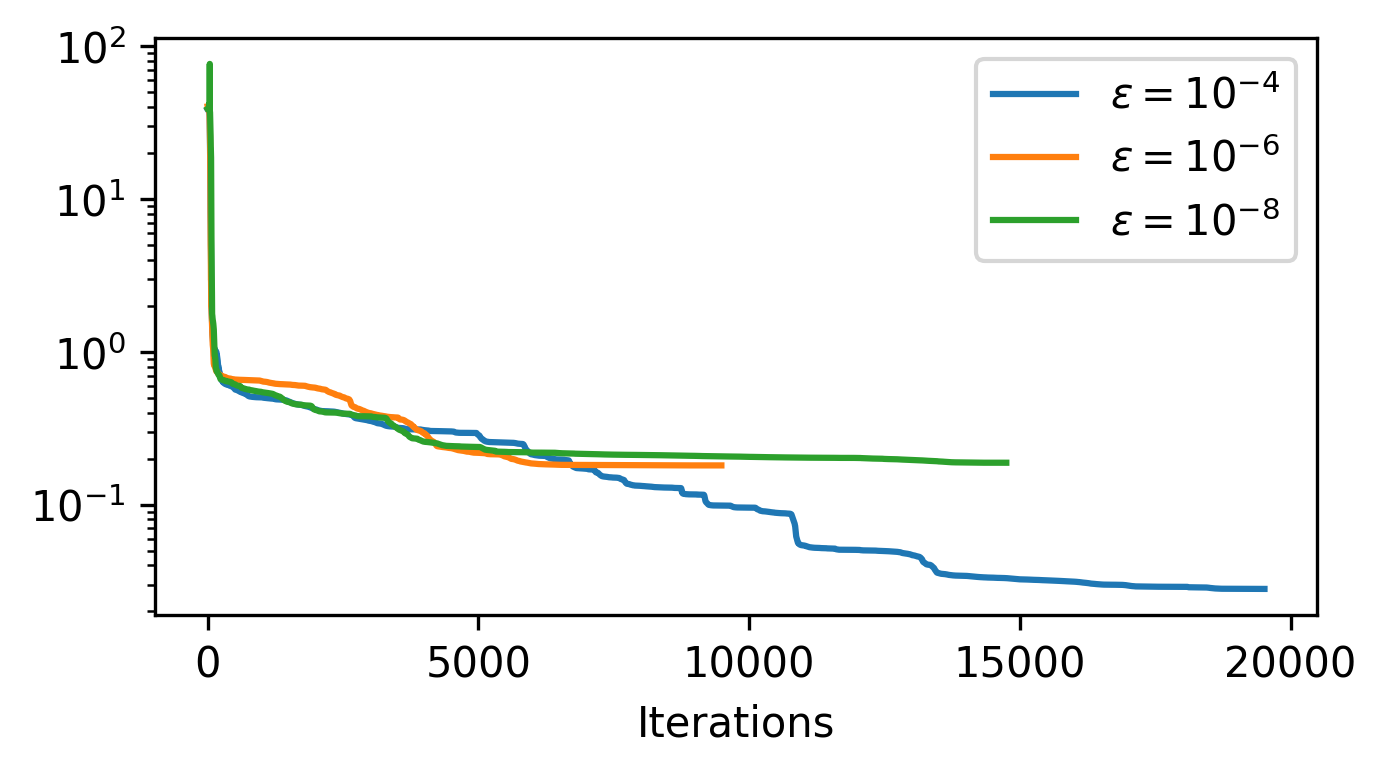}
        \caption{PINNs}
    \end{subfigure}%
    \begin{subfigure}{0.5\textwidth}
        \centering
        \includegraphics[width=\linewidth]{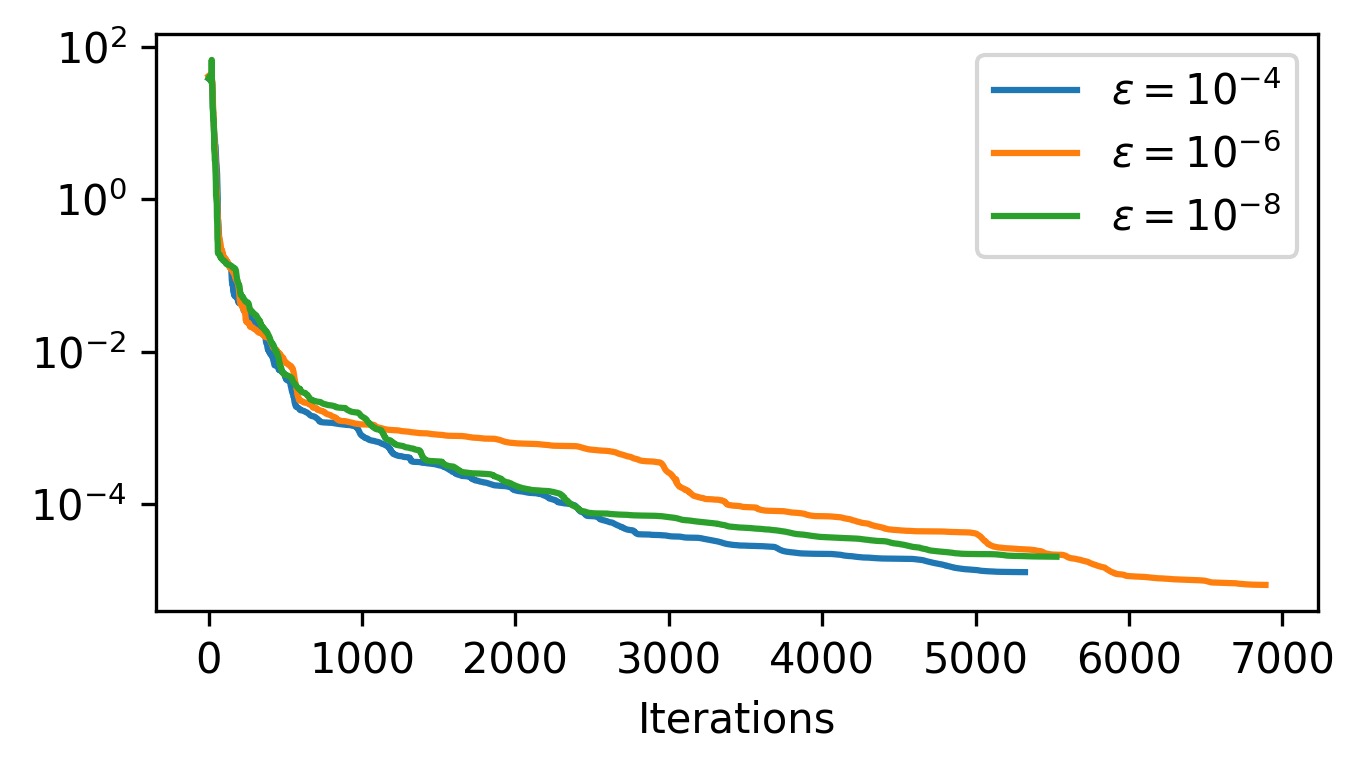}
        \caption{{\em sl-PINNs}}
    \end{subfigure}
    \caption{The loss values during the training process for $u_{2}^{\ep}$}
    \label{plotloss2}
\end{figure}

We aim to verify the high-performance 
of our {\em sl-PINNs} for a small (but fixed) viscosity $0 < \ep \ll 1$ without excessive computational cost. Specifically, our goal is to sustain the efficiency of our learning method by preserving low computational costs while maintaining accurate approximations. 
To this end, we use two single-hidden-layer Neural Networks (NNs) with 20 neurons each to approximate $u_{1}^{\ep}$ and $u_{2}^{\ep}$ respectively. 
In addition, we select the number of training points $N=25$ uniformly distributed inside the domain $\Omega_{z}=[0, 1]$ for training the loss $\mathcal{L}_{1}$ (\ref{loss1}) and $N=50^2$ uniformly distributed inside the domain $\Omega=[0, 1] \times [0, 1]$ for training the loss $\mathcal{L}_{2}$ (\ref{loss2}). The L-BFGS optimizer with learning rate $0.1$ is adopted to find a local minimum of $\mathcal{L}_{1}$ and $\mathcal{L}_{2}$ respectively. The maximum iteration is set as $50,000$, but the automatic termination is applied if tolerance is less than $10^{-8}$. In our experiments, we train the cases for $\epsilon=10^{-3}, 10^{-4}, 10^{-5}, 10^{-6}, 10^{-7}$, and $10^{-8}$ with the same aforementioned settings. Five independent runs are conducted for each case to present the most accurate result. For the comparison, the conventional PINNs (PINNs) is also performed under the same settings for each case. We present here the loss values during the training process in Figures \ref{plotloss1} and \ref{plotloss2} for $\ep=10^{-4}$, $10^{-6}$, and $10^{-8}$. 

After completing the training process, we investigate the accuracy by plotting the testing points uniformly on the three locations in $z$ direction: $N_{test}=50$ points on the left boundary layer $\{ z: 0\leq z \leq \sqrt{\ep} \}$, $N_{test}=50$ points on the right boundary layer $\{ z: 1-\sqrt{\ep} \leq z \leq 1 \}$ and $N_{test}=500$ points outside the boundary layer $\{ z: \sqrt{\ep} < z < 1-\sqrt{\ep} \}$. In the $x$-direction, $N_{test}=500$ testing points are simply drawn uniformly in $\{ x: 0\leq x\leq 1 \}$. The results are shown in Tables \ref{tableu1} and \ref{tableu2} for approximating $u_{1}^{\ep}$ and $u_{2}^{\ep}$ respectively. 

\begin{figure}
    \centering

    \begin{subfigure}[t]{0.3\textwidth}
        \raggedleft
        \includegraphics[width=\linewidth]{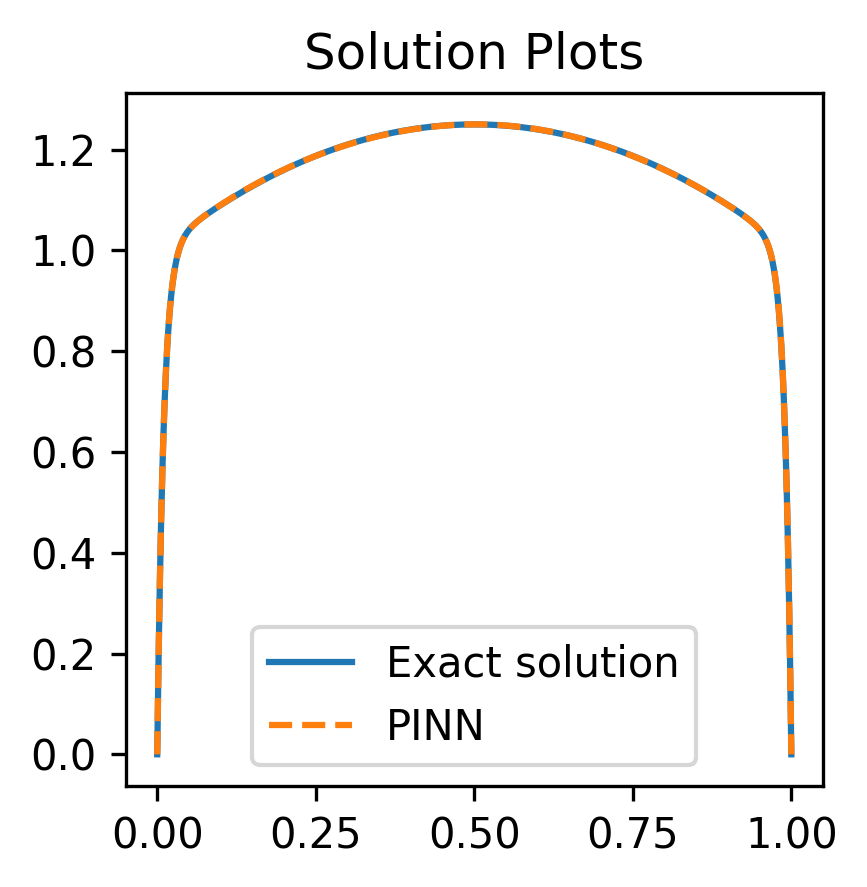}
        \includegraphics[width=\linewidth]{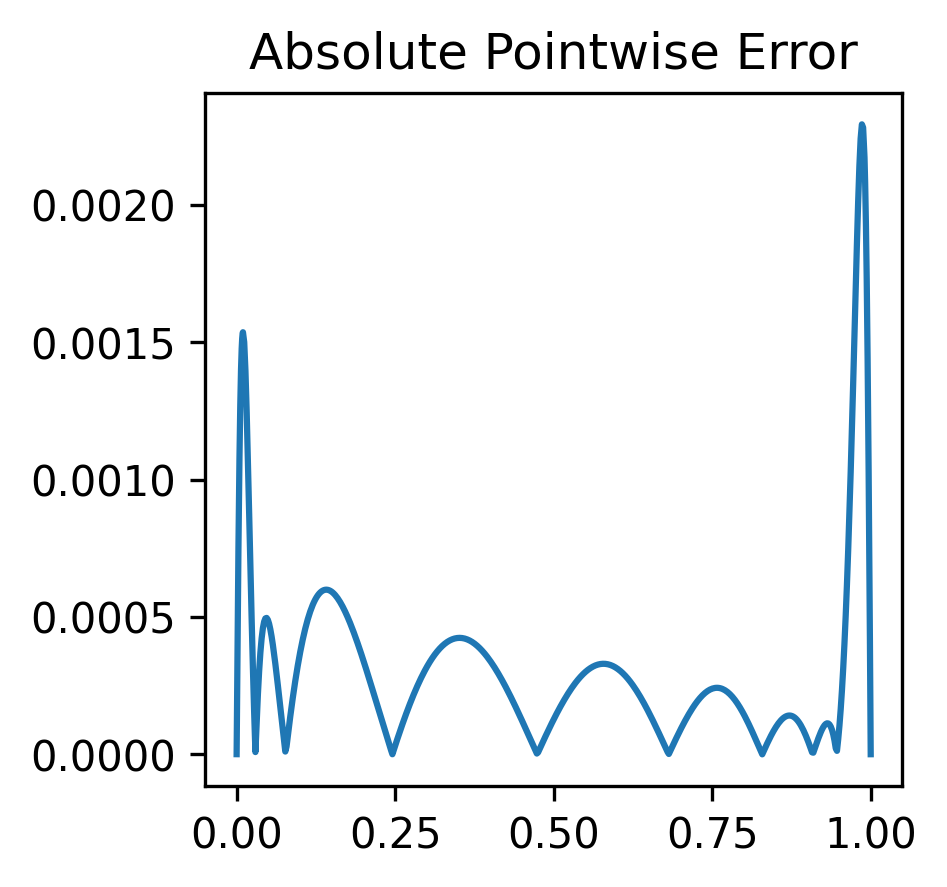}
        \caption{$\ep=10^{-4}$}
    \end{subfigure}%
    \hfill
    \begin{subfigure}[t]{0.3\textwidth}
        \raggedleft
        \includegraphics[width=\linewidth]{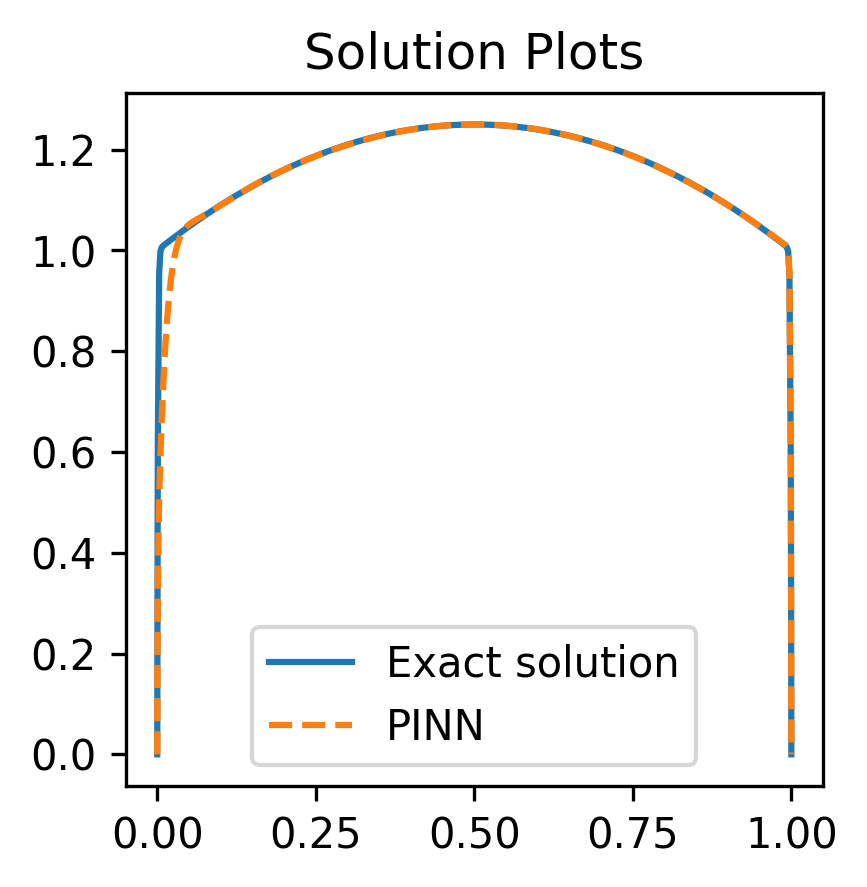}
        \includegraphics[width=\linewidth]{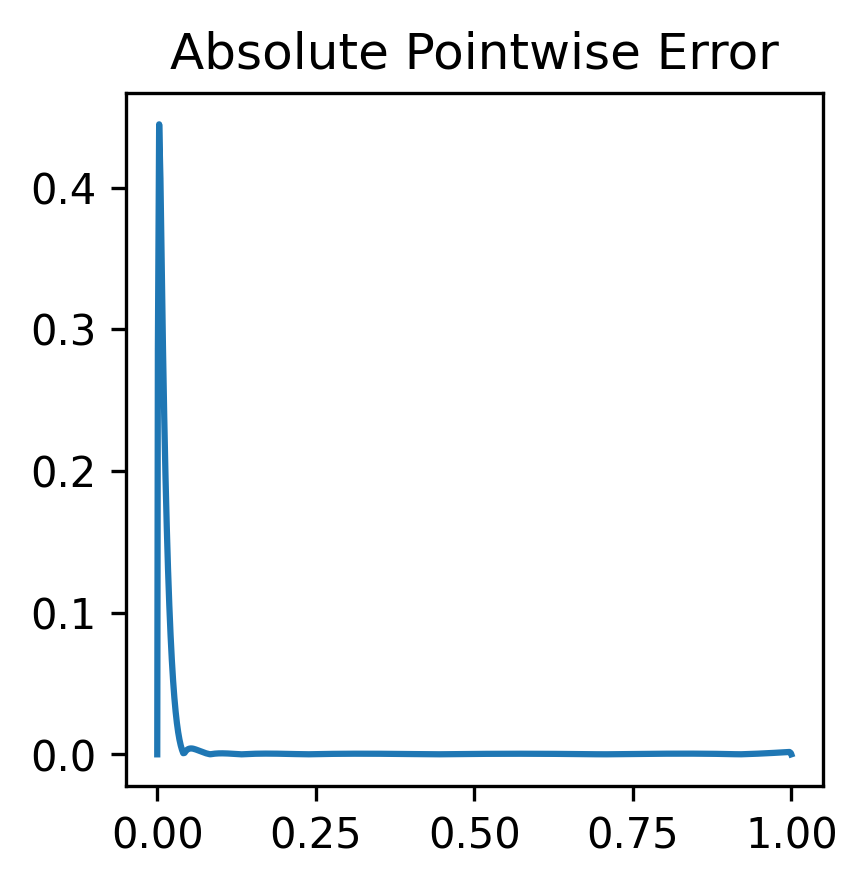}
        \caption{$\ep=10^{-6}$}
    \end{subfigure}%
    \hfill
    \begin{subfigure}[t]{0.3\textwidth}
        \raggedleft
        \includegraphics[width=\linewidth]{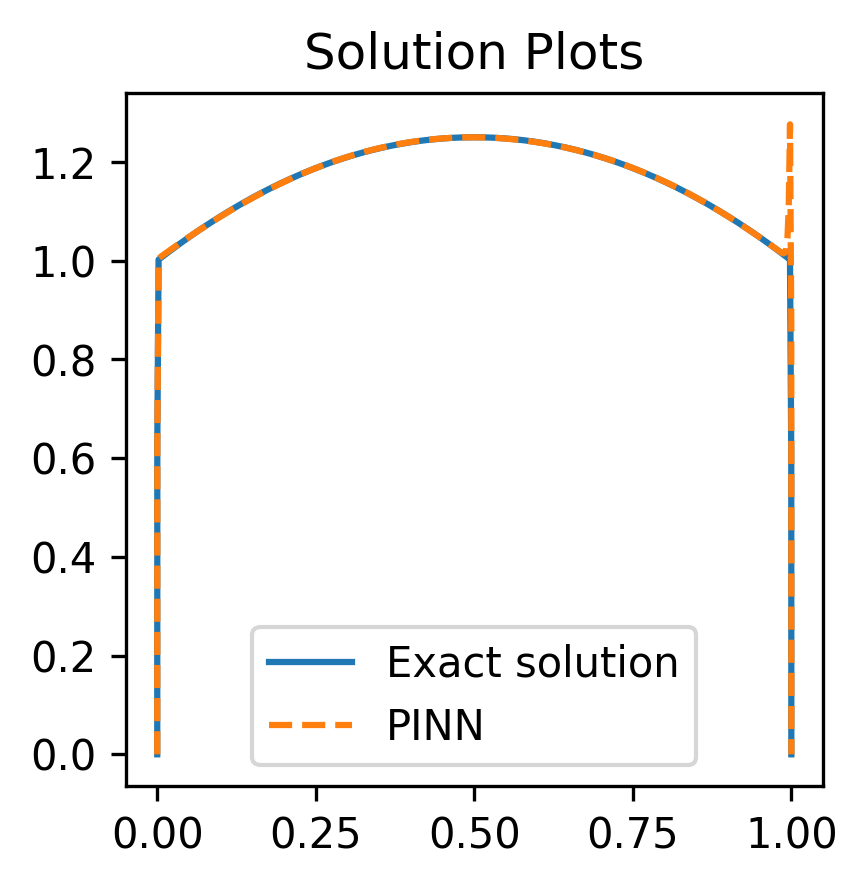}
        \includegraphics[width=\linewidth]{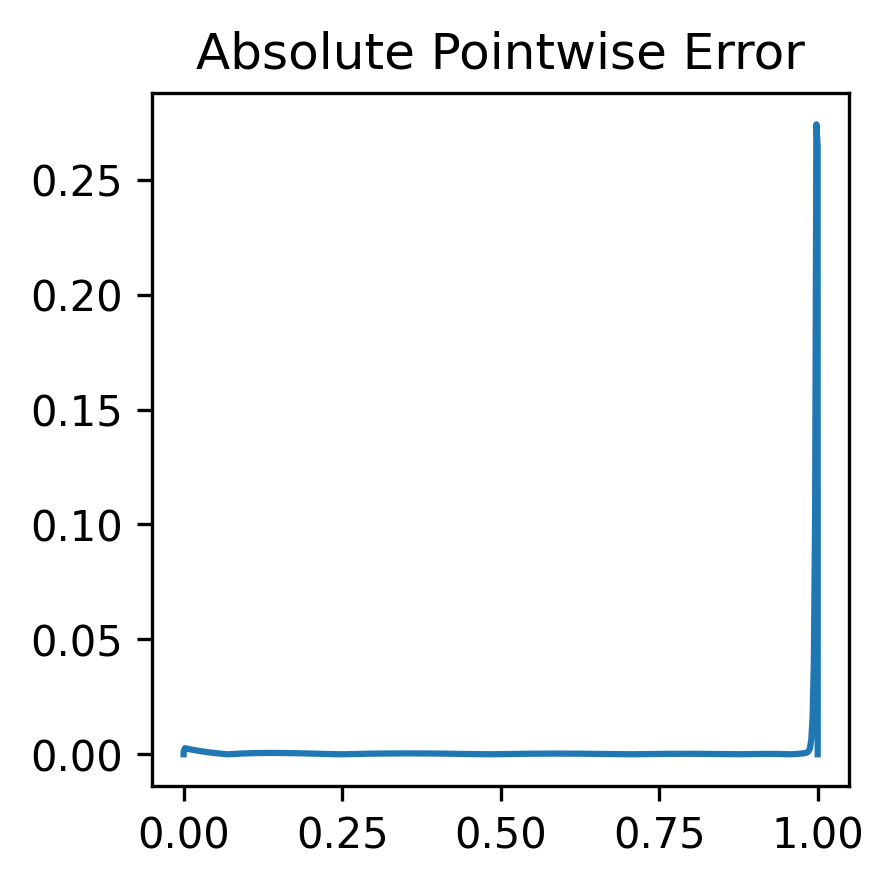}
        \caption{$\ep=10^{-8}$}
    \end{subfigure}

    \caption{Exact solutions and 
    PINN predictions 
    of $u_{1}^{\ep}$. }
    \label{plotu1PINN}
\end{figure}

\begin{figure}
    \centering
    
    \begin{subfigure}[t]{0.3\textwidth}
        \raggedleft
        \includegraphics[width=\linewidth]{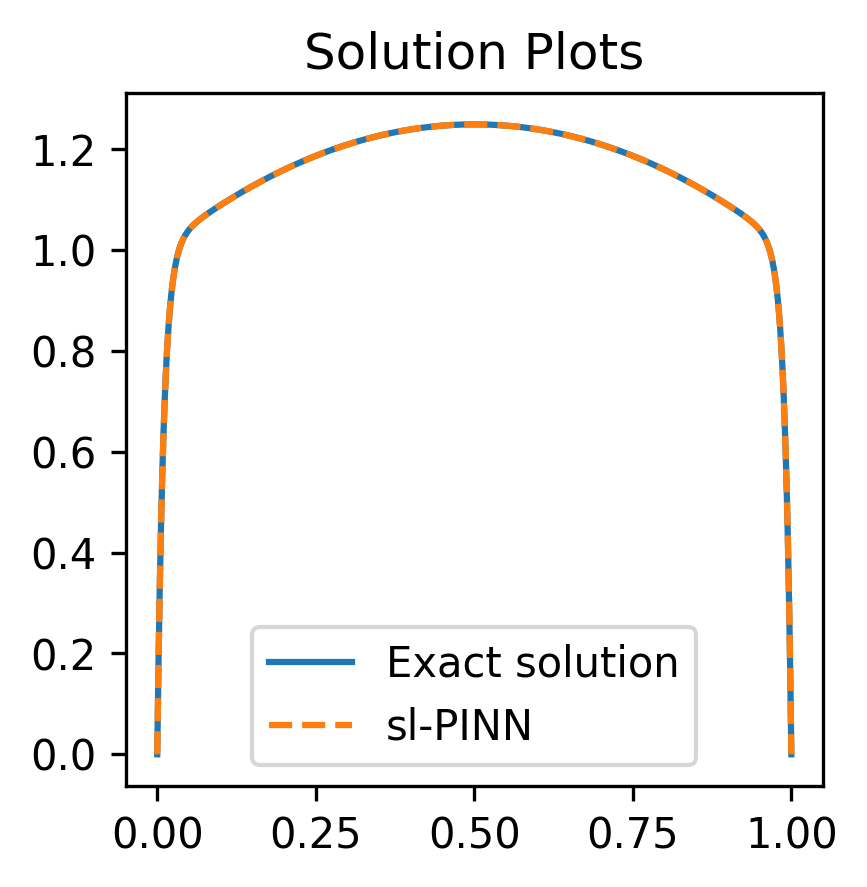}
        \includegraphics[width=\linewidth]{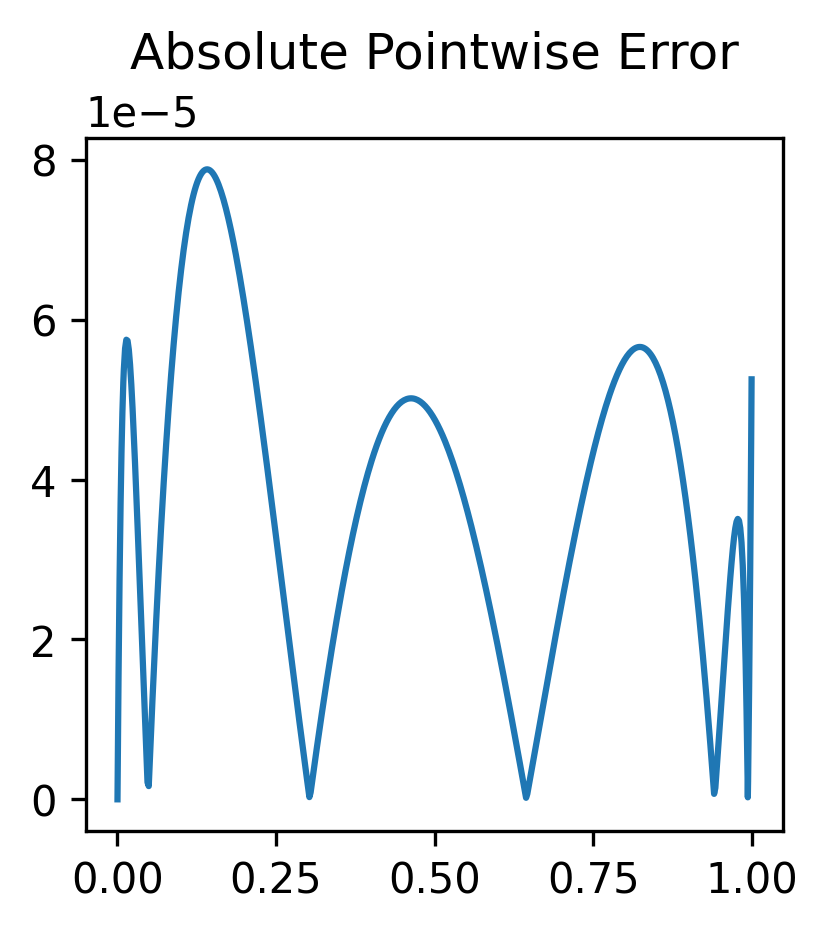}
        \caption{$\ep=10^{-4}$}
    \end{subfigure}%
    \hfill
    \begin{subfigure}[t]{0.3\textwidth}
        \raggedleft
        \includegraphics[width=\linewidth]{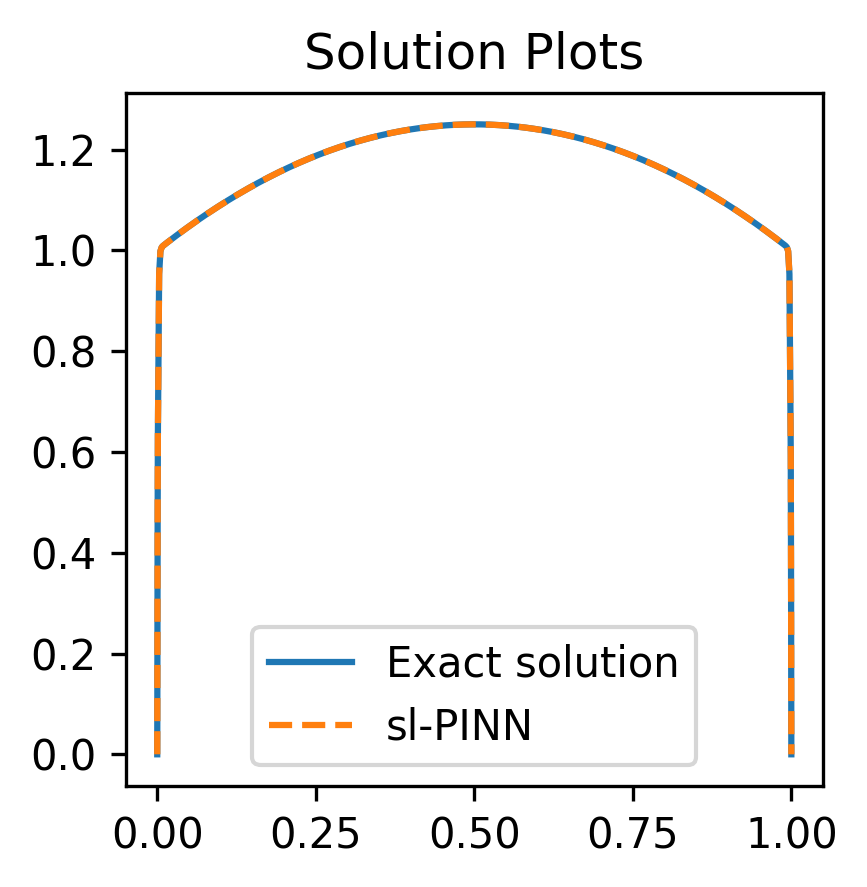}
        \includegraphics[width=\linewidth]{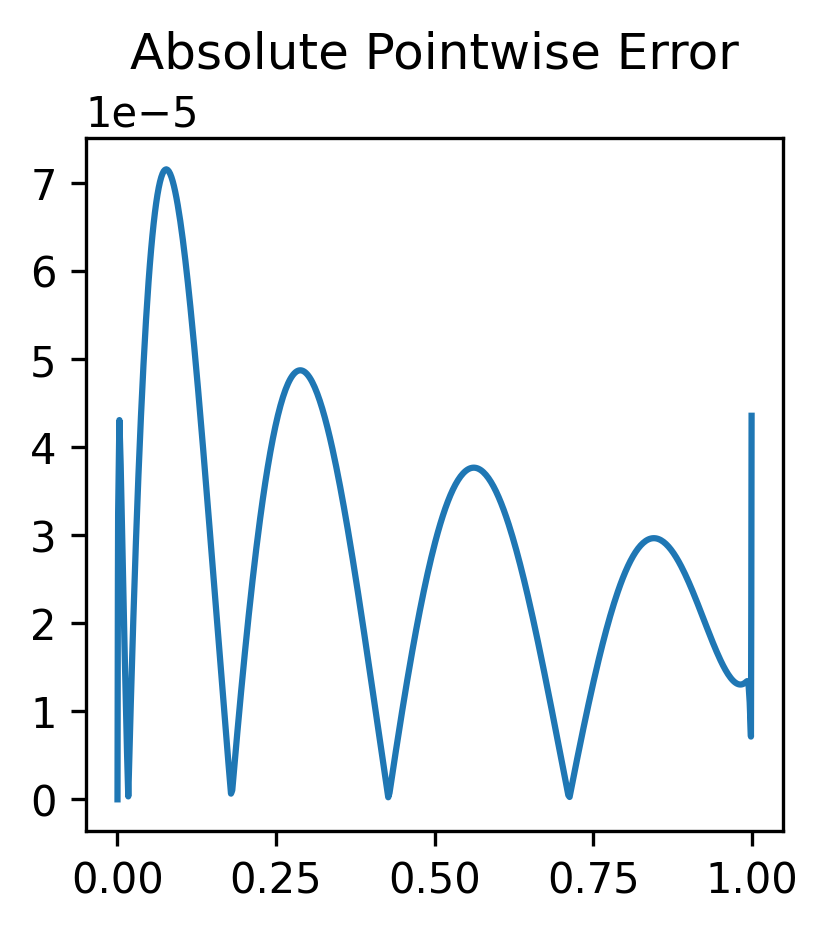}
        \caption{$\ep=10^{-6}$}
    \end{subfigure}%
    \hfill
    \begin{subfigure}[t]{0.3\textwidth}
        \raggedleft
        \includegraphics[width=\linewidth]{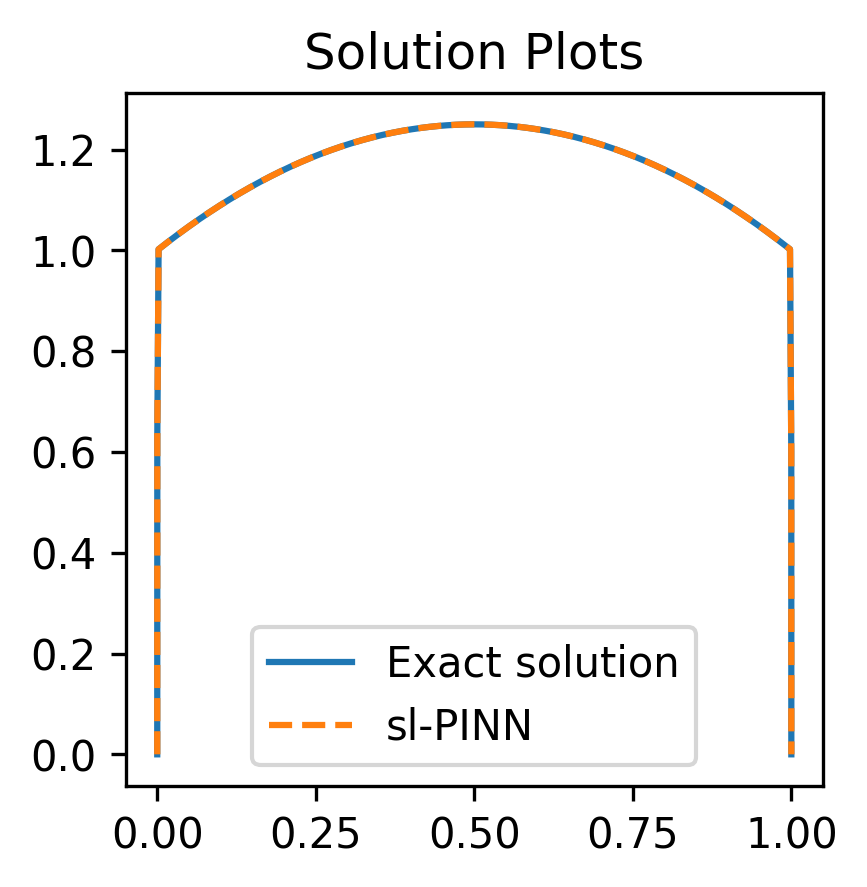}
        \includegraphics[width=\linewidth]{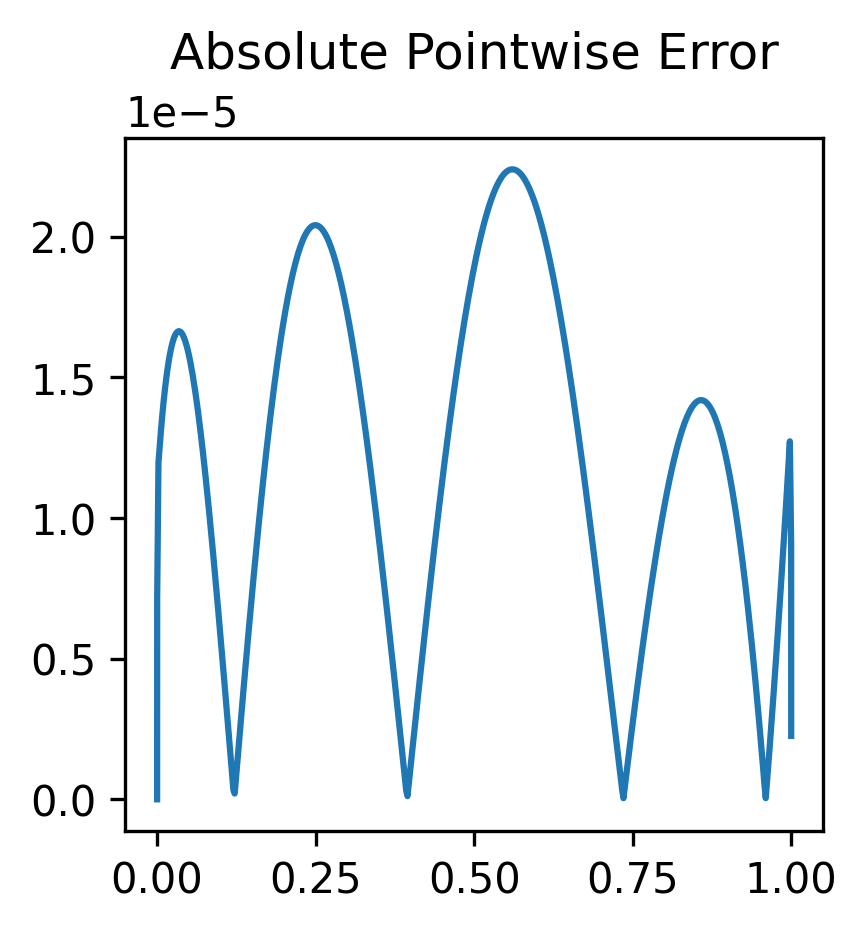}
        \caption{$\ep=10^{-8}$}
    \end{subfigure}

    \caption{Exact solutions and 
    {\em sl-PINN} predictions  
    of $u_{1}^{\ep}$.}
    \label{plotu1EPINN}
\end{figure}

\begin{figure}
    \centering
    \begin{subfigure}{\textwidth}
        \centering
        \includegraphics[width=\linewidth]{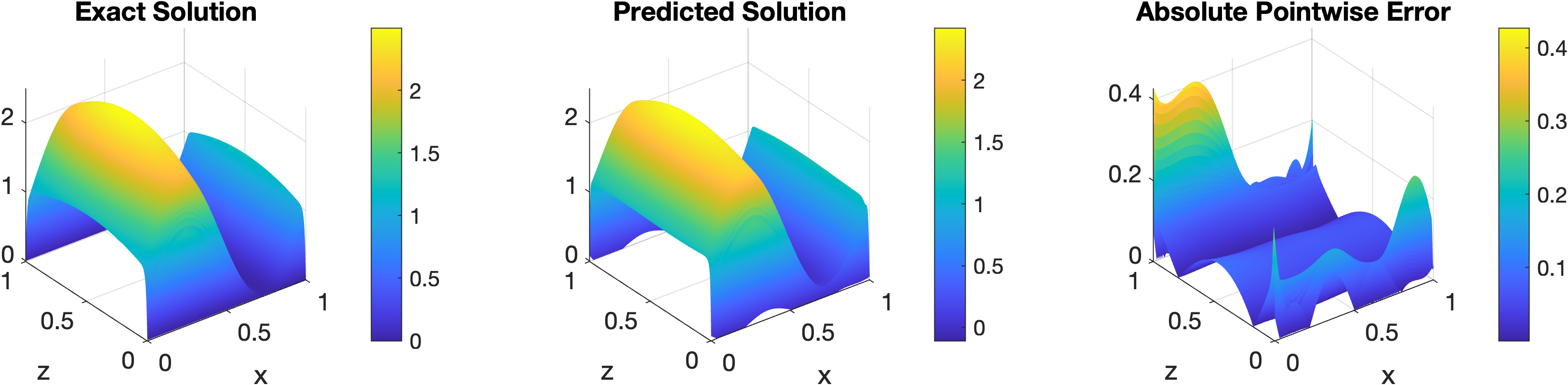}
        \caption{$\ep = 10^{-4}$}
    \end{subfigure}
    \begin{subfigure}{\textwidth}
        \centering
        \includegraphics[width=\linewidth]{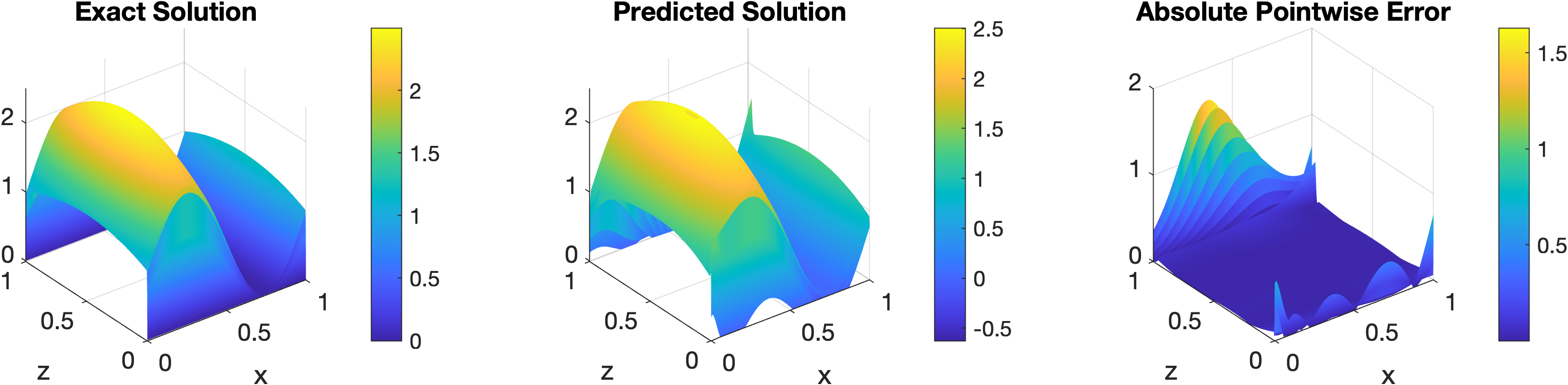}
        \caption{$\ep = 10^{-6}$}
    \end{subfigure}
    \begin{subfigure}{\textwidth}
        \centering
        \includegraphics[width=\linewidth]{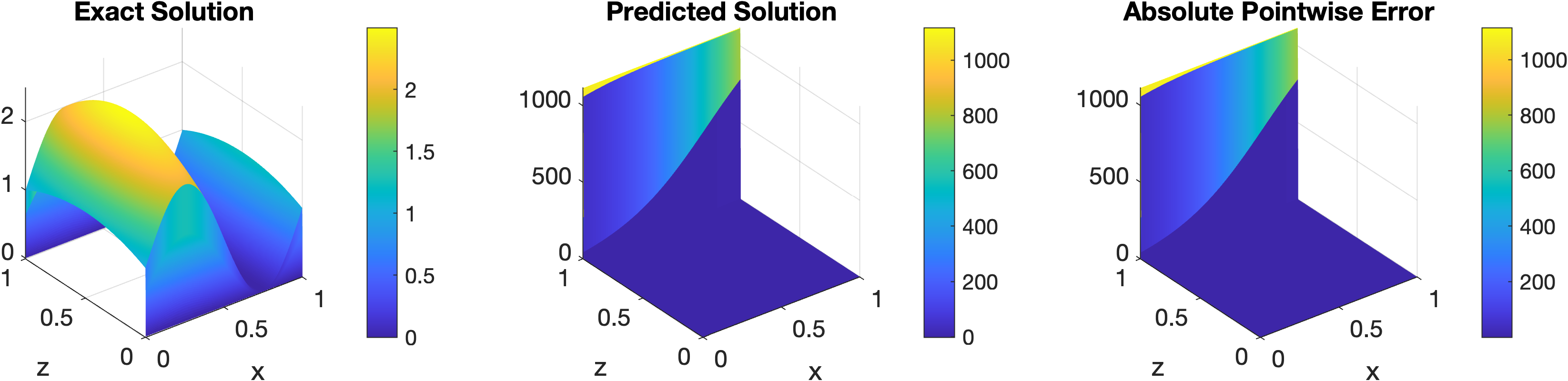}
        \caption{$\ep = 10^{-8}$}
    \end{subfigure}
    \caption{Exact solutions and 
    PINN predictions of $u_{2}^{\ep}$.}
    \label{plotu2PINN}
\end{figure}

\begin{figure}
    \centering
    \begin{subfigure}{\textwidth}
        \centering
        \includegraphics[width=\linewidth]{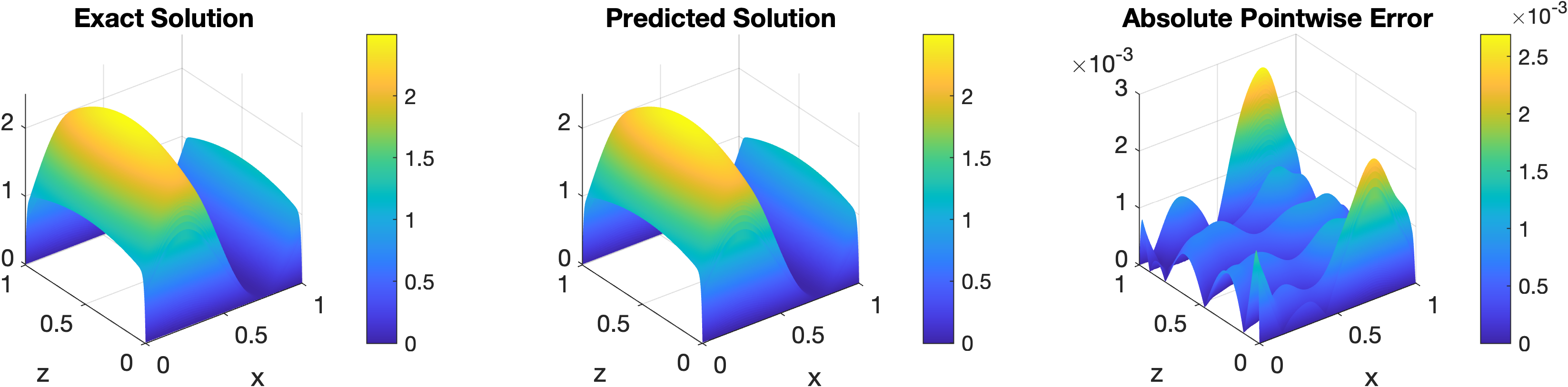}
        \caption{$\ep = 10^{-4}$}
    \end{subfigure}
    \begin{subfigure}{\textwidth}
        \centering
        \includegraphics[width=\linewidth]{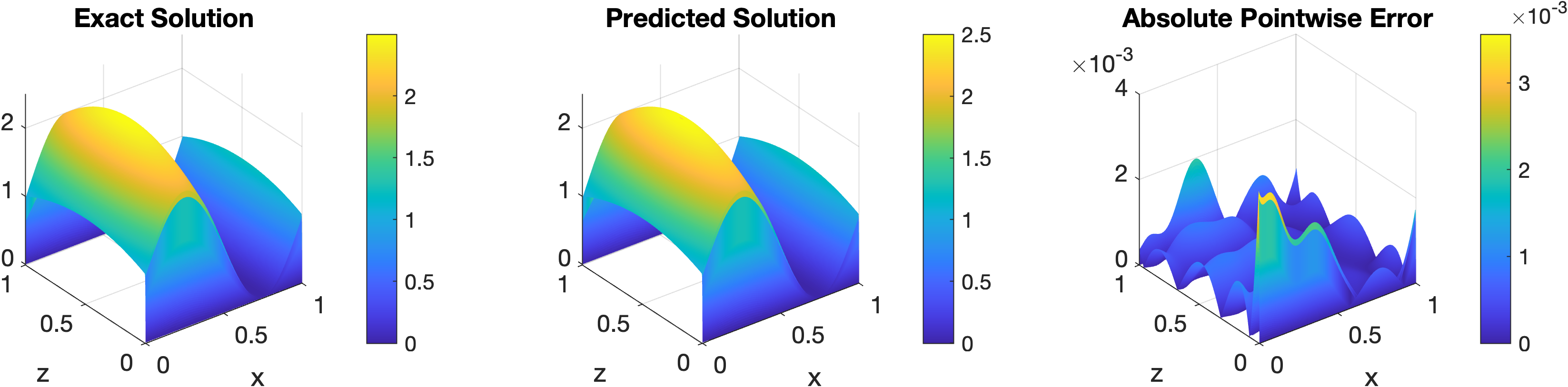}
        \caption{$\ep = 10^{-6}$}
    \end{subfigure}
    \begin{subfigure}{\textwidth}
        \centering
        \includegraphics[width=\linewidth]{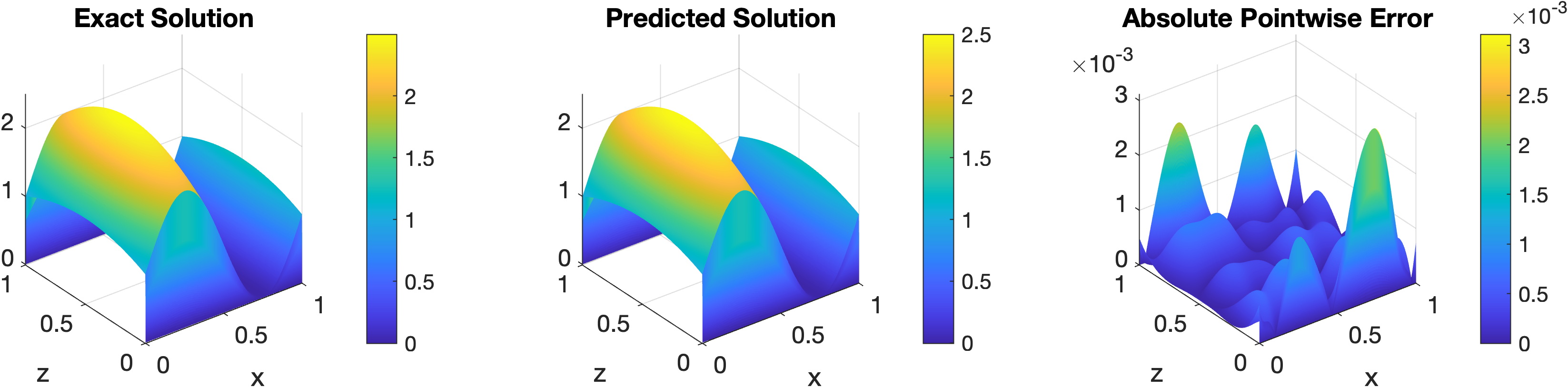}
        \caption{$\ep = 10^{-8}$}
    \end{subfigure}
    \caption{Exact solutions and 
    {\em sl-PINN} predictions of $u_{2}^{\ep}$.}
    \label{plotu2EPINN}
\end{figure}

We observe that our {\em sl-PINNs}  produce accurate predictions 
for every small viscosity $\ep$, while 
the (usual) PINNs fail to predict the solution as the solution becomes stiffer. 
For in-depth understanding the solutions' behavior, 
the profiles of the predicted solutions are shown in Figures \ref{plotu1PINN} and \ref{plotu1EPINN} for $u_{1}^{\ep}$ and in 
Figures \ref{plotu2PINN} and \ref{plotu2EPINN} for $u_{2}^{\ep}$. 
We notice that the smaller $\ep>0$ gets, the more challenging it becomes for the PINNs to capture the (singular) behavior around the boundary layers. 
In contrast, the {\em sl-PINNs} predict well the (singular) solution for every small $\ep$. Moreover, the training process of the {\em sl-PINN} is much faster than the PINN, as shown in the loss plots in Figures \ref{plotloss1} and \ref{plotloss2}. 
These observations indicate that 
our {\em sl-PINNs} are not only accurate but also efficient for  %learning methods for 
solving the stiff problems with a small viscosity $\ep$.

\begin{table}[]
\begin{tabular}{|l|ll|ll|}
\hline
                   & PINNs                                       &                             & {\em sl-PINNs}                                      &                             \\ \hline
                   & \multicolumn{1}{l|}{Relative $L^{2}$ error} & Relative $L^{\infty}$ error & \multicolumn{1}{l|}{Relative $L^{2}$ error} & Relative $L^{\infty}$ error \\ \hline
$\epsilon=10^{-3}$ & \multicolumn{1}{l|}{5.8741e-04}             & 1.5534e-03                  & \multicolumn{1}{l|}{1.0452e-05}             & 2.8015e-05                  \\
$\epsilon=10^{-4}$ & \multicolumn{1}{l|}{6.1017e-04}             & 1.8356e-03                  & \multicolumn{1}{l|}{3.8054e-05}             & 6.3043e-05                  \\
$\epsilon=10^{-5}$ & \multicolumn{1}{l|}{2.9297e-02}             & 1.4135e-01                  & \multicolumn{1}{l|}{1.0672e-05}             & 2.3559e-05                  \\
$\epsilon=10^{-6}$ & \multicolumn{1}{l|}{5.7345e-02}             & 3.5576e-01                  & \multicolumn{1}{l|}{2.8798e-05}             & 5.7216e-05                  \\
$\epsilon=10^{-7}$ & \multicolumn{1}{l|}{3.4765e-02}             & 1.6841e-01                  & \multicolumn{1}{l|}{1.2835e-05}             & 2.7431e-05                  \\
$\epsilon=10^{-8}$ & \multicolumn{1}{l|}{4.3961e-02}             & 2.1925e-01                  & \multicolumn{1}{l|}{1.1757e-05}             & 1.7917e-05                  \\ \hline
\end{tabular}
\caption{Comparison between PINNs and {\em sl-PINNs} for   $u_{1}^{\ep}$.}
\label{tableu1}
\end{table}

\begin{table}[]
\begin{tabular}{|l|ll|ll|}
\hline
                   & PINNs                                       &                             & {\em sl-PINNs}                                      &                             \\ \hline
                   & \multicolumn{1}{l|}{Relative $L^{2}$ error} & Relative $L^{\infty}$ error & \multicolumn{1}{l|}{Relative $L^{2}$ error} & Relative $L^{\infty}$ error \\ \hline
$\epsilon=10^{-3}$ & \multicolumn{1}{l|}{2.2551e-02}             & 1.1572e-01                  & \multicolumn{1}{l|}{5.3896e-04}             & 7.0597e-04                  \\
$\epsilon=10^{-4}$ & \multicolumn{1}{l|}{6.7786e-02}             & 1.7088e-01                  & \multicolumn{1}{l|}{5.0758e-04}             & 1.0775e-03                  \\
$\epsilon=10^{-5}$ & \multicolumn{1}{l|}{3.7465e-01}             & 3.8308e+00                  & \multicolumn{1}{l|}{4.6892e-04}             & 6.4305e-04                  \\
$\epsilon=10^{-6}$ & \multicolumn{1}{l|}{1.2475e-01}             & 6.5216e-01                  & \multicolumn{1}{l|}{3.9502e-04}             & 1.4244e-03                  \\
$\epsilon=10^{-7}$ & \multicolumn{1}{l|}{1.8177e-01}             & 1.9354e+00                  & \multicolumn{1}{l|}{6.8371e-04}             & 1.7341e-03                  \\
$\epsilon=10^{-8}$ & \multicolumn{1}{l|}{6.9807e+01}             & 4.4609e+02                  & \multicolumn{1}{l|}{3.7568e-04}             & 1.2459e-03                  \\ \hline
\end{tabular}
\caption{Comparison between PINNs and {\em sl-PINNs} for   $u_{2}^{\ep}$.}
\label{tableu2}
\end{table}

\section{ Numerical experiments for the vorticity}\label{sec_numerics_vor}

In this section, we present the framework and the experiments of our semi-analytic singular layer PINN ({\em sl-PINN}) methods for 
solving the vorticity equations (\ref{e:NSE_PP_Vor_eqn}), especially when   
the viscosity $\ep$ is small. 

Following the neural network (NN) architecture, introduced in Section \ref{sec_numerics}, 
we recall the conventional PINN methods for the vorticity $\blds{\omega}^\ep$ by constructing one network for each component of $\blds{\omega}^\ep$. 
% Similar to Section \ref{comparion_velocity}, 
% a comparison with conventional PINNs (PINNs) is also presented. 
% The neural network (NN) architecture was described in Section \ref{sec_numerics}. Three different NNs are used to approximate the solution of $\blds{\omega}^{\ep}=(\omega^{\ep}_{1}, \omega^{\ep}_{2}, \omega^{\ep}_{3})$ for each component. 
More precisely, we build a NN  for $\omega^{\ep}_{2}$ (sol. of (\ref{e:NSE_PP_Vor_eqn})$_2$ with (\ref{e:NSE_PP_Vor_BC})$_2$) first.  
Then, using this predicted solution for $\omega^{\ep}_{2}$ as well as the one for $u^\ep_1$, built by the conventional PINN in the previous section \ref{sec_numerics}, 
we construct the PINN predicted solution for $\omega^{\ep}_{3}$ (sol. of (\ref{e:NSE_PP_Vor_eqn})$_3$ with (\ref{e:NSE_PP_Vor_BC})$_3$), and, in a sequel,  
one for $\omega^{\ep}_{1}$ 
(sol. of (\ref{e:NSE_PP_Vor_eqn})$_1$ with (\ref{e:NSE_PP_Vor_BC})$_1$). 
% as 
% $\omega^{\ep}_{1}$ is obtained after solving $\omega^{\ep}_{2}$ and $\omega^{\ep}_{3}$. 
We let 
$\hat{\omega}_{1}=\hat{\omega}_{1}(x, z; \blds{\theta}_{1})$, $\hat{\omega}_{2}=\hat{\omega}_{2}(z; \blds{\theta}_{2})$, and $\hat{\omega}_{3}=\hat{\omega}_{3}(x, z; \blds{\theta}_{3})$ 
denote 
the outputs of each NNs respectively 
where each $\blds{\theta}_{i}$ is the collection of the weights and bias for $\hat{\omega}_{i}$, $i=1,2,3$. 

Using the equations and the boundary conditions in (\ref{e:NSE_PP_Vor_eqn}) - (\ref{e:NSE_PP_Vor_BC}), 
we train the conventional PINN by 
minimizing the loss function,  
$\mathcal{L}_{i}$ for $\hat{\omega}_{i}$, $i=1,2,3$ below:  
\begin{align} \label{loss2vorPINN}
    \mathcal{L}_{2}(\blds{\theta}_{2}; \mathcal{T}_{2})
    = \frac{1}{|\mathcal{T}_{2}|} \sum_{z \in \mathcal{T}_{2}} \left|  
        \hat{\omega}_{2} - \ep d^{2}_{z} \hat{\omega}_{2} - d_{z} f_{1}
    \right|^{2}
    + \left| d_{z} \hat{\omega}_{2}(0)+\frac{1}{\ep} f_{1}(0) \right|^{2}
    + \left| d_{z} \hat{\omega}_{2}(1)+\frac{1}{\ep} f_{1}(1) \right|^{2},
\end{align}
where $\mathcal{T}_{2} \subset [0,1]$ is the set of training points, 
\begin{equation} \label{loss3vorPINN}
    \begin{split}
    \mathcal{L}_{3}(\blds{\theta}_{3}; \mathcal{T}_{3}, \mathcal{T}_{\Gamma}, \mathcal{T}_{B})
    &= \frac{1}{|\mathcal{T}_{3}|} \sum_{(x,z)\in \mathcal{T}_{3}} \left|
        \hat{\omega}_{3} -\ep \Delta \hat{\omega}_{3} + \hat{u}_{1} \pa_x \hat{\omega}_{3} - \pa_{x} f_{2}
    \right|^{2} \\
    &+ \frac{1}{|\mathcal{T}_{\Gamma}|} \sum_{(x,z)\in \mathcal{T}_{\Gamma}} \left| \hat{\omega}_{3}(x, z=0) \right|^{2} + \left| \hat{\omega}_{3}(x, z=1) \right|^{2} \\
    &+ \frac{1}{|\mathcal{T}_{B}|} \sum_{(x,z)\in \mathcal{T}_{B}} \left| \hat{\omega}_{3}(x=0, z) - \hat{\omega}_{3}(x=L, z)\right|^{2},
    \end{split}
\end{equation}
where $\mathcal{T}_{3} \subset \Omega$, $\mathcal{T}_{\Gamma} \subset \Gamma$ and $\mathcal{T}_{B} \subset \partial\Omega \setminus \Gamma$ are the sets of training points, and
\begin{equation} \label{loss1vorPINN}
    \begin{split}
    \mathcal{L}_{1}(\blds{\theta}_{1}; \mathcal{T}_{1}, \mathcal{T}_{\Gamma}, \mathcal{T}_{B})
    &= \frac{1}{|\mathcal{T}_{1}|} \sum_{(x,z)\in \mathcal{T}_{1}} \left|
        \hat{\omega}_{1} - \ep \Delta \hat{\omega}_{1} - \hat{\omega}_{2}\hat{\omega}_{3} + \hat{u}_{1} \pa_{x} \hat{\omega}_{1} + \pa_{z} f_{2} 
    \right|^{2} \\
    &+ \frac{1}{|\mathcal{T}_{\Gamma}|} \sum_{(x,z)\in \mathcal{T}_{\Gamma}} \left| \pa_{z} \hat{\omega}_{1} (x, 0) - \frac{1}{\ep} f_{2}(x, 0) \right|^{2} + \left| \pa_{z} \hat{\omega}_{1} (x, 1) - \frac{1}{\ep} f_{2}(x, 1) \right|^{2} \\
    &+ \frac{1}{|\mathcal{T}_{B}|} \sum_{(x,z)\in \mathcal{T}_{B}} \left| \hat{\omega}_{1}(x=0, z) - \hat{\omega}_{1}(x=L, z)\right|^{2},
    \end{split}
\end{equation}
where $\mathcal{T}_{1} \subset \Omega$ is the set of training points. 

Notice that each of the loss functions $\mathcal{L}_2$ in (\ref{loss2vorPINN}) and $\mathcal{L}_1$ in (\ref{loss1vorPINN}) contains  a term of order $\ep^{-1}$, 
and hence they are not (expected to be) bounded uniformly in $\ep$ as $\ep \rightarrow 0$. 
% are unbounded as $\ep \to 0$ since the non-homogeneous Neumann conditions in order of $\ep^{-1}$ are involved. 
In fact,  because of this unbounded nature of the loss functions, 
the conventional PINN methods just fail to produce good predictions 
for the vorticity 
as we observe from the below numerical experiments in Section \ref{comparion_vor}. 
% This gives much more challenge for PINNs in minimizing the loss functions during the training process and can easily fall to an incorrect minimizer. The numerical experiments for PINNs can be found in Section \ref{comparion_vor}. 

In the following section, 
we propose our new semi-analytic singular layer PINNs ({\em sl-PINNs}) 
to predict well the vorticity $\blds{\omega}^\ep$. 
By embedding the singular profile of $\blds{\omega}^\ep$ in the {\em sl-PINNs} structure, 
we build the corresponding NNs to train the smooth part of $\blds{\omega}^\ep$ and finally obtain the accurate predicted solutions for 
$\blds{\omega}^\ep$. 
As it appears in Section \ref{comp_vor}, 
our {\em sl-PINNs} produce accurate predictions for 
$\blds{\omega}^\ep$ independent of a small $\ep$. 

% overcome the difficulty from the conventional PINNs. In particular, we construct our solutions motivated by considering the corrector to learn the behaviour near the boundary layers. The details of such solutions are formulated in the next subsection \ref{comp_vor}.

\subsection{{\em sl-PINN} predictions for the vorticity} \label{comp_vor}

To train predicted solutions for the vorticity 
$
\blds{\omega}^{\ep}
    =
        \text{curl } \blds{u}^\ep
    = \big(
			-\pa_z u^\ep_2, \,
			d_z u^\ep_1, \,
			\pa_x u^\ep_2
			\big)
$ 
(sol. of (\ref{e:NSE_PP_Vor_eqn}) - (\ref{e:NSE_PP_Vor_BC})), 
we follow the methodology introduced in Sections \ref{comp_u1} and \ref{comp_u2} and 
construct our 
{\em sl-PINNs} sequentially for the components 
${\omega}^{\ep}_2$, ${\omega}^{\ep}_3$, and then ${\omega}^{\ep}_1$. 
In order to incorporate 
the singular behavior near the boundary of $\blds{\omega}^{\ep}$ 
in our {\em sl-PINN} structure, 
we use below the profiles of the curl of the corrector $\tilde{\blds{\varphi}}_* = (\tilde{{\varphi}}_{1, *}, \, \tilde{{\varphi}}_{2, *}, \, 0)$, $* = L, R$, defined in (\ref{e:cor_PINN_1L}), (\ref{e:cor_PINN_1R}), (\ref{tilde tht2L}), and (\ref{tilde tht2R}). 
 
Note that 
the first component of the corrector $\tilde{{\varphi}}_{1, *}$, $*=L, R$,  is given in its explicit form as in (\ref{e:cor_PINN_1L}), (\ref{e:cor_PINN_1R}), and hence 
we embed the $d_z \tilde{{\varphi}}_{1, *}$, $*=L, R$, 
in our {\em sl-PINN} structure for ${\omega}^{\ep}_2$; see below in (\ref{tilde omega_2}) and (\ref{tilde omega_2-2}). 
On the other hand, 
because the explicit expression of second component $\tilde{{\varphi}}_{2, *}$, $*=L, R$,  
is not available,   
we train the viscous solution ${\omega}^{\ep}_3$ (and then ${\omega}^{\ep}_1$) together with training simultaneously 
the correctors 
$\pa_x \tilde{{\varphi}}_{2, *}$, $*=L, R$ (and then 
$-\pa_z \tilde{{\varphi}}_{2, *}$, $*=L, R$). 
This process is exactly what we employed for building 
our {\em sl-PINN} approximations for $u^\ep_2$ in (\ref{tilde u2}) - (\ref{tilde tht2R}), but here 
we need to work with a different type of boundary condition, e.g., Neumann boundary condition for ${\omega}^{\ep}_1$. 
\smallskip

We introduce the {\em sl-PINN} predictions  
$\widetilde{\blds{\omega}}
= (\widetilde{\omega}_1, \widetilde{\omega}_2, \widetilde{\omega}_3)$ enriched by the correctors. 
For the $\tilde{\omega}_{2}$, 
using the fact that 
the corrector for $\omega^{\ep}_{2} = d_z u^\ep_1$ is 
$d_z \varphi_1$ and 
its explicit expression is available, i.e.,  
$\ep^{-1/2} e^{-z/\sqrt{\ep}}$ near $z = 0$,  % known as the derivative of the corrector of $u^\ep_1$, 
we write the {\em sl-PINN} structure as
\begin{equation} \label{tilde omega_2}
    \tilde{\omega}_{2}(z; \blds{\tht}_{2}) 
        = 
            \hat{\omega}_{2}(z; \blds{\tht}_{2}) 
            + \tilde{\psi}_{2,L} (z; \blds{\tht}_{2}) 
            + \tilde{\psi}_{2,R} (z; \blds{\tht}_{2}) ,
\end{equation}
where the correctors are defined by
\begin{equation}\label{tilde omega_2-2}
    \begin{cases}
        \spacer
        \tilde{\psi}_{2,L}(z; \blds{\tht}_{2}) 
            = 
            \displaystyle 
            \big( 
                f_1(0) + \ep d_z \hat{\omega}_{2}(0; \blds{\tht}_{2}) 
            \big) 
            \dfrac{1}{\sqrt{\ep}}e^{-z/\sqrt{\ep}}, \\
        \tilde{\psi}_{2,R}(z; \blds{\tht}_{2}) 
            = 
            \displaystyle 
            \big( 
                - f_1(1) - {\ep} d_z \hat{\omega}_{2}(1; \blds{\tht}_{2}) 
            \big)
            \dfrac{1}{\sqrt{\ep}}e^{-(1-z)/\sqrt{\ep}}.
    \end{cases}
\end{equation}
Thanks to this setting of the correctors $\tilde{\psi}_{2,*}$, $* = L, R$, 
we notice that  
\begin{itemize}
    \item 
        $\tilde{\psi}_{2,L} \to 0$ as $z \to \infty$ 
        and 
        $\tilde{\psi}_{2,R} \to 0$ as $z \to -\infty$,   \vspace{2mm}
    \item   
        $d_z \tilde{\omega}_2
                                        = -\dfrac{1}{\ep} f_1$, up to an $e.s.t.$, at $z = 0, 1$.
\end{itemize} 
\smallskip

For the $\tilde{\omega}_{3}$, 
we use exactly the same construction as for the $u^\ep_2$ in (\ref{tilde u2}) - (\ref{tilde tht2R}), 
and write
% For the $\tilde{\omega}_{3}$, since the explicit form of the corrector is not well-known, we write the EPINNs solution of $\omega_{3}^{\ep}$, similar to the one (\ref{tilde u2}) of $u^{\ep}_{2}$, as
\begin{equation} \label{tilde omega_3}
    \tilde{\omega}_{3}(x, z; \blds{\Theta}_{3}) 
        = 
        \hat{\omega}_{3}(x, z; \blds{\tht}_{3}) 
        + 
        \tilde{\psi}_{3, L}(x, z; \blds{\tht}_{3}, \blds{\tht}_{3,s,L}) 
        + 
        \tilde{\psi}_{3, R}(x, z; \blds{\tht}_{3}, \blds{\tht}_{3,s,R}),
\end{equation}
with 
\begin{equation} 
\begin{cases}
    \spacer
    \tilde{\psi}_{3, L}(x, z; \blds{\tht}_{3}, \blds{\tht}_{3,s,L}) 
        =
        -\hat{\omega}_{3}(x, 0; \blds{\tht}_{3}) e^{-z/\sqrt{\ep}} 
        + 
        z e^{-z/\sqrt{\ep}} 
        \hat{\psi}_{3, L}(x, z; \blds{\tht}_{3, s, L}), \\
    \tilde{\psi}_{3, R}(x, z; \blds{\tht}_{3}, \blds{\tht}_{3,s,R}) 
        =
        -\hat{\omega}_{3}(x, 1; \blds{\tht}_{3}) e^{-(1-z)/\sqrt{\ep}} 
        + 
        (1-z) e^{-(1-z)/\sqrt{\ep}} 
        \hat{\psi}_{3, R}(x, z; \blds{\tht}_{3, s, R}).
\end{cases}
\end{equation}
That is, 
we 
use a neural network  
for training $\hat{\omega}_3$ and   $\hat{\psi}_3$ at the same time.
By our construction (and the analysis in (\ref{e:Lp_est_Tht_2_LR})), 
we observe that 
\begin{itemize}
    \item 
        $\tilde{\psi}_{3,L} \to 0$ as $z \to \infty$ 
        and 
        $\tilde{\psi}_{3,R} \to 0$ as $z \to -\infty$,   \vspace{2mm}
    \item   
        $\tilde{\omega}_3
                                        = 0$, up to an $e.s.t.$, at $z = 0, 1$.
\end{itemize} 
\smallskip

For the first component $\tilde{\omega}_{1}$, 
we combine the methods we employed for the other two components $\tilde{\omega}_{2}$ and $\tilde{\omega}_{3}$. 
Namely, we write our {\em sl-PINN} structure for $\tilde{\omega}_{1}$ by training  
the $\hat{\omega}_1$ (smooth part) and the $\hat{\psi}_1$ (singular part) simultaneously 
and by 
enforcing the Neumann boundary condition, 
in the form, 
\begin{equation} \label{tilde omega_1}
    \tilde{\omega}_{1}(x, z; \blds{\Theta}_{1}) 
        = 
        \hat{\omega}_{1}(x, z; \blds{\tht}_{1}) 
        + 
        \tilde{\psi}_{1, L}(x, z; \blds{\tht}_{1}, \blds{\tht}_{1,s,L} ) 
        + 
        \tilde{\psi}_{1, R}(x, z; \blds{\tht}_{1}, \blds{\tht}_{1,s,R} ),
\end{equation}
with 
\begin{equation} \label{tilde phi1}
\begin{cases}
    \tilde{\psi}_{1, L}(x, z; \blds{\tht}_{1}, \blds{\tht}_{1,s,L} ) 
        =&
        \spacer
        \Big(
                - f_2(x, 0) 
                + 
                {\ep} \pa_z \hat{\omega}_{1}(x, 0; \blds{\tht}_1)
            \Big)
            \dfrac{1}{\sqrt{\ep}} e^{-z/\sqrt{\ep}} \\
        &
        \spacer
        + 
        z^{2} e^{-z/\sqrt{\ep}} 
        \hat{\psi}_{1, L}(x, z; \blds{\tht}_{1, s, L}), \\
    \tilde{\psi}_{1, R}(x, z; \blds{\tht}_{1}, \blds{\tht}_{1,s,R} ) 
        =&
        \spacer
        \Big(
                f_2(x, 1) 
                - 
                {\ep} \pa_z \hat{\omega}_{1}(x, 1; \blds{\tht}_1)
            \Big)
            \dfrac{1}{\sqrt{\ep}} e^{-(1-z)/\sqrt{\ep}} \\
        &+ 
        (1-z)^{2} e^{-(1-z)/\sqrt{\ep}} 
        \hat{\psi}_{1, R}(x, z; \blds{\tht}_{1, s, R}).
\end{cases}
\end{equation}
% Notice that we take square: $z^{2}$ and $(1-z)^{2}$ on the last terms of the right-hand side of (\ref{tilde phi1}) so that (\ref{tilde omega_1}) satisfies the Neumann boundary condition.
We then observe that 
\begin{itemize}
    \item 
        $\tilde{\psi}_{1,L} \to 0$ as $z \to \infty$ 
        and 
        $\tilde{\psi}_{1,R} \to 0$ as $z \to -\infty$,   \vspace{2mm}
    \item   
        $\pa_z \tilde{\omega}_1
                                        =  \dfrac{1}{\ep} f_2$, up to an $e.s.t.$, at $z = 0, 1$.
\end{itemize} 
\smallskip

In the constructions of $\tilde{\omega}_i$, $i = 1, 2,3,$ above, we used the notation for the learning parameters, 
\begin{equation}\label{e:parameters_r_s}
    \blds{\Theta}_{i} 
        =
        \{
            \blds{\tht}_{i}, \, 
            \blds{\tht}_{i, s, L}, \, 
            \blds{\tht}_{i, s, R}
        \}, 
    \quad
        i = 1, 2,3. 
\end{equation}

% Note that the enriched PINN approximation 
% $\widetilde{\blds{\omega}}$ satisfies the boundary condition (\ref{e:NSE_PP_Vor_BC}), 
% that is, 
% \begin{equation*} 
%                                         \pa_z \widetilde{\omega}_1
%                                         = \dfrac{1}{\ep} f_2,
%                                         \quad
%                                         d_z \widetilde{\omega}_2
%                                         = -\dfrac{1}{\ep} f_1,
%                                         \quad
%                                         \widetilde{\omega}_3
%                                         = 0,
%                                         \quad
%                                             \text{on } \Gamma.
% \end{equation*}

Because the boundary condition for $\tilde{\blds{\omega}}$ is already taken into account in the construction of our {\em sl-PINNs},   
we define the loss functions for {\em sl-PINN}s as follows:
\begin{align} \label{loss2vorEPINN}
    \mathcal{L}_{2}(\blds{\theta}_{2}; \mathcal{T}_{2})
    = \frac{1}{|\mathcal{T}_{2}|} \sum_{z \in \mathcal{T}_{2}} \left|  
        \tilde{\omega}_{2} - \ep d^{2}_{z} \tilde{\omega}_{2} - d_{z} f_{1}
    \right|^{2},
\end{align}
where $\mathcal{T}_{2} \subset [0,1]$ is the set of training points, 
\begin{equation} \label{loss3vorEPINN}
    \begin{split}
    \mathcal{L}_{3}(\blds{\Theta}_{3}; \mathcal{T}_{3}, \mathcal{T}_{\Gamma}, \mathcal{T}_{B})
    &= \frac{1}{|\mathcal{T}_{3}|} \sum_{(x,z)\in \mathcal{T}_{3}} \left|
        \tilde{\omega}_{3} -\ep \Delta \tilde{\omega}_{3} + \tilde{u}_{1} \pa_x \tilde{\omega}_{3} - \pa_{x} f_{2}
    \right|^{2} \\
    &+ \frac{1}{|\mathcal{T}_{B}|} \sum_{(x,z)\in \mathcal{T}_{B}} \left| \tilde{\omega}_{3}(x=0, z) - \tilde{\omega}_{3}(x=L, z)\right|^{2},
    \end{split}
\end{equation}
where $\mathcal{T}_{3} \subset \Omega$, $\mathcal{T}_{\Gamma} \subset \Gamma$ and $\mathcal{T}_{B} \subset \partial\Omega \setminus \Gamma$ are the sets of training points, and
\begin{equation} \label{loss1vorEPINN}
    \begin{split}
    \mathcal{L}_{1}(\blds{\Theta}_{1}; \mathcal{T}_{1}, \mathcal{T}_{\Gamma}, \mathcal{T}_{B})
    &= \frac{1}{|\mathcal{T}_{1}|} \sum_{(x,z)\in \mathcal{T}_{1}} \left|
        \tilde{\omega}_{1} - \ep \Delta \tilde{\omega}_{1} - \tilde{\omega}_{2}\tilde{\omega}_{3} + \tilde{u}_{1} \pa_{x} \tilde{\omega}_{1} + \pa_{z} f_{2} 
    \right|^{2} \\
    &+ \frac{1}{|\mathcal{T}_{B}|} \sum_{(x,z)\in \mathcal{T}_{B}} \left| \tilde{\omega}_{1}(x=0, z) - \tilde{\omega}_{1}(x=L, z)\right|^{2},
    \end{split}
\end{equation}
where $\mathcal{T}_{1} \subset \Omega$ is the set of training points. 

% Here, the boundary condition (\ref{e:NSE_PP_Vor_BC}) are no need to enforce in the EPINNs loss functions since they are automatically satisfied by the EPINNs approximation. Moreover, one can verify 

It is noteworthy to point out that 
the loss functions in (\ref{loss2vorEPINN}), (\ref{loss3vorEPINN}) and (\ref{loss1vorEPINN}) 
stay bounded as the viscosity gets small, i.e., as $\ep \to 0$. 
Thanks to this important feature, as we shall see below in Section \ref{comparion_vor}, 
 our {\em sl-PINN}s produce stable and accurate predictions 
 for the vorticity, independent of the viscosity $\ep$.
% This make the training process easier to attain an expected minimizer. The numerical experiments can be found in the next subsection \ref{comparion_vor}.

\subsection{Comparison between the conventional PINNs and the {\em sl-PINNs} for the vorticity $\blds{\omega}^\ep$}\label{comparion_vor}

To perform the numerical experiments, 
we first write the exact solution for $\blds{\omega}^{\ep}$ 
by taking curl of $\blds{u}^{\ep}$ in (\ref{exactu1}) and (\ref{exactu2}):
\begin{equation}\label{e:vor_exact}
    \begin{cases}
        \spacer
        \omega^{\ep}_{1}(x, z) = -\omega^{\ep}_{2}(z) (1+\sin(2\pi x)), \\
        \spacer
        \displaystyle \omega^{\ep}_{2}(z) = -(1-2\ep) \frac{1-e^{-\frac{1}{\sqrt{\ep}}}}{1-e^{-\frac{2}{\sqrt{\ep}}}} \left( \frac{-1}{\sqrt{\ep}} e^{-\frac{z}{\sqrt{\ep}}} + \frac{1}{\sqrt{\ep}} e^{-\frac{1-z}{\sqrt{\ep}}} \right) + 1-2z, \\
        \omega^{\ep}_{3}(x, z) = 2\pi u^{\ep}_{1}(z) \cos(2\pi x).
    \end{cases}
\end{equation}

We infer from (\ref{e:vor_exact}) and the boundary layer analysis results in Theorem \ref{t:PPF} that 
\begin{itemize}
    \item 
        Asymptotic behavior of $\omega^\ep_3$ is as much singular as that of $u^\ep_1$ as $\ep \rightarrow 0$. That is, as small as 
        $e^{-\eta/\sqrt{\ep}}$, $\eta = z$ or $1-z$, near the boundary. 
    \item   
        The $\omega^\ep_1$ and $\omega^\ep_2$ behave near the boundary as the approximation of identity, $\ep^{-1/2}\, e^{-\eta/\sqrt{\ep}}$, $\eta = z$ or $1-z$, which converges as $\ep \rightarrow 0$ to 
        a positive measure $\delta_\Gamma$ supported on the boundary $\Gamma$.  
\end{itemize}
The asymptotic behavior at a small viscosity of the vorticity ($\omega^\ep_1$ and $\omega^\ep_2$) is more singular than that of the velocity vector field, and hence 
it is more challenging to approximate/predict  
the vorticity than the velocity at a small viscosity 
by employing PINN or any other classical numerical methods, e.g., Finite Elements, Finite Differences, Finite Volumes, or discontinuous Galerkin methods, and so on.  
 % Notice that $\omega^{\ep}_{2}$ and $\omega^{\ep}_{1}$ blow up near the boundary at $z=0, 1$ as fast as the approximation of identity, e.g. $1/\sqrt{\ep} e^{-z/\sqrt{\ep}}$ near $z=0$. The problem becomes stiffer than the one for velocity. 

We aim to verify the accurate performance of our {\em sl-PINN} methods to approximate the solution $\blds{\omega}^{\ep}$ to (\ref{e:NSE_PP_Vor_eqn}) without causing an excessive computation cost. 
The experimental settings are the same as those for the velocity in Section \ref{comparion_velocity}. 
We recall, in particular, that we use single-layer NNs with a small amount of training data to remain a low computational cost while maintaining sharp accuracy.  
The training data set $\mathcal{T}_{1}\subset \Omega$ for $\omega^{\ep}_{1}$ is the same as $\mathcal{T}_{3}\subset \Omega$ for $\omega^{\ep}_{3}$. 
The experiments for conventional PINNs are conducted for comparison with our {\em sl-PINNs}. 
We present here the results of $\omega^{\ep}_{2}$ and $\omega^{\ep}_{1}$ only since the equation for $\omega^{\ep}_{3}$ is identical to that of $u^{\ep}_{2}$, whose numerical simulations are well-studied in Section \ref{comparion_velocity}.

\begin{figure}
    \begin{subfigure}{0.5\textwidth}
        \centering
        \includegraphics[width=\linewidth]{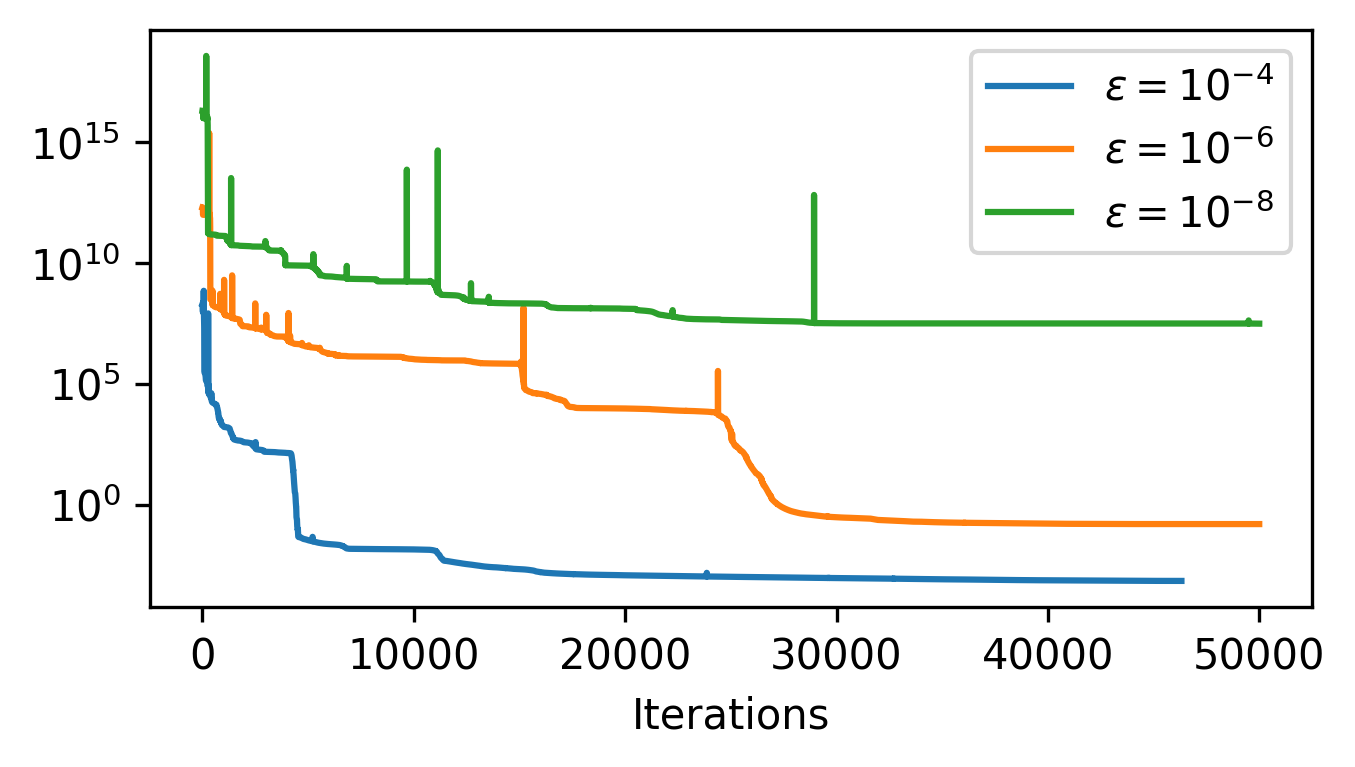}
        \caption{PINNs}
    \end{subfigure}%
    \begin{subfigure}{0.5\textwidth}
        \centering
        \includegraphics[width=\linewidth]{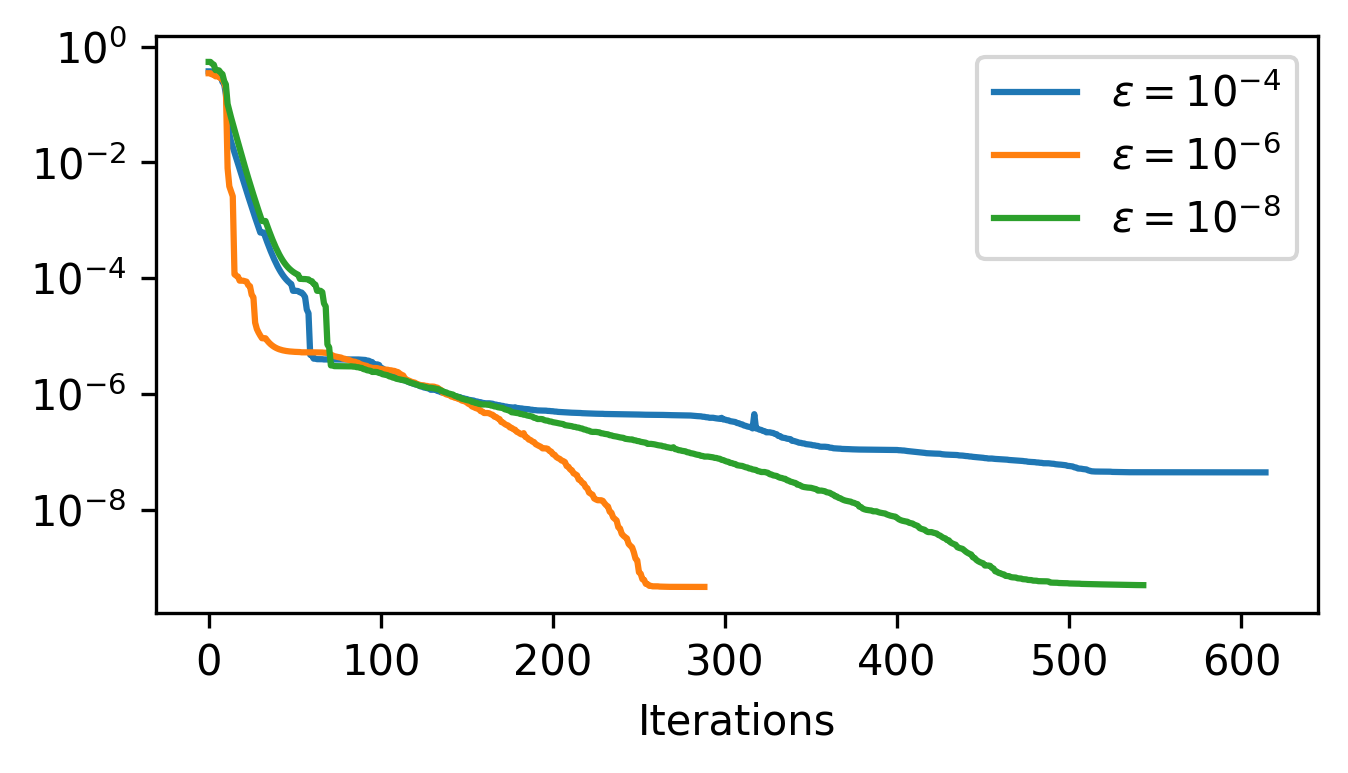}
        \caption{{\em sl-PINNs}}
    \end{subfigure}
    \caption{The loss values during the training process for $\omega_{2}^{\ep}$.}
    \label{plotloss1vor}
\end{figure}

\begin{figure}
    \begin{subfigure}{0.5\textwidth}
        \centering
        \includegraphics[width=\linewidth]{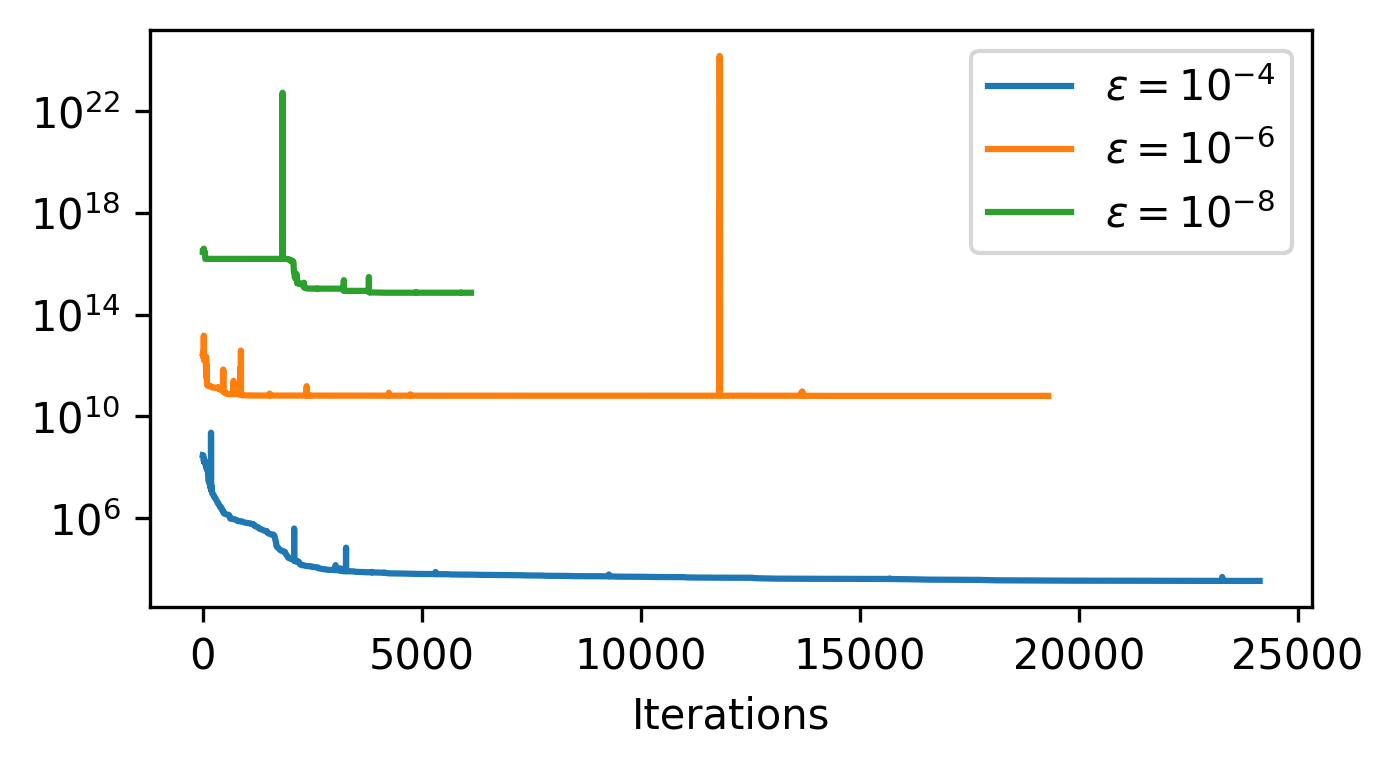}
        \caption{PINNs}
    \end{subfigure}%
    \begin{subfigure}{0.5\textwidth}
        \centering
        \includegraphics[width=\linewidth]{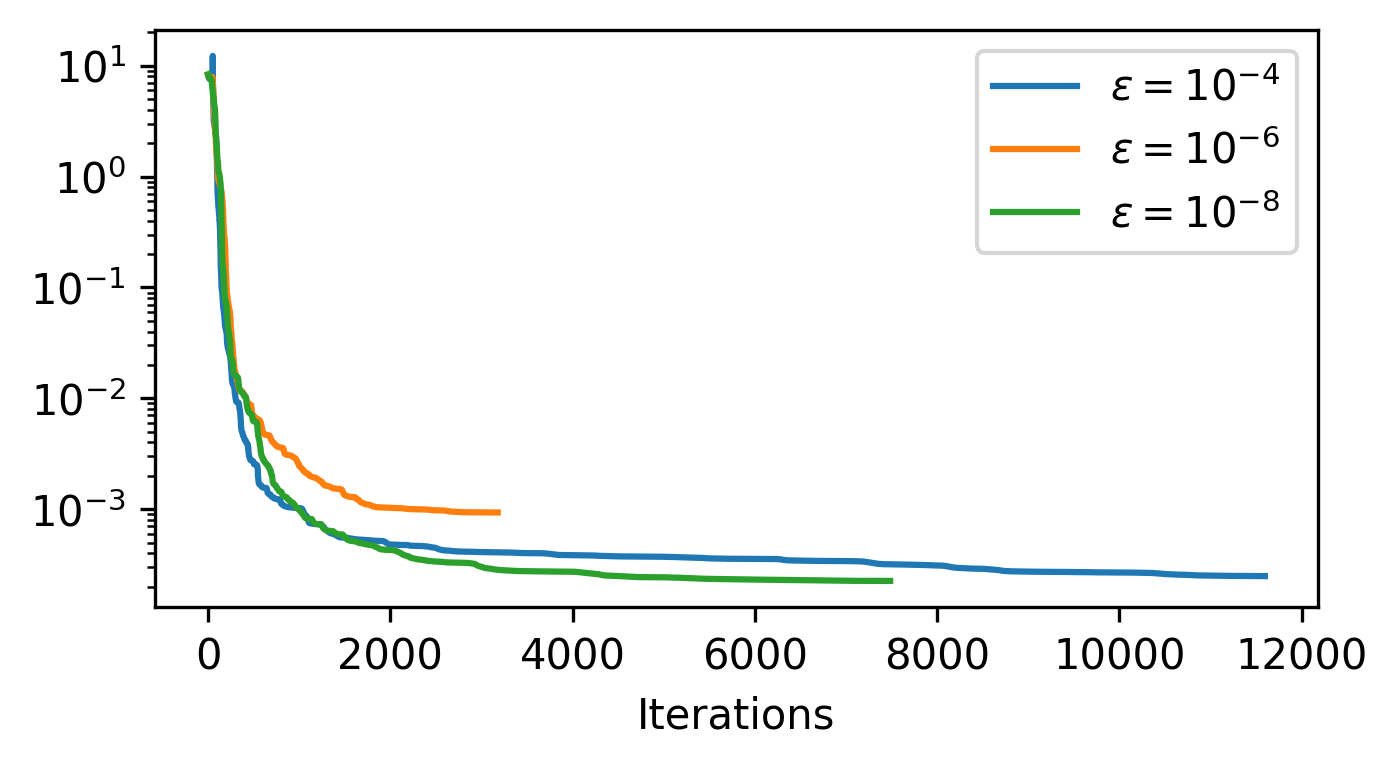}
        \caption{{\em sl-PINNs}}
    \end{subfigure}
    \caption{The loss values during the training process for $\omega_{1}^{\ep}$}
    \label{plotloss2vor}
\end{figure}

First, the loss value during the training process for $\omega^{\ep}_{2}$ and $\omega^{\ep}_{1}$ are shown in Figures \ref{plotloss1vor} and \ref{plotloss2vor} respectively. 
As expected, the conventional PINNs just fail to minimize the loss functions (\ref{loss2vorPINN}) and (\ref{loss1vorPINN}) 
because they started with an extremely large value in their minimization, results in an incorrect minimizer. 
On the contrary, our new {\em sl-PINNs}' loss functions (\ref{loss2vorEPINN}) and (\ref{loss1vorEPINN}) are bounded independent of $\ep$,  
and thus they are well minimized for every small $\ep$. 
Moreover, the efficiency of {\em sl-PINNs} 
outperforms that of PINNs, 
as it reaches a minimum with far fewer iterations.

\begin{figure}
    \centering
    
    \begin{subfigure}[t]{0.3\textwidth}
        \raggedleft
        \includegraphics[width=\linewidth]{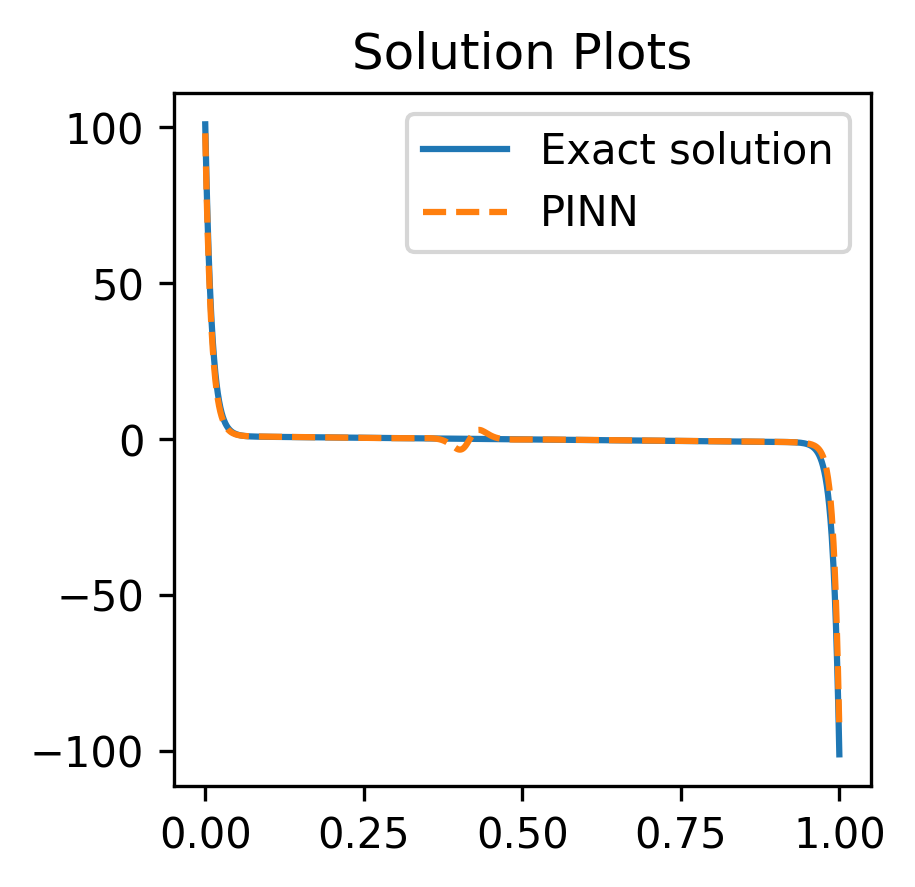}
        \includegraphics[width=0.9\linewidth]{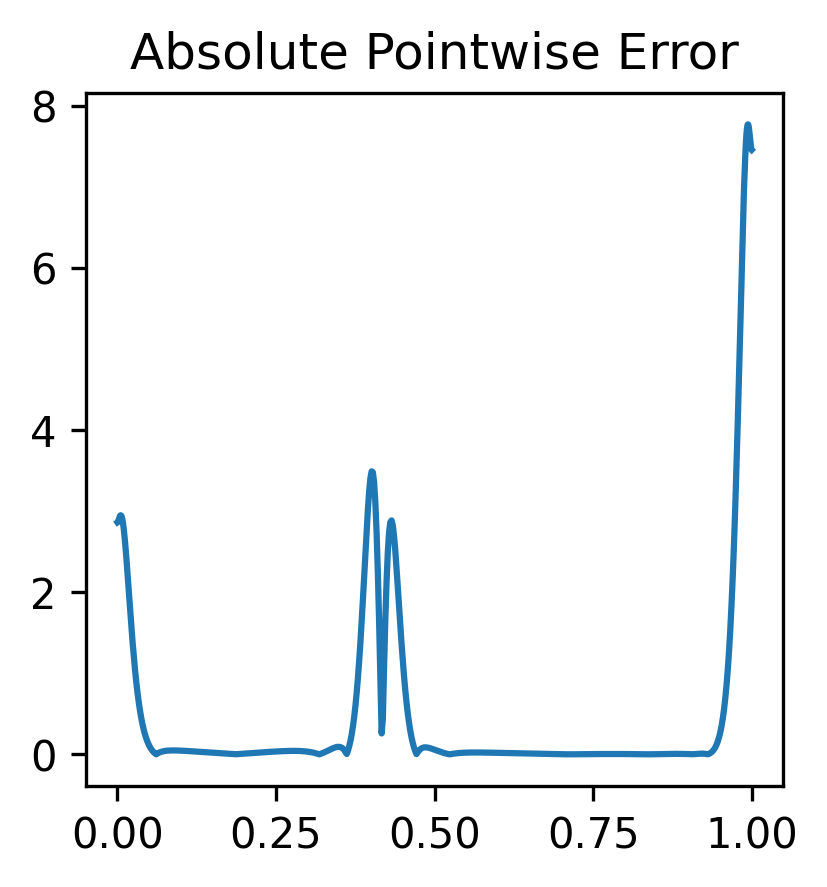}
        \caption{$\ep=10^{-4}$}
    \end{subfigure}%
    \hfill
    \begin{subfigure}[t]{0.3\textwidth}
        \raggedleft
        \includegraphics[width=\linewidth]{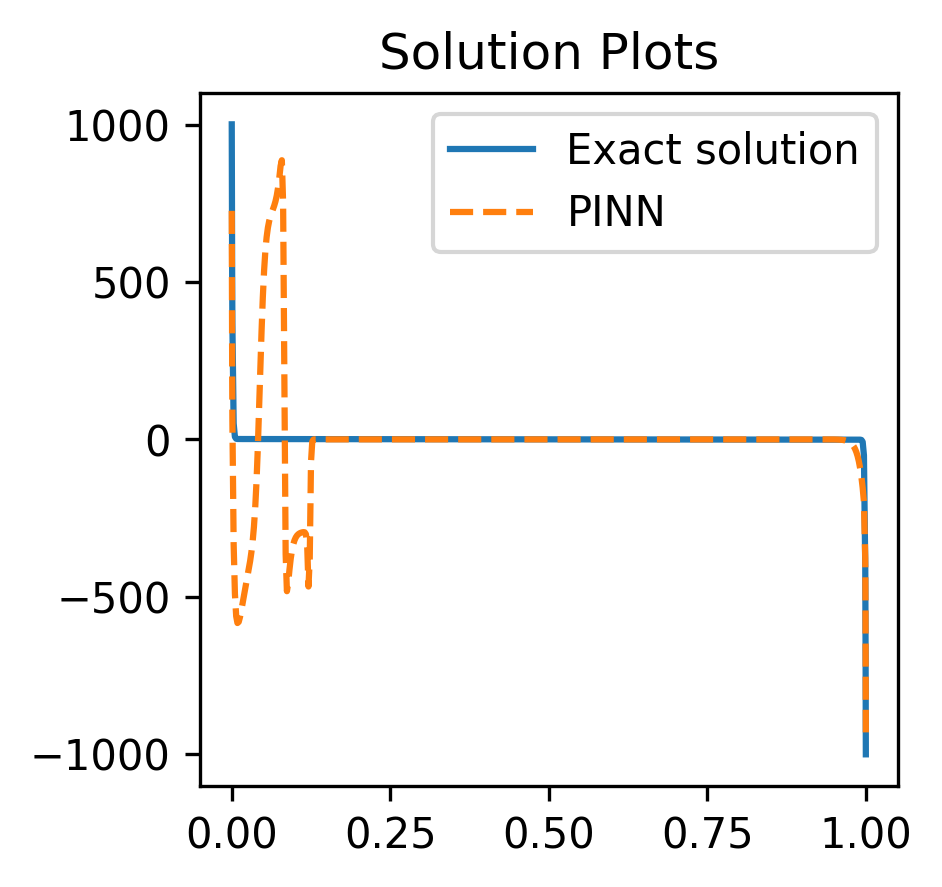}
        \includegraphics[width=0.95\linewidth]{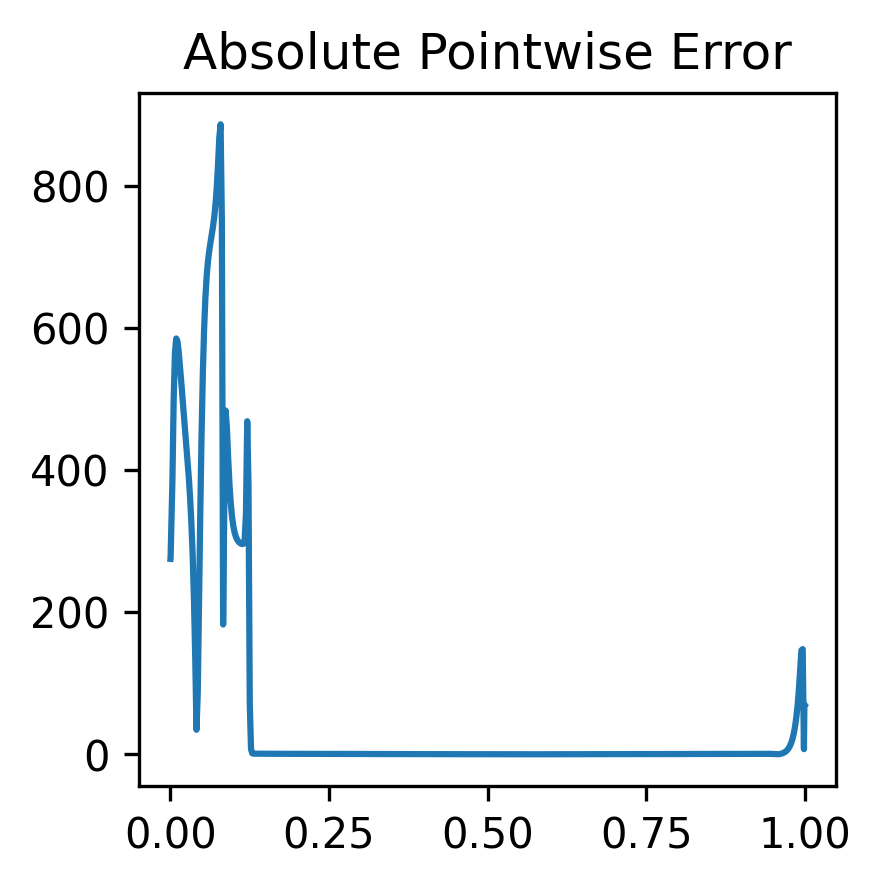}
        \caption{$\ep=10^{-6}$}
    \end{subfigure}%
    \hfill
    \begin{subfigure}[t]{0.3\textwidth}
        \raggedleft
        \includegraphics[width=\linewidth]{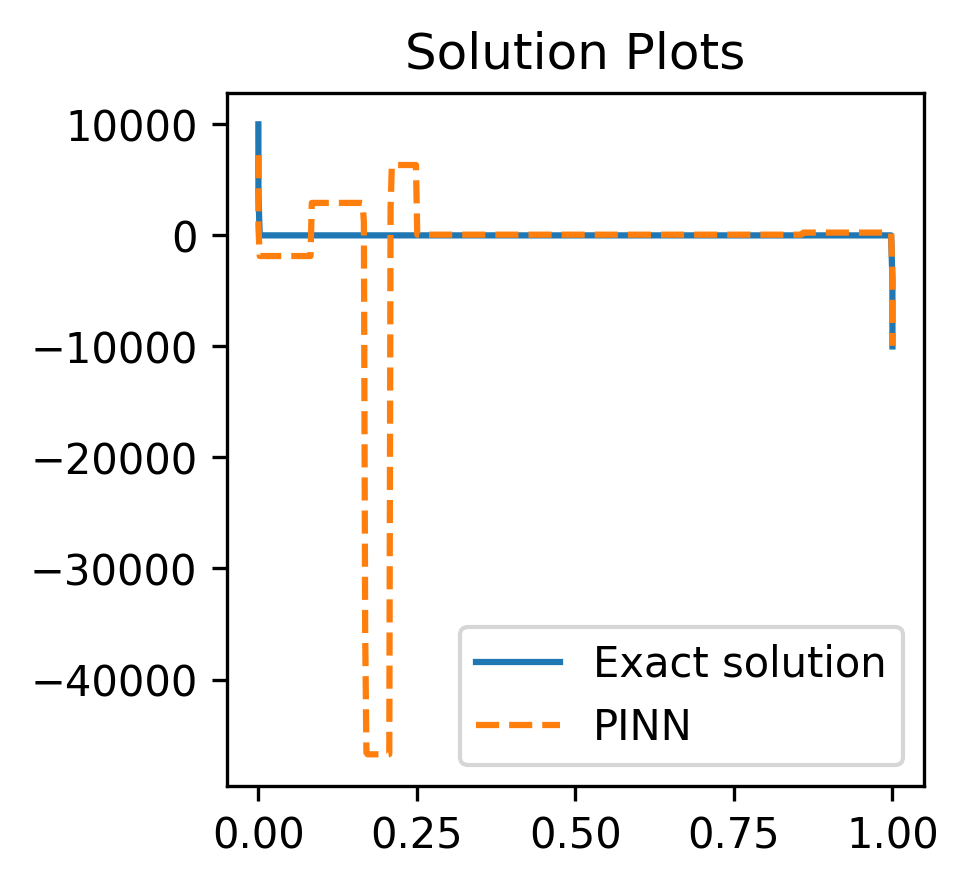}
        \includegraphics[width=\linewidth]{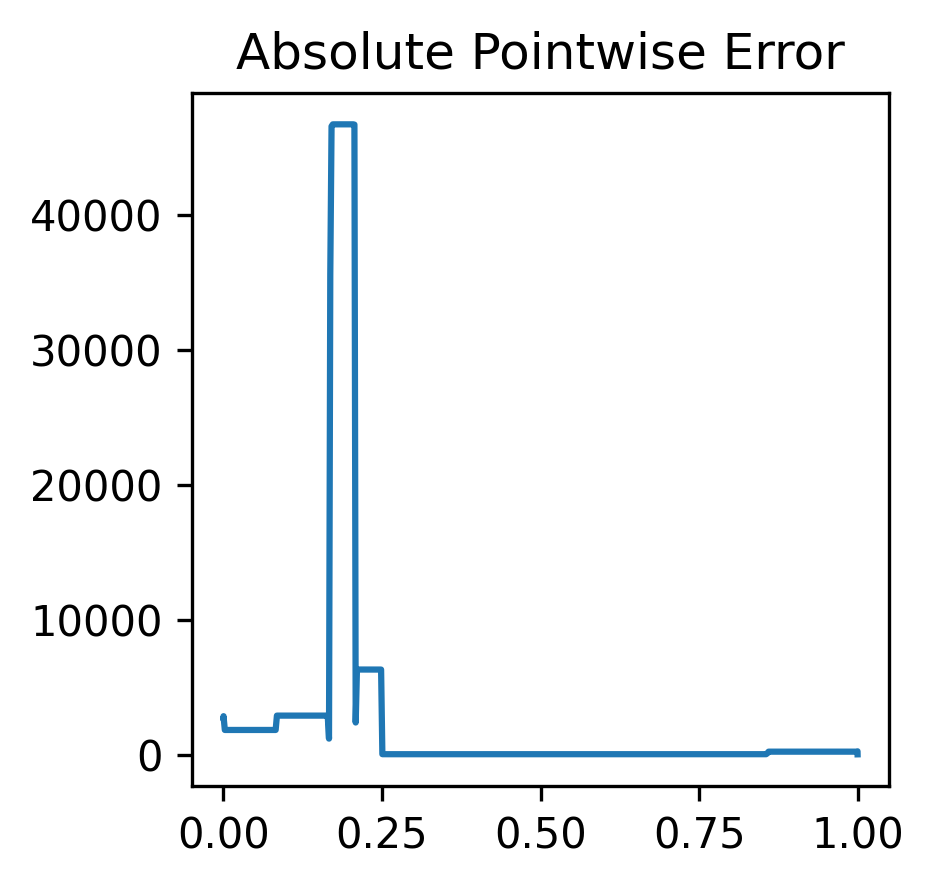}
        \caption{$\ep=10^{-8}$}
    \end{subfigure}

    \caption{Exact solutions and PINN predictions of $\omega_{2}^{\ep}$.}
    \label{plotw2PINN}
\end{figure}

\begin{figure}
    \centering
    
    \begin{subfigure}[t]{0.3\textwidth}
        \raggedleft
        \includegraphics[width=\linewidth]{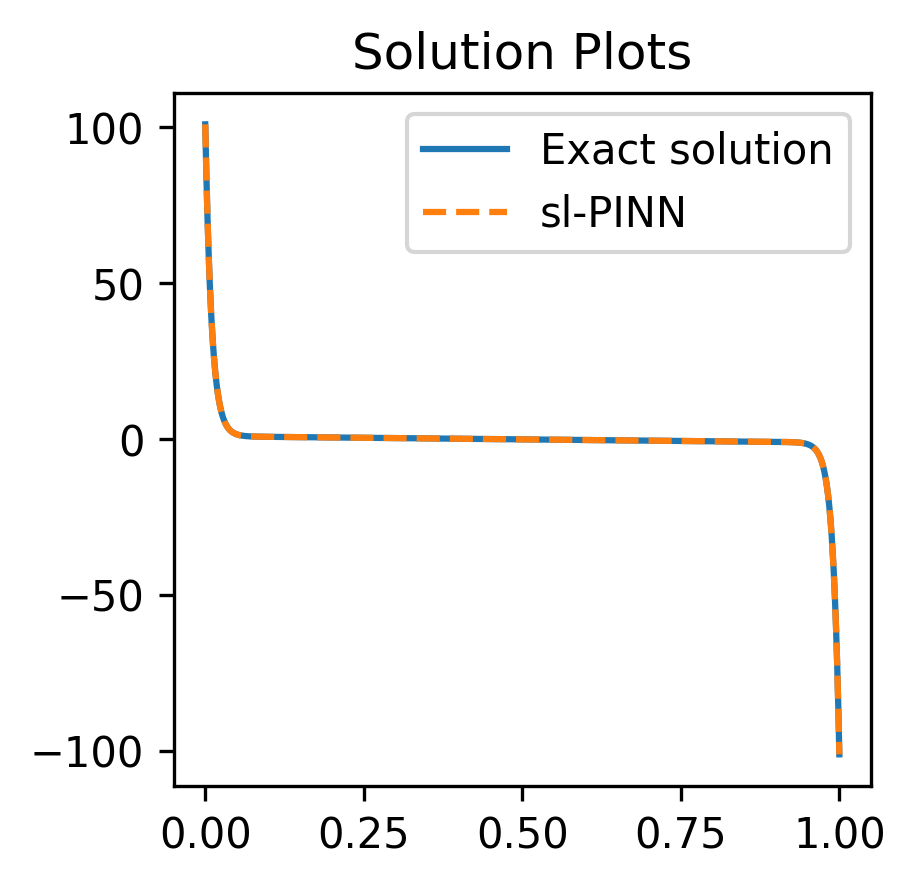}
        \includegraphics[width=\linewidth]{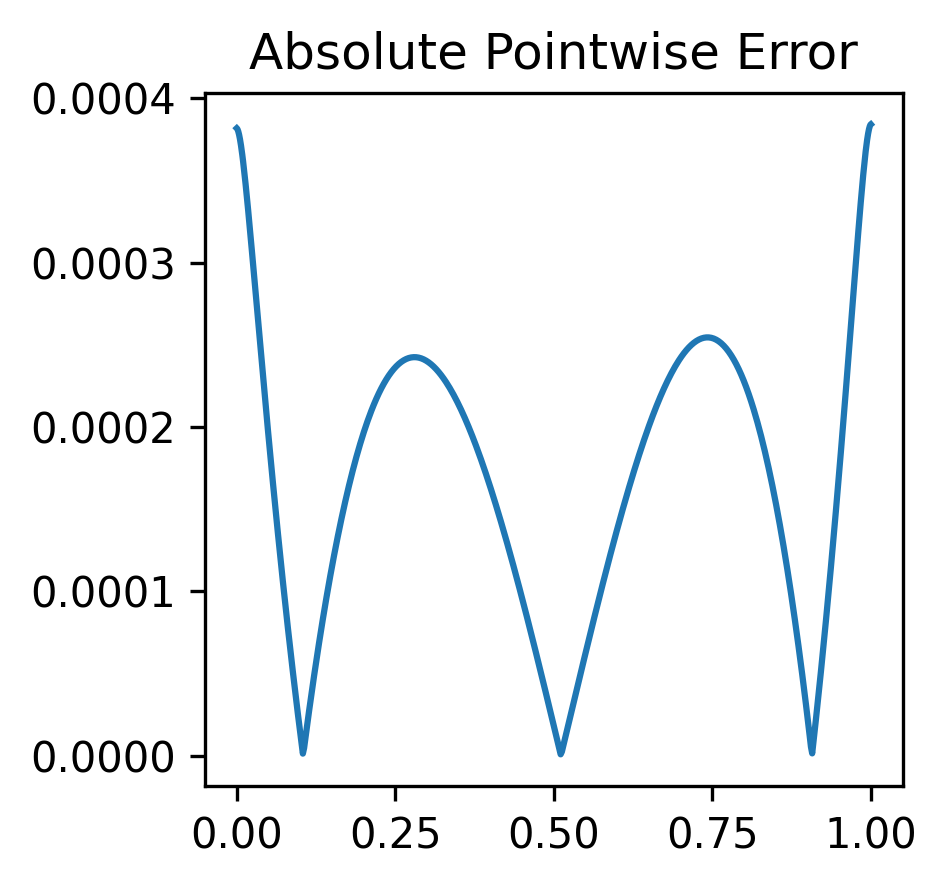}
        \caption{$\ep=10^{-4}$}
    \end{subfigure}%
    \hfill
    \begin{subfigure}[t]{0.3\textwidth}
        \raggedleft
        \includegraphics[width=\linewidth]{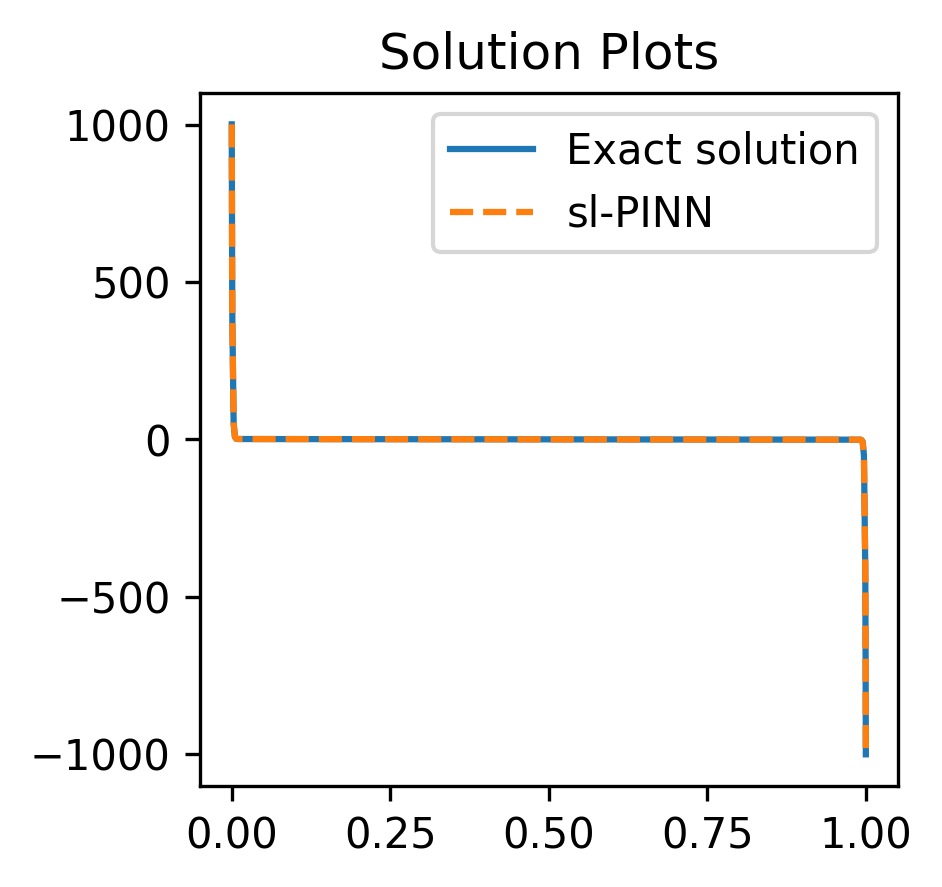}
        \includegraphics[width=0.85\linewidth]{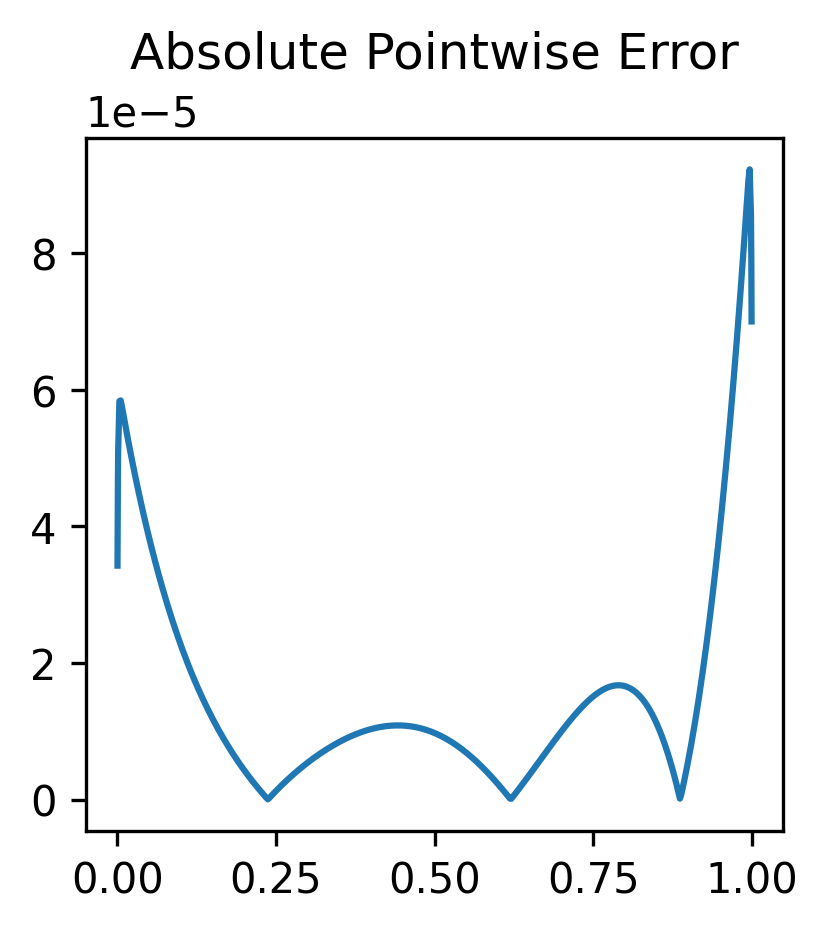}
        \caption{$\ep=10^{-6}$}
    \end{subfigure}%
    \hfill
    \begin{subfigure}[t]{0.3\textwidth}
        \raggedleft
        \includegraphics[width=\linewidth]{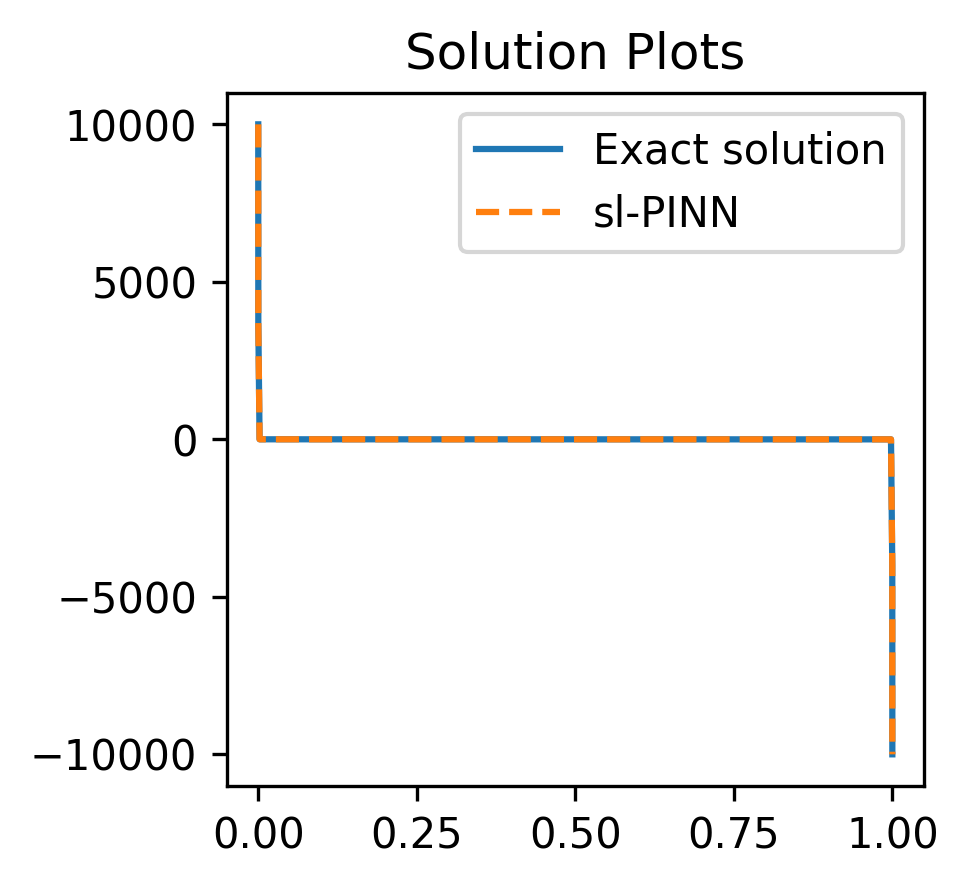}
        \includegraphics[width=\linewidth]{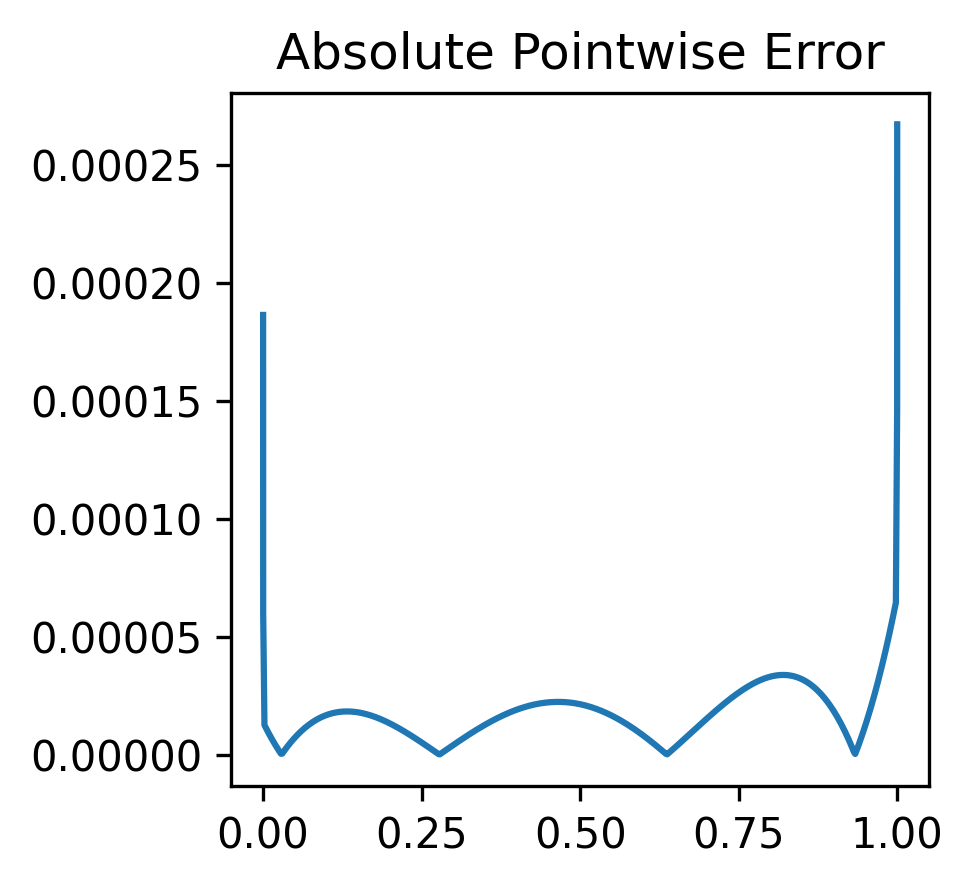}
        \caption{$\ep=10^{-8}$}
    \end{subfigure}

    \caption{Exact solutions and {\em sl-PINNs} predictins of $\omega_{2}^{\ep}$.}
    \label{plotw2EPINN}
\end{figure}

\begin{figure}
    \centering
    \begin{subfigure}{\textwidth}
        \centering
        \includegraphics[width=\linewidth]{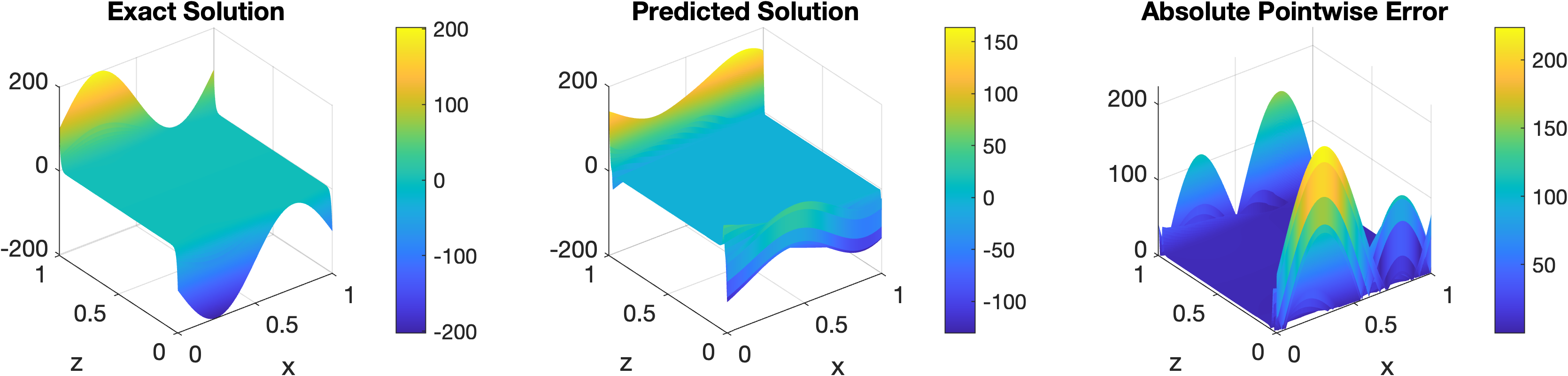}
        \caption{$\ep = 10^{-4}$}
    \end{subfigure}
    \begin{subfigure}{\textwidth}
        \centering
        \includegraphics[width=\linewidth]{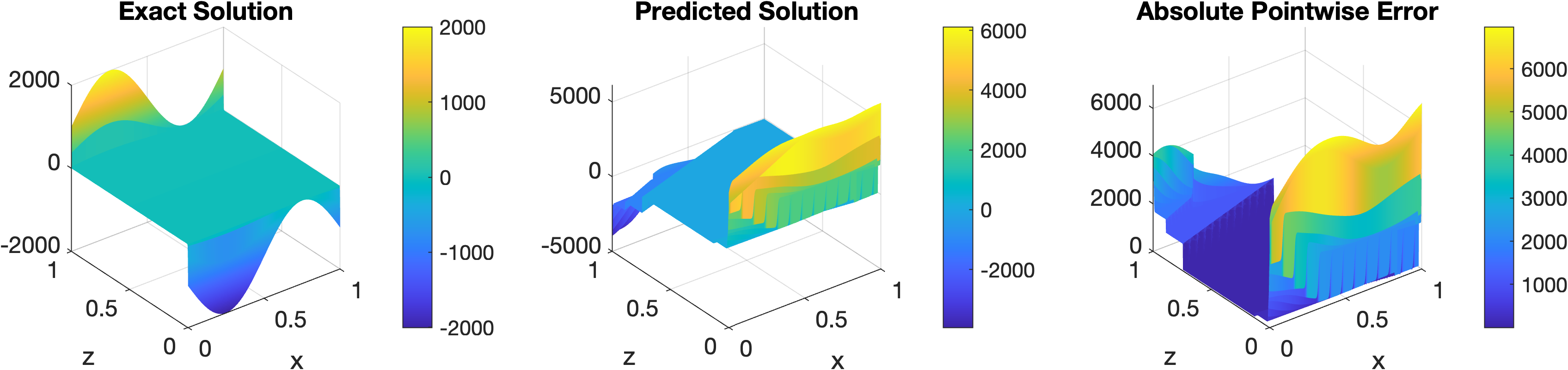}
        \caption{$\ep = 10^{-6}$}
    \end{subfigure}
    \begin{subfigure}{\textwidth}
        \centering
        \includegraphics[width=\linewidth]{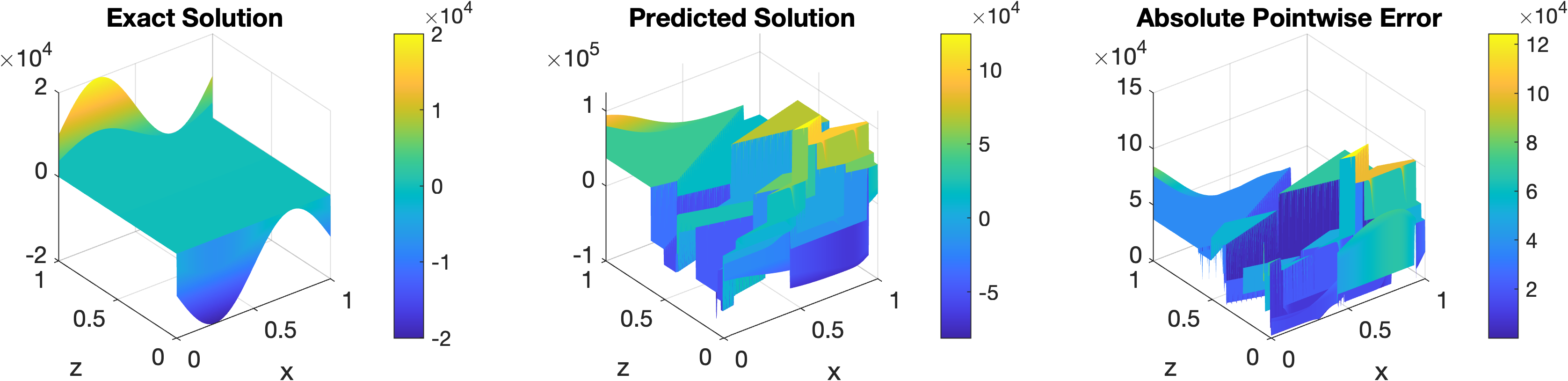}
        \caption{$\ep = 10^{-8}$}
    \end{subfigure}
    \caption{Exact solutions and PINN predictins of  $\omega_{1}^{\ep}$.}
    \label{plotw1PINN}
\end{figure}

\begin{figure}
    \centering
    \begin{subfigure}{\textwidth}
        \centering
        \includegraphics[width=\linewidth]{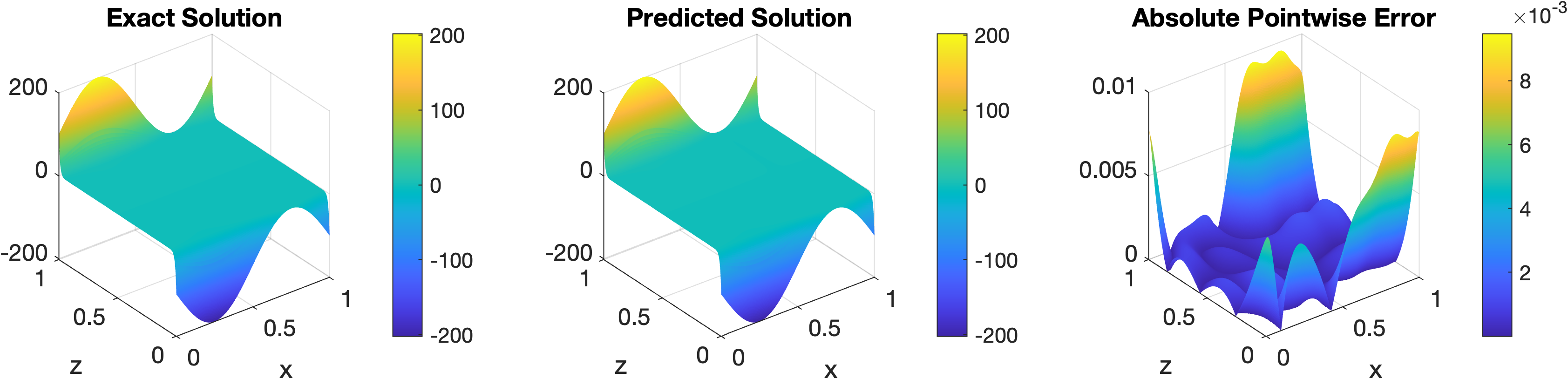}
        \caption{$\ep = 10^{-4}$}
    \end{subfigure}
    \begin{subfigure}{\textwidth}
        \centering
        \includegraphics[width=\linewidth]{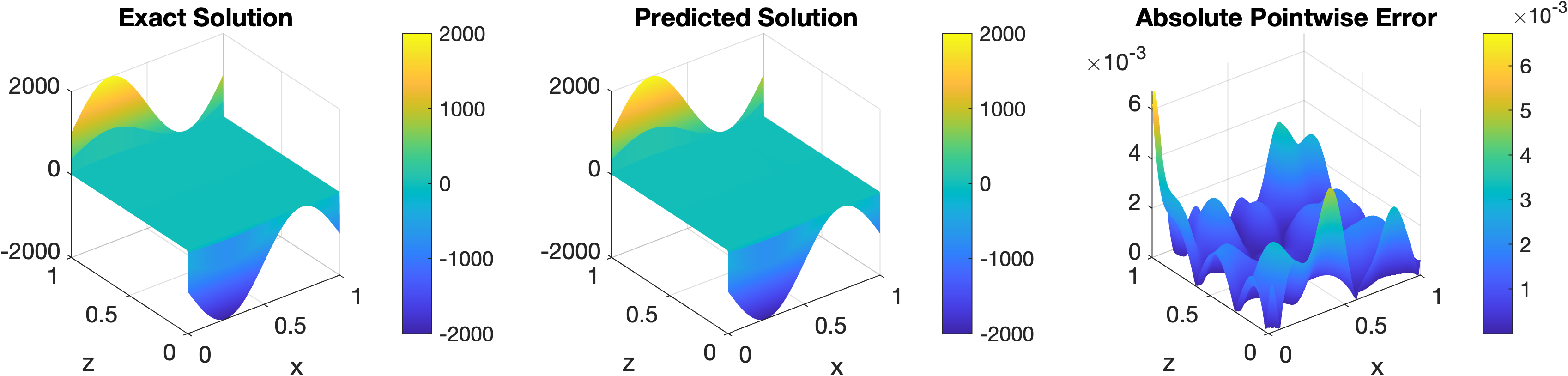}
        \caption{$\ep = 10^{-6}$}
    \end{subfigure}
    \begin{subfigure}{\textwidth}
        \centering
        \includegraphics[width=\linewidth]{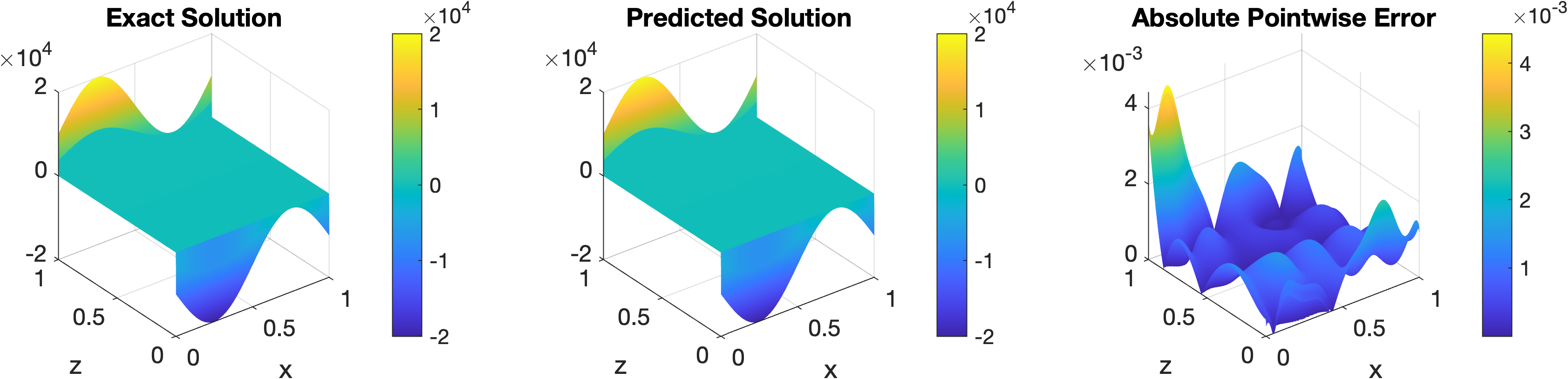}
        \caption{$\ep = 10^{-8}$}
    \end{subfigure}
    \caption{Exact solutions and {\em sl-PINNs} predictins of $\omega_{1}^{\ep}$.}
    \label{plotw1EPINN}
\end{figure}

The computational errors are shown in Tables \ref{tablew2} and \ref{tablew1} for $\omega^{\ep}_{2}$ and $\omega^{\ep}_{1}$ respectively. The {\em sl-PINNs} obtain remarkable accuracy for every small $\ep$, while the PINNs fail to  approximate the solutions. 
More detailed information about the predicted solutions appear in 
Figures \ref{plotw2PINN} and \ref{plotw2EPINN} for $\omega^{\ep}_{2}$ and Figures \ref{plotw1PINN} and \ref{plotw1EPINN} for $\omega^{\ep}_{1}$. 
We conclude from the simulations that 
our {\em sl-PINNs} produce a sharp prediction of the vorticity at every small value of the viscosity $\ep$, while 
the PINNs completely fail to predict the solution.

\begin{table}[]
\begin{tabular}{|l|rr|rr|}
\hline
                   & \multicolumn{1}{l}{PINNs}                   & \multicolumn{1}{l|}{}                            & \multicolumn{1}{l}{{\em sl-PINNs}}                  & \multicolumn{1}{l|}{}                            \\ \hline
                   & \multicolumn{1}{l|}{Relative $L^{2}$ error} & \multicolumn{1}{l|}{Relative $L^{\infty}$ error} & \multicolumn{1}{l|}{Relative $L^{2}$ error} & \multicolumn{1}{l|}{Relative $L^{\infty}$ error} \\ \hline
$\epsilon=10^{-3}$ & \multicolumn{1}{r|}{2.1347E-03}             & 2.3484E-03                                       & \multicolumn{1}{r|}{2.7971E-05}             & 1.2545E-05                                       \\
$\epsilon=10^{-4}$ & \multicolumn{1}{r|}{9.1606E-02}             & 7.7010E-02                                       & \multicolumn{1}{r|}{8.2069E-06}             & 3.8037E-06                                       \\
$\epsilon=10^{-5}$ & \multicolumn{1}{r|}{2.8411E-01}             & 6.4230E-01                                       & \multicolumn{1}{r|}{8.6951E-08}             & 5.3314E-08                                       \\
$\epsilon=10^{-6}$ & \multicolumn{1}{r|}{6.7560E-01}             & 8.8517E-01                                       & \multicolumn{1}{r|}{1.2268E-07}             & 9.2237E-08                                       \\
$\epsilon=10^{-7}$ & \multicolumn{1}{r|}{9.0956E-01}             & 1.4045E+00                                       & \multicolumn{1}{r|}{8.3588E-08}             & 5.2581E-08                                       \\
$\epsilon=10^{-8}$ & \multicolumn{1}{r|}{3.1965E+00}             & 4.6734E+00                                       & \multicolumn{1}{r|}{2.5566E-08}             & 2.6715E-08                                       \\ \hline
\end{tabular}
\caption{Comparison between PINNs and {\em sl-PINNs} for   $\omega_{2}^{\ep}$.}
\label{tablew2}
\end{table}

\begin{table}[]
\begin{tabular}{|l|rr|rr|}
\hline
                   & \multicolumn{1}{l}{PINNs}                   & \multicolumn{1}{l|}{}                            & \multicolumn{1}{l}{{\em sl-PINNs}}                  & \multicolumn{1}{l|}{}                            \\ \hline
                   & \multicolumn{1}{l|}{Relative $L^{2}$ error} & \multicolumn{1}{l|}{Relative $L^{\infty}$ error} & \multicolumn{1}{l|}{Relative $L^{2}$ error} & \multicolumn{1}{l|}{Relative $L^{\infty}$ error} \\ \hline
$\epsilon=10^{-3}$ & \multicolumn{1}{r|}{7.5808E-01}             & 8.5224E-01                                       & \multicolumn{1}{r|}{1.9381E-04}             & 1.0485E-04                                       \\
$\epsilon=10^{-4}$ & \multicolumn{1}{r|}{8.5395E-01}             & 1.1084E+00                                       & \multicolumn{1}{r|}{8.7578E-05}             & 4.6931E-05                                       \\
$\epsilon=10^{-5}$ & \multicolumn{1}{r|}{8.2919E-01}             & 7.1802E-01                                       & \multicolumn{1}{r|}{9.1033E-06}             & 5.1658E-06                                       \\
$\epsilon=10^{-6}$ & \multicolumn{1}{r|}{5.8776E+00}             & 3.4830E+00                                       & \multicolumn{1}{r|}{4.6694E-06}             & 3.3595E-06                                       \\
$\epsilon=10^{-7}$ & \multicolumn{1}{r|}{1.3017E+02}             & 4.2646E+01                                       & \multicolumn{1}{r|}{3.1714E-06}             & 2.6672E-06                                       \\
$\epsilon=10^{-8}$ & \multicolumn{1}{r|}{1.2058E+01}             & 6.2136E+00                                       & \multicolumn{1}{r|}{2.8830E-07}             & 2.2172E-07                                       \\ \hline
\end{tabular}
\caption{Comparison between PINNs and {\em sl-PINNs} for  $\omega_{1}^{\ep}$.}
\label{tablew1}
\end{table}

\newpage
\section*{Acknowledgments}

\noindent 
Teng-Yuan Chang gratefully acknowledges the financial support provided by the Graduate Students Study Abroad Program grant funded by the National Science and Technology Council (NSTC) in Taiwan.
Gie was partially supported by 
Ascending Star Fellowship, Office of EVPRI, University of Louisville; 
Simons Foundation Collaboration Grant for Mathematicians; 
Research R-II Grant, Office of EVPRI, University of Louisville; 
Brain Pool Program through the National Research Foundation of Korea (NRF) (2020H1D3A2A01110658).  
Hong was supported by Basic Science Research Program through the National Research Foundation of Korea (NRF) funded by the Ministry of Education (NRF-2021R1A2C1093579). 
Jung was supported by the National Research Foundation of Korea(NRF) grant
funded by the Korea government(MSIT) (No. 2023R1A2C1003120).
% \smallskip
% \smallskip
% \newline

% {\bf Data availability} 
% Enquiries about data availability should be directed to the authors.
% \smallskip

% {\bf Declarations}

% {\bf 
% Conflict of interest} 

%\input{figs}

\end{document}